\documentclass[12pt]{book}

\usepackage{fullpage}

\usepackage[english]{babel}

\usepackage{amsmath,amsthm}
\usepackage{amsfonts}
\usepackage{mathtools}
\usepackage{amssymb}
\usepackage{mathrsfs}
\usepackage{stmaryrd} 

\usepackage{bbm}

\usepackage{enumerate}
\usepackage{mdwlist}
\usepackage[normalem]{ulem}
\usepackage{crossreference}
\usepackage{color}

\usepackage{tikz}

\usepackage{perpage}

\newtheorem{thm}{Theorem}[section]
\newtheorem*{thm*}{Theorem}
\newtheorem{cor}[thm]{Corollary}
\newtheorem*{cor*}{Corollary}
\newtheorem{lem}[thm]{Lemma}
\newtheorem*{lem*}{Lemma}
\newtheorem{prop}[thm]{Proposition}
\newtheorem*{prop*}{Proposition}
\newtheorem*{claim}{Claim}

\newtheorem{defn}[thm]{Definition}
\newtheorem*{defn*}{Definition}
\newtheorem{notation}[thm]{Notation}
\newtheorem*{notation*}{Notation}
\newtheorem{fact}[thm]{Fact}
\newtheorem*{fact*}{Facts}

\newtheorem*{question}{Question}
\newtheorem{exe}{Exercise}
\theoremstyle{definition}
\newtheorem*{hint}{Hint}
\newtheorem{exa}[thm]{Example}
\newtheorem{rem}[thm]{Remark}
\newtheorem*{rem*}{Remark}

\MakePerPage{footnote} 
\usepackage[colorlinks=true, allcolors=blue]{hyperref}
\begin{document}

\newcommand{\cupdot}{\mathbin{\mathaccent\cdot\cup}}
\makeatletter
\def\moverlay{\mathpalette\mov@rlay}
\def\mov@rlay#1#2{\leavevmode\vtop{
   \baselineskip\z@skip \lineskiplimit-\maxdimen
   \ialign{\hfil$\m@th#1##$\hfil\cr#2\crcr}}}
\newcommand{\charfusion}[3][\mathord]{
    #1{\ifx#1\mathop\vphantom{#2}\fi
        \mathpalette\mov@rlay{#2\cr#3}
      }
    \ifx#1\mathop\expandafter\displaylimits\fi}
\makeatother
\newcommand{\bigcupdot}{\charfusion[\mathop]{\bigcup}{\cdot}}

\def \b {\begin}
\def \e {\end}

\def \l {\left}
\def \m {\middle}
\def \r {\right}
\def \mi {\;\middle|\;}

\def \mb {\mathbb}
\def \mbm {\mathbbm}
\def \mc {\mathcal}
\def \mf {\mathfrak}
\def \mr {\mathrm}
\def \ms {\mathscr}
\def \bf {\mathbf}

\def \f {\varphi}
\def \ep {\varepsilon}

\def \ss {\subseteq}
\def \sm {\setminus}
\def \fr {\frac}
\def \too {\longrightarrow}
\def \pr {\prime}
\def \ceq {\coloneqq}

\def \t {\tilde}

\def \x {\!\times\!}
\def \cd {\cdot}
\def \ci {\circ}

\def \i {\item}

\def \then {\Longrightarrow}
\def \neht {\Longleftarrow}

\def \pa {\path}
\def \dr {\draw}
\def \at {coordinate}

\def \pmtb {\begin{pmatrix}}
\def \pmte {\end{pmatrix}}
\def \bsmt {\begin{smallmatrix}}
\def \esmt {\end{smallmatrix}}

\def \bml {\begin{align*}}
\def \eml {\end{align*}}

\def \g {\gamma}
\def \G {\Gamma}
\def \S {\Sigma}
\def \s {\sigma}
\def \Om {\Omega}

\def \acts {\!\curvearrowright\!}
\def \nrml {\trianglelefteqslant}
\def \Ad {\mr {Ad}}
\def \sph {\mb S^2}
\def \sphn {\mb S^n}

\def \le {\leqslant}
\def \ge {\geqslant}

\def \lig {\ell_\infty\l(\Gamma\r)}
\def \leg {\ell_1\l(\Gamma\r)}
\def \ltg {\ell_2\l(\Gamma\r)}
\def \lpg {\ell_p\l(\Gamma\r)}
\def \lqg {\ell_q\l(\Gamma\r)}

\def \lpgg {\ell_p\l(G\r)}

\def \liz {\ell_\infty\l(\mb Z\r)}
\def \loz {\ell_1\l(\mb Z\r)}
\def \ltz {\ell_2\l(\mb Z\r)}

\def \lin {\ell_\infty\l(\mb N\r)}
\def \lon {\ell_1\l(\mb N\r)}

\def \rot {\mr{SO}\l(3\r)}
\def \son {\mr{SO}\l(n\r)}
\def \psl {\mr{PSL}_2 \l(\mb R\r)}
\def \pslz {\mr{PSL}_2 \l(\mb Z\r)}
\def \wkly {\overset{w}{\too}}
\def \sltq {\mr{SL}_2\l(\mb Q_p\r)}
\def \sltz {\mr{SL}_2\l(\mb Z_p\r)}
\def \sltk {\mr{SL}_2\l(k\r)}
\def \glnc {\mr{GL}_n\l(\mb C\r)}

\def \matI {\l(\bsmt1 &0 &0\\0 &1 &0\\0 &0 &1\esmt\r)}
\def \mat1 {\l(\bsmt1 &0 &0\\0 &1 &0\\0 &0 &{-1}\esmt\r)}

\def \L {\mr{LIM}}

\def \bd {\partial}

\title{An Invitation to Analytic Group Theory}
\date{}
\author{Tal Cohen and Tsachik Gelander}

\maketitle

\pagebreak
\setcounter{tocdepth}{2}

\chapter*{}
\begin{center}
{\it Dedicated to the memory of Nicolas Bergeron}
\end{center}

\chapter*{Preface}

This book is based on a course that I gave at the Weizmann Institute of Science in 2015. The original title of the course was `Geometric and Combinatorial Group Theory' but while teaching it I drifted towards the analytic sides of the theory. Much attention is devoted to the study of Amenability and Kazhdan's property (T) which are
perhaps the most important analytic properties of a group, but we also touch other analytic notions. I tried to introduce tricks, ideas and lemmas which repeatedly appear useful in various situations. My main guideline was to expose the beauty of the theory and touch many different aspects while keeping the text short, simple and accessible, sometimes on the expense of diving deep or providing thorough expositions.

The formation of the book relied on lecture notes taken by
Tal Cohen with whom, a few years later, we improved and expanded the text.
Tal's enthusiasm and persistence were crucial for the process, and Tal's writing style is seeded in the text.
Hopefully this book could serve as a smooth entry to Analytic Group Theory.

\vspace{2cm}

Tsachik Gelander

\tableofcontents

\let\stdpart\part
\renewcommand\part{\newpage\stdpart}
\cleardoublepage
\phantomsection
\addcontentsline{toc}{chapter}{~~~~~Introduction}
\chapter*{Introduction}

This book is concerned with analytic approaches to study groups and their actions. 
It could serve as a graduate course textbook or as a lightweight book to take for a one month trek in the Nepali Himalayas.

\medskip

The first chapter is dedicated to the famous Hausdorff–Banach–Tarski paradox, which says that one can decompose the sphere into finitely many pieces, and, by rotating these pieces, build two new spheres of the same size. This is done via ‘paradoxical decompositions’ of the free group, and sets the stage for one of the most important properties of analytic group theory: amenability. Amenable groups are, intuitively speaking, exactly those groups that do not allow such ‘paradoxical’ behaviour. They are a generalisation of abelian groups as well as finite groups. In Chapter 2, we explore the basic features of amenability and give a few equivalent definitions.

Chapter 3 is about actions of amenable groups. We discuss two additional equivalent definitions of amenability, defined via actions on compact Hausdorff spaces and on ‘compact-convex subsets’ of topological vector spaces.

Chapter 4 is dedicated to two famous applications of amenability: Banach limits and the von Neumann ergodic theorem. Chapter 5 introduces the notion of the growth rate of finitely generated groups, and discusses its relation to amenability.

In Chapter 6, we introduce the notion of topological groups, generalising the notion of abstract (or discrete) groups. 
We present the two most important classes: Lie groups and totally disconnected locally compact groups, which, (very) roughly speaking, are used to study symmetries of connected spaces and combinatorial objects respectively. We define the Haar measure of locally compact groups, and generalise the definition of amenability to topological groups. 

Chapter 7 is dedicated to the study of another important analytic property of groups: property (T). While this notion is somewhat opposite to amenability,
it takes longer to develop an intuition for it. Amenable groups are in some sense small, but groups with property (T) are not exactly big, but rather more rigid.

Chapter 8 is about the so-called classical groups, focusing on $\mr{SL}_n(\mb R)$. The main result is to establish a metric version of the Moore ergodicity theorem, and to deduce that the group $\mr{SL}_n(\mb R)$ has property (T) when $n\ge 3$. We then define property FH in chapter 9, and show that it is equivalent (for most interesting groups) to property (T).

Our stage is then set to demonstrate the power of the tools we established, 
and prove various classical theorems.
In Chapter 10 we prove Zassenhaus’ theorem and the Margulis lemma, and we then deduce three famous theorems: the Jordan theorem, the Bieberbach theorem and a theorem of Zimmer (a baby version of the Zimmer program). Along the way we develop further important tools, such as the $p$-adic numbers.
Chapter 11 is dedicated to the celebrated Tits theorem, which says that every finitely generated subgroup of $\mr{GL}_n(\mb R)$ either contains a free subgroup or is virtually solvable.

\chapter{The Hausdorff--Banach--Tarski Paradox}
Consider the following question:
\begin{question}
	Does there exist a finitely additive probability measure defined on
	all subsets of $\mathbb{S}^{n}$ that is invariant under rotations?
	In other words, a map $m:\mathcal{P}(\mb S^n)\to\l[0,1\r]$ that is
	\begin{enumerate}
		\item \emph{Finitely additive}: if $A,B\subseteq\mathbb{S}^{n}$
		and $A\cap B=\varnothing$, then $m(A\cup B)=m\l(A\r)+m\l(B\r)$.
		\item A \emph{probability} measure: $m(\mathbb{S}^{n})=1$.
		\item \emph{Invariant under rotations}:  $m\l(A\r)=m\l(gA\r)$
		for every rotation $g$ and every $A\subseteq\mathbb{S}^{n}$.
	\end{enumerate}
\end{question}

Formally, a \textit{rotation} is an element in the \textit{Special Orthogonal Group} $\mr{SO}(n+1,\mb R)$.

\begin{thm*}[Banach]
	There exists such a measure on $\mathbb{S}^{1}$. 
\end{thm*}

\begin{thm*}[Hausdorff]
	There doesn't exist such a measure on $\mathbb{S}^{n}$ for $n\geqslant2$. 
\end{thm*}
In this chapter we shall prove the second theorem.
\begin{rem*}
	Observe that, usually, the word measure is more restrictive: it is
	required to be \emph{$\sigma$-additive}. That is, to satisfy the
	following: if $\left\{ A_{i}\right\} _{i\in\mathbb{N}}$ is a countable
	collection of measurable subsets that are pairwise disjoint ($A_{i}\cap A_{j}\neq\varnothing$
	for $i\neq j$), then  $m\left(\bigcupdot_{i=1}^{\infty}A_{i}\right)=\sum_{i=1}^{\infty}m(A_{i})$.
	It is easier to see that there is no $\sigma$-additive probability
	measure that is defined on all subsets of $\mathbb{S}^{n}$ and that
	is invariant under rotations.
\end{rem*}

\section{Equidecomposability}

From now on, whenever we write $G\acts X$, it means
that $G$ is a group that acts on some set $X$ (or, later on, a set
$X$ endowed with some structure, like a topological space). In this
case we also say that $X$ is a \emph{$G$-set }(or a \emph{$G$-space}).
\begin{defn}
	If $G\acts X$, two subsets $A,B\subseteq X$ are called \textbf{\emph{$G$-equidecomposable}}
	if they admit a decomposition into disjoint subsets,
	\[
	A=\bigcupdot_{i=1}^{n}A_{i},B=\bigcupdot_{i=1}^{n}B_{i},
	\]
	and there are elements $g_{1},\dots,g_{n}\in G$ such that $g_{i}A_{i}=B_{i}$
	for $i=1,\dots,n$. This is denoted by $A\sim B$. 
\end{defn}

\begin{defn}
	If $G\acts X$, a subset $A\subseteq X$ is \textbf{\emph{$G$-embeddable}}
	in a subset $B\subseteq X$ if it admits a decomposition into disjoint
	subsets $A=\bigcupdot_{i=1}^{n}A_{i}$ and there are elements $g_{1},\dots,g_{n}\in G$
	such that $\bigcupdot g_{i}A_{i}\subseteq B$ (a disjoint union). This is denoted by
	$A\prec B$. 
\end{defn}

\begin{rem}
	It is immediate that $\sim$ is an equivalence relation, and that
	$A\sim B$ implies that both $A\prec B$ and $B\prec A$. It is also immediate that $A\prec B$ if and only if $A$ is $G$-equidecomposable with some subset of $B$.
\end{rem}

We may now state the theorem we want to prove: 
\begin{thm*}[Hausdorff]
	For $n\geqslant2$, there are disjoint subsets $A,B\subseteq\mathbb{S}^{n}$
	such that $\mathbb{S}^{n}\prec A$, $\mathbb{S}^{n}\prec B$ with
	respect to the action of $\mathrm{SO}\l(n+1\r)$. 
\end{thm*}

The theorem we stated before follows directly from this one;  if there were such a measure, two such sets could not exist (as we would get $1=m\l(\mb S^n\r)\ge m\l(A\r)+m\l(B\r)\ge 2$).

Before we get to the proof, we make a few more statements and definitions:

\begin{thm}[Banach--Tarski]
\label{BT}
Let $G$ be a group acting on a set $X$, and let $A,B\ss X$ be subsets. If $A\prec B$ and $B\prec A$, then $A\sim B$.
\end{thm}
\begin{proof}
By definition, there is a decomposition $\l(A_i\r)_{i=1}^n$ of $A$ and elements $\l(f_i\r)_{i=1}^n$ in $G$ such that $\bigcupdot f_i A_i\subseteq B$, and a decomposition $\l(B_j\r)_{j=1}^m$ of $B$ and elements $\l(g_j\r)_{j=1}^m$ of $G$ such that $\bigcupdot g_jB_j\subseteq A$. 

Define a map $f:A\to B$ by setting $f(x)= f_i.x$ for $x\in A_i$, and a map $g:B\to  A$ by setting $f(y)= g_j.y$ for $y\in B_j$. Both maps $f:A\to B, g:B\to A$ are injective, but not necessarily onto. 

We can now define for each $x\in A$ a sequence of points, alternating between $A$ and $B$, using $f$ and $g$; the first point is $x\in A$, the second is $f\l(x\r)\in B$, the third is $g\l(f\l(x\r)\r)\in A$, and so forth. This sequence may also go backwards, if $x$ is in the image of $g$. The point before $x$ would be $g^{-1}\l(x\r)$, and then, if it exists, $f^{-1}\l(g^{-1}\l(x\r)\r)$, and so forth. The sequence looks like this: $$\cdots \rightarrow  f^{-1}(g^{-1}(a)) \rightarrow g^{-1}(a) \rightarrow   a  \rightarrow  f(a) \rightarrow  g(f(a)) \rightarrow \cdots $$
We have four options:

\paragraph{}

\begin{tikzpicture}


\draw (0,3) node[above] {\uline{Option 1}};

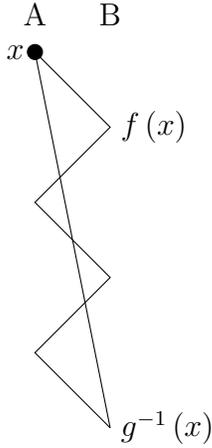
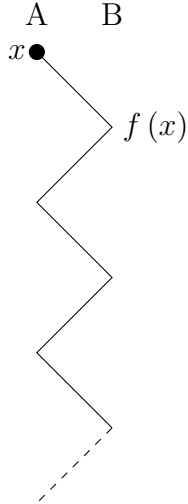
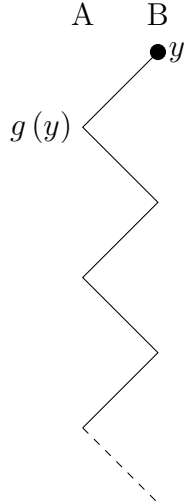
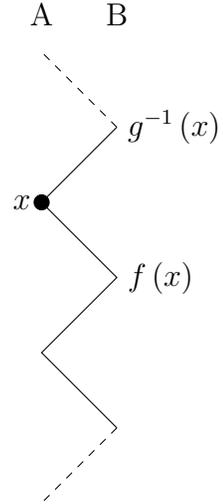
\draw (0,2.5) node {A};
\draw (1,2.5) node {B};

\path (1,-3) coordinate (A1);
\path (0,-2) coordinate (A2);
\path (1,-1) coordinate (A3);
\path (0, 0) coordinate (A4);
\path (1,1) coordinate (A5);
\path (0,2) coordinate (A6);
\draw (A1)node[right]{$g^{-1}\l(x\r)$} -- (A2) -- (A3) -- (A4) -- (A5) node[right] {$f\l(x\r)$} -- (A6) node[left] {$x$};
\draw (A6) -- (A1);

\draw [fill, color=black] (0,2) circle (0.1cm);


\draw (4,3) node[above] {\uline{Option 2}};

\draw (4,2.5) node {A};
\draw (5,2.5) node {B};

\path (5,-3) coordinate (B1);
\path (4,-2) coordinate (B2);
\path (5,-1) coordinate (B3);
\path (4, 0) coordinate (B4);
\path (5,1) coordinate (B5);
\path (4,2) coordinate (B6);
\path (4,-4) coordinate (B7);
\draw (B1) -- (B2) -- (B3) -- (B4) -- (B5)node[right] {$f\l(x\r)$} -- (B6) node[left] {$x$};

\draw[dashed] (B1) -- (B7);

\draw [fill, color=black] (B6) circle (0.1cm);


\draw (8,3) node[above] {\uline{Option 3}};

\draw (8,2.5) node {A};
\draw (9,2.5) node {B};

\path (8,-3) coordinate (C1);
\path (9,-2) coordinate (C2);
\path (8,-1) coordinate (C3);
\path (9, 0) coordinate (C4);
\path (8,1) coordinate (C5);
\path (9,2) coordinate (C6);
\path (9,-4) coordinate (C7);
\draw (C1) -- (C2) -- (C3) -- (C4)  -- (C5) node[left]{$g\l(y\r)$} -- (C6)node[right] {$y$};

\draw[dashed] (C1) -- (C7);

\draw [fill, color=black] (C6) circle (0.1cm);


\draw (12,3) node[above] {\uline{Option 4}};

\draw (12,2.5) node {A};
\draw (13,2.5) node {B};

\path (13,-3) coordinate (D1);
\path (12,-2) coordinate (D2);
\path (13,-1) coordinate (D3);
\path (12, 0) coordinate (D4);
\path (13,1) coordinate (D5);
\path (12,2) coordinate (D6);
\path (12,-4) coordinate (D7);
\draw (D1) -- (D2) -- (D3)node[right] {$f\l(x\r)$} -- (D4) node[left] {$x$} -- (D5)node[right]{$g^{-1}\l(x\r)$};

\draw[dashed] (D1) -- (D7);
\draw[dashed] (D5) -- (D6);
\draw [fill, color=black] (12,0) circle (0.1cm);

\end{tikzpicture}
\begin{description}
\item[Option 1:] The sequence is cyclic.
\item[Option 2:] The sequence has a starting point in $A$, and is infinite.
\item[Option 3:] The sequence has a starting point in $B$, and is infinite.
\item[Option 4:] The sequence is infinite in both directions.
\end{description}

Each $x \in A$ induces a unique sequence, and this sequence satisfies exactly one of these options. We can define a bijection $h:A\leftrightarrow B$: if $x$ is in a sequence of options $1$, $2$ or $4$ then $h\l(x\r)=f\l(x\r)$; otherwise, $h\l(x\r)=g^{-1}\l(x\r)$. You can check this is indeed a bijection.

Now we can decompose $A$ into $\l(\Lambda_i\r)_{i=1}^{n+m}$ as follows: $\Lambda _i$ for $i=1,\dots,n$ is the set of all $x$ such that $h\l(x\r)=f_ix$, and $\Lambda_{n+j}$ for $j=1,\dots,m$ is the set of all $x$ such that $h\l(x\r)=g_j^{-1}x$. Hence by choosing the elements $f_1,\dots,f_n,g_1^{-1},\dots,g_m^{-1}$ in $G$, we get $A\sim B$.
\end{proof}
\begin{rem}
	Observe that this is exactly the proof of the Cantor--Bernstein theorem in set theory. We just needed to observe that the map obtained by this proof is still ``piecewise elemental". To be more precise:
	
	We say that a map $f:A\to B$, between two subsets $A,B\subseteq X$
	of a $G$-set $X$, is \emph{piecewise elemental }if $A$ can be written
	as a finite union, $A=\bigcup_{i=1}^{n}A_{i}$, such that the restriction
	of $f$ to each $A_{i}$ is given by the action of a single element
	of $G$; that is, if there are $g_{1},\dots,g_{n}\in G$ such that
	$f(a)=g_{i}.a$ for every $a\in A_{i}$ (for every $i=1,\dots,n$).\emph{
	}	Then $A$ is $G$-embeddable in $B$ if and only if there is an injective,
	piecewise elemental map $f:A\to B$. Similarly, $A$ and $B$ are
	$G$-equidecomposable if and only if there is a bijective, piecewise
	elemental map $f:A\to B$.

\end{rem}
\begin{exa}
Consider $X=\mb{R}^2$, $G$ the group of similarities on $X$ (isometries and scalings). Then the unit square and unit disk are $G$-equidecomposable. This is easily proved using the preceding theorem: the unit square can be shrunk to fit inside the unit disk, and vice versa.
\end{exa}

\begin{exe}
Draw such an equivalence.
\end{exe}
\begin{exe}
Find a bijection $f:\mb{N}\leftrightarrow\mb{N}$ such that, for every $a\in\mb{N}$, $f\l(a\r)=b$ implies $a=b^2$ or $a^2=b$.
\end{exe}
The following notion is useful:
\begin{defn}
Let $G\acts X$, and let $A,B\subseteq X$. We denote $mA\sim nB$ if $m$ copies of $A$ can construct $n$ copies of $B$, and $mA\prec nB$ if $m$ copies of $A$ can be embedded into $n$ copies of $B$.

More rigorously: $mA\sim nB$ if there are $m$ decompositions of $A$ into disjoint subsets, $\{\l(A_i^k\r)_{i=1}^{a_k}\}_{k=1}^m$, $n$ decompositions of $B$ into disjoint subsets, $\{\l(B_j^\ell\r)_{j=1}^{b_\ell}\}_{\ell=1}^n$, such that $a_1+a_2+\cdots+a_m=b_1+\cdots+b_n\eqqcolon r$, and such that there exist $r$ elements $\l({}_qg\r)_{q=1}^r$ in $G$ that satisfy the following: ${}_qg_.{}_qA={}_qB$ for each $q=1,\dots,r$, while ${}_qA$ and ${}_qB$ are defined to be the $q^{\mathrm{th}}$ element of (a permutation of) the sequences $A_1^1,\dots,A_{a_1}^1,\dots,A^m_1,\dots,A_{a_m}^m$ and $B_1^1,\dots,B_{a_n}^n$ respectively. Embeddability is similar.
\end{defn}
\begin{exe}
\label{can}
Prove that, if $nA\sim nB$, then $A\sim B$.
\end{exe}
\begin{hint}
        This follows from \textit{the Marriage theorem}:
If $G\l(A,B\r)$ is an $n$-regular bipartite graph with sides in $A$\ and $B$, then there exist $n$ perfect matchings.
\end{hint}

\section{Paradoxical Sets and Groups}
\begin{defn}
Given $G\acts X$, the set $X$ is said to be \emph{\textbf{$G$-paradoxical}} if $X\sim 2X$ with respect to this action.
\end{defn}
And now what we want to prove is the following:
\begin{thm}[Hausdorff]
\label{Hausdorff}
For $n\ge2$, $\mb S^n$ is paradoxical with respect to the natural action of $\mr{SO}\l(n+1\r)$.
\end{thm}

This is equivalent to the previous version of the theorem; if the first version holds, then $2\mb S^n\prec \mb S^n$ and hence $2\mb S^n\sim \mb S^n$, and if this version holds, then by definition of $2\mb S^n\sim \mb S^n$ it's easy to see the required $A,B$ exist.

We will only prove this for $n=2$ (the proof of the general case is very similar). In order to do this, we will use the free group on $2$ generators.
\begin{exa}[Important]\label{free-gp}
The free group $F_2$ on two generators is paradoxical (with respect to its action on itself by left translation; i.e., $g.x=gx$).
\end{exa}
\begin{proof}
Clearly, $F_2\prec 2F_2$, so we only need to show that $2F_2\prec F_2$.

Denote the generators of $F_2$ by $a,b$ (i.e., the elements in $F_2$ are words in $a,b$), and define  $F_x$ for $x\in \l\{a,a^{-1},b,b^{-1}\r\}$ to be the set of all words in $F_2$ starting with $x$. Clearly, we have $$F=F_a\cupdot F_b\cupdot F_{a^{-1}}\cupdot F_{b^{-1}}\cupdot \l\{1\r\}.$$ Direct computation shows that $a^{-1}F_a=F_2\setminus F_{a^{-1}}$ and $b^{-1}F_b=F_2\setminus F_{b^{-1}}$, so \begin{align*}F_2&=F_{a^{-1}}\cupdot a^{-1}F_a,\\F_2&=F_{b^{-1}}\cupdot b^{-1}F_b,\end{align*}and if we take $A_1=F_{a^{-1}}$, $g_1=1$, $A_2=a^{-1}F_a$, $g_2=a$, $A_3=F_{b^{-1}}$, $g_3=1$, $A_4=b^{-1}F_b$, $g_4=b$, we get $2F_2\prec F_2$, as desired.
\end{proof}
It is as easy as \ref{free-gp} to show $F_n$ is paradoxical for every $n\ge 2$.

For the next and main proposition, we need to make yet another definition:
\begin{defn}
We say that an action of a group $G$ on a set $X$ is \emph{\textbf{free}}, or that $G$ \emph{\textbf{acts freely}} on $X$, if $g.x\neq x$ for every $g\in G\setminus\{1\}$  and $x\in X$.
\end{defn}
\begin{rem}
If $G\acts X$ and $x\in X$, the \textit{\textbf{stabiliser}} of $x$ (with respect to this action) is $\mr{Stab}_G(x)=\{g\in G|g.x=x\}$. This is a subgroup of $G$ for every $x$. An action is free if all its stabilisers are trivial (that is, $\mr{Stab}(x)=\{1\}$ for every $x\in X$).
\end{rem}
And now we can state
\begin{prop}
\label{a}
If $F_2$ acts freely on a space $X$, then $X$ is paradoxical (with respect to this action).
\end{prop}
\begin{proof}
Let $A\subseteq X$ be a set of representatives of the $F_2$-orbits in $X$. Observe that $X=\bigcupdot _{\gamma\in F_2}\gamma A=F_2A$.
(The union is disjoint because the action is free.) 

Define a decomposition of $X$ as follows: $$X=F_aA\cupdot F_{a^{-1}}A\cupdot F_bA\cupdot F_{b^{-1}}A\cupdot A,$$ where $F_a,F_b,F_{a^{-1}},F_{b^{-1}}$ are defined as in the proof of Example~\ref{free-gp}. Now we can decompose $X$ twice, in the exact same manner as before: \begin{align*}X&=F_{a^{-1}}A\cupdot a^{-1}F_aA,\\X&=F_{b^{-1}}A\cupdot b^{-1}F_bA,\end{align*} and hence, again choosing $g_1=1$, $g_2=a$, $g_3=1$, $g_4=b$, we get $2X\prec X$, and since it is clear that $X\prec 2X$, it follows that $X\sim 2X$ and $X$ is paradoxical.
\end{proof}
\begin{rem}
\label{paradoxical}
This proposition could easily be generalised by replacing $F_2$ with any other paradoxical group.
\end{rem}

\section{The Proof}
Now we can get back to $\mb S^2$ and $\mr{SO}\l(3\r)$. 
\begin{exe}
\label{f2:so3}
Let $g_{x,\frac{1}{3}},g_{y,\frac{1}{3}}\in \mr{SO}\l(3\r)$ be the rotations by $\arctan \frac{1}{3}$ around the $x$- and $y$-axis respectively. Show that $\langle g_{x,\frac{1}{3}},g_{y,\frac{1}{3}}\rangle$ is isomorphic to the free group on two generators.
\end{exe}

It seems our proof that $\mb S^2$ is paradoxical is done; $F_2$ acts on it as a subgroup of $\mr{SO}\l(3\r)$, and hence by Proposition~\ref{a} $\mb S^2$ is paradoxical. Of course, this isn't true, because we haven't proved $F_2$ acts \emph{freely} on $\mb S^2$, and we haven't done so because it isn't true. 

What we will in fact do is to find a subset $\hat S$ of $\mb S^2$ such that $F_2$ act freely on $\hat S$ and such that $\mb S^2\sim \hat S$, and then it will follow $\mb S^2$ is paradoxical (because $\sim$ is an equivalence relation).
\begin{exe}
Show that any element in $\mr{SO}\l(3\r)$ has an axis around which it rotates.
\begin{hint}
If $n$ is odd, every $A\in M_n\l(R\r)$ has at least one real eigenvalue, and if $A$ is orthogonal, this eigenvalue can either be $1$ or $-1$. Recall also that all matrices in $\mr{SO}\l(n\r)$ have determinant $1$.
\end{hint}
\end{exe}
Using this exercise, we can prove
\begin{lem}
\label{thm:b}
If $D\ss \mb S^2$ is countable, then $\mb S^2\sim \mb S^2\setminus D$ with respect to the action of $\mr{SO}\l(3\r)$.
\end{lem}
For the proof, we take two antipodal points in $\mb S^2\setminus D$ and call them $S$ (South) and $N$ (North)\footnote{We can do so because there are uncountably many such pairs, and $D$ is only countable.}. Now we need a rotation $g$ around the North-South axis $N-S$ that satisfies $g^nD\cap D=\varnothing$ for all $n\in \mb{N}$.
\begin{claim}
\label{c}
There is such a rotation.
\end{claim}
\begin{proof}[Proof of the claim]
If $x,y\ne N,S$, there is at most one rotation around $N-S$ that sends $x$ to $y$. Therefore the set $$\l\{R\in \mr{SO}\l(3\r) : R(N)=N,\exists x\in D\text{ such that } R\l(x\r)\in D\r\}$$ (of all the rotations around 
$N-S$ that send some point in $D$ into some point in $D$) is countable, and since there are uncountably many rotations around $N-S$, there must be a rotation $g$ for which $g^nD\cap D=\varnothing$ for every $n\in\mb N$.
\end{proof}
\begin{proof}[Proof of Lemma~\ref{thm:b}]
By the claim, there is a rotation $g$ such that $g^nD\cap D=\varnothing$ for every $n\in \mb N$. Set $A=\bigcup_{n=0}^\infty g^n D$ for such a $g$. Then $$\mb S^2=A\cupdot \l(\mb S^2\setminus A\r)\sim gA\cupdot \l(\mb S^2\setminus A\r).$$Now, we claim that $$\mb S^2\setminus D=gA\cupdot \l(\mb S^2\setminus A\r).$$ 
This is actually just a simple calculation:
$$gA\cup \mb (S^2\setminus A)=\l(g\bigcup_{n=0}^\infty g^n D\r)\cup \l(\mb S^2\setminus \bigcup_{n=0}^\infty g^n D\r) =\l(\bigcup_{n=1}^\infty g^n D\r)\cup \l(\mb S^2\setminus \bigcup_{n=0}^\infty g^n D\r)=\mb S^2\setminus g^0 D.$$
Since $g^0$ is by definition the identity element, $g^0 D$ is just $D$, so we are done.
\end{proof}
Now we may finally prove our main theorem.
\begin{proof}[Proof of Theorem~\ref{Hausdorff} (Hausdorff's theorem)]\
Let $\Gamma\le \mr{SO}\l(3\r)$ be a subgroup isomorphic to $F_2$ (see Exercise~\ref{f2:so3}). Let $D$ be the set of fixed points of nontrivial elements in $\Gamma$; then $D$ is countable, since each element in $\Gamma$ has precisely two fixed points (being a rotation).
\begin{exe}
Prove $\Gamma.D=D$, and deduce $\Gamma.\l(\mb S^2\setminus D\r)=\mb S^2\setminus D$.
\end{exe}
Hence, $\Gamma$ acts on $\mb S^2\setminus D$, and it acts freely because all the fixed points of its elements are in $D$. So, $\mb S^2\setminus D$ is paradoxical. By our last theorem, $\mb S^2\sim \mb S^2\setminus D$. Thus $$\sph\sim\sph\!\sm\!D\sim2\!\cdot\!\l(\sph\!\sm\!D\r)\sim 2\!\cdot\!\sph,$$ so $\mb S^2$ is also paradoxical.
\end{proof}

An interesting result of the Hausdorff–Banach–Tarski paradox is the following:
\begin{thm}
If $A,B\subseteq \mb S^2$ have nonempty interiors, then $A\sim B$.
\end{thm}
\begin{proof}
It's enough to show that $A\sim \mb S^2$ whenever $A$ has nonempty interior.

$\mathrm{Int}A\ne\varnothing$, so by rotating some open set inside $A$ you can cover $\mb S^2$. But $\mb S^2$ is compact, so this cover has some finite subcover, and hence $\mb S^2\prec k\cdot A$ for some $k$. Moreover, obviously $k\cdot A\prec k\cdot \mb S^2$, and (by the Hausdorff--Banach--Tarski paradox) $k\cdot \mb S^2\prec \mb S^2$, so $$k\cdot \mb S^2\prec \mb S^2\prec k\cdot A \prec k\cdot \mb S^2.$$ Therefore, by Theorem~\ref{BT}, we have $k\cdot A\sim k\cdot \mb S^2$, and hence $A\sim \mb S^2$ (by Exercise~\ref{can}).
\end{proof}
\begin{exe}
	Prove the Hausdorff--Banach--Tarski paradox for $\mb S^n$ for $n\geqslant 3$.
\end{exe}
\begin{hint}
    One can use the fact that $\mr{SU}(2)$ is diffeomorphic to $\mb S^3$ (in fact, when considered as a subset of $\mathbb{C}^2\cong \mathbb{R}^4$, is \textit{is} $\mb S^3$). Its action on itself (by left or right multiplication) is by isometries, which makes it a subgroup of $\mr{SO}(4)$ acting freely on $\mb S^3$; there, it is enough to find a free subgroup of $\mr{SU}(2)$ to complete the proof. Alternatively, one may use induction and the theorem we have already proved.
\end{hint}
\begin{rem} 
The sphere $\mb S^n$ can easily be replaced here with the ball $\mb B^{n+1}$; the only problem is the point $0$, which is fixed by $\rot$. However, if you add translations to the mix, this is easily fixed. That is, considering the action of $\mr{Isom}(\mb R^n)$ on $\mb R^n$, one may show that $\mb B^n$ is a paradoxical subset.
\end{rem}
\begin{exe}
	Prove that $\mb B^n$ is a paradoxical subset of $\mb R^n$ with respect to the action of $\mr{Isom}(\mb R^n)$.
\end{exe}
Thus, we have shown that one can decompose the $n$-sphere (or $(n+1)$-ball) to finitely many (necessarily non-measurable) subsets, and build from them -- using rotations (and translations) -- two spheres (or balls).
\begin{exe}
	Follow the proof of the Hausdorff--Banach--Tarski paradox, and compute how many subsets were needed for this decomposition. 
\end{exe}
\subsection{Summary of the Proof}
For $G\acts X$ and $A,B\ss X$, we denoted $A\sim B$ if there is a decomposition $\l(A_i\r)_{i=1}^n$ of $A$ into disjoint sets, and elements $\l(g_i\r)_{i=1}^n$ in $G$, such that $\bigcupdot g_iA_i=B$. 
We called $X$ paradoxical if $X\sim{2X}$, which means there is a decomposition $\l(X_i\r)_{i=1}^{n+m}$ of $X$ and elements $\l(g_i\r)_{i=1}^{n+m}$ such that $\bigcupdot_{i=1}^n g_iX_i=\bigcupdot_{j=1}^{m} g_{n+j} X_{n+j}=X$. Then we proved
\begin{thm*}
For $n\ge 2$, $\mb S^n$ is paradoxical (with respect to the action of $\mr{SO}\l(n+1\r)$).
\end{thm*}

The outline of the proof was:
\begin{enumerate}
\i $F_2$ is paradoxical (with respect to its action on itself by left translations).
\i If $F_2$ acts freely on some set $X$, then $X$ is paradoxical (with respect to this action).
\i $\mr{SO}\l(3\r)$ admits a subgroup isomorphic to $F_2$, but this subgroup doesn't act freely on $\mb S^2$.
\i $\mb S^2\sim \mb S^2\setminus D$ for every countable $D\ss \mb S^2$.

\i $F_2$ acts freely on $\sph\sm D$, where $D$ is the (countable) set of fixed points of (nontrivial) elements of $F_2$, and hence $\sph\sim\sph\!\sm\!D\sim2\!\cdot\!\l(\sph\!\sm\!D\r)\sim 2\!\cdot\!\sph$.

\end{enumerate}

\section{More on \texorpdfstring{$F_2$}{F\_2} inside \texorpdfstring{$\rot$}{SO(3)}}

In our proof of the Hausdorff--Banach--Tarski paradox, we used the
fact $\mathrm{SO}_{\mathbb{R}}(3)$ contains the free group $F_{2}$
on two generators. Indeed, the rotations by $\arctan\frac{1}{3}$
around the $x$-axis and around the $y$-axis generate a free group. However,
there is a cooler way to prove this fact. The proof
uses a bit of algebraic geometry, but we encourage the reader to attempt
to read it even if they don't know much about it, and take what they
can get out of it. We will also use the Haar measure, which we only
define formally a bit later on.

In order to prove $\mathrm{SO}_{\mathbb{R}}(3)$ contains a free group,
we will actually prove something stronger: that almost every pair
of elements in $\mathrm{SO}_{\mathbb{R}}(3)$ generates a free group.
Using the fact $\mathrm{SO}_{\mathbb{R}}(3)$ is ``Zariski dense''
in $\mathrm{SO}_{\mathbb{C}}(3)$, we will show that it is enough
to prove the same thing about $\mathrm{SO}_{\mathbb{C}}(3)$. The
funny thing about the proof is that, in order to show almost every
pair of elements in $\mathrm{SO}_{\mathbb{C}}(3)$ generates a free
group, we will see it is actually enough to prove that \emph{one }pair
of elements does this, which is the opposite of what we did in the
case of $\mathrm{SO}_{\mathbb{R}}(3)$. The point is that there is a very nice
way to find a free group in $\mathrm{SO}_{\mathbb{C}}(3)$, using what is
called ``the ping-pong lemma".

Recall that a subset of a topological space is \emph{nowhere dense
}if its closure has empty interior, and \emph{meagre }if it is a countable
union of nowhere dense subsets.
\begin{thm}
	Almost every pair of elements in $\mathrm{SO}_{\mathbb{R}}(3)$ generates
	a free group. More precisely, $\left\{ (x,y):\left\langle x,y\right\rangle \ncong F_{2}\right\} $
	is meagre and of Haar measure zero inside $\mathrm{SO}_{\mathbb{R}}(3)\times\mathrm{SO}_{\mathbb{R}}(3)$.
\end{thm}

\begin{rem}
	We will in fact see that almost every pair of elements $(x,y)$ ``satisfies
	no relations'', in the sense that the natural homomorphism $F_{2}\to\mathrm{SO}_{\mathbb{R}}(3)$
	sending the generators of $F_{2}$ to $x$ and $y$ is injective.
\end{rem}

\begin{proof}
	Denote the free generators of $F_{2}$ by $a$ and $b$. Every element
	$w=w(a,b)\in F_{2}$ gives us a map
	\begin{align*}
		w:\mathrm{SO}_{\mathbb{R}}(3)\times\mathrm{SO}_{\mathbb{R}}(3) & \to\mathrm{SO}_{\mathbb{R}}(3)\\
		(x,y) & \mapsto w(x,y),
	\end{align*}
	sending $(x,y)$ to the word obtained by replacing $a$ and $b$ in
	$w(a,b)$ with $x$ and $y$ respectively. Consider 
	\[
	\Omega_{w}\coloneqq w^{-1}(1)=\left\{ (x,y)\in\mathrm{SO}_{\mathbb{R}}(3)^{2}:w(x,y)=1\right\} .
	\]
	For every $w\in F_{2}$, $\Omega_{w}$ is a closed subset of $\mathrm{SO}_{\mathbb{R}}(3)\times\mathrm{SO}_{\mathbb{R}}(3)$.
	We will show that, if $w\neq 1$, then $\Omega_w$ has empty interior and is of measure zero.
	Since $F_{2}$ is countable, it will follow that $\bigcup_{w\in F_{2}\backslash\left\{ 1\right\}}\Omega_{w}$
	is meagre and of measure zero, and the proof will be done.
	
	Now comes the part where we use a bit of algebraic geometry. Since
	the group $\mathrm{SO}_{\mathbb{R}}(3)$ is defined by algebraic equations,
	it is a real algebraic group, and its topology as an algebraic group is
	called the Zariski topology. The map $w:\mathrm{SO}_{\mathbb{R}}(3)^{2}\to\mathrm{SO}_{\mathbb{R}}(3)$
	is also an algebraic map, so it follows that $\Omega_{w}$ is an algebraic
	subvariety of $\mathrm{SO}_{\mathbb{R}}(3)^{2}$. It is therefore
	either everything or of  lower dimension. If it is of lower dimension, then it is nowhere dense and of measure zero. Thus,
	it is enough for us to show that $\Omega_{w}\neq\mathrm{SO}_{\mathbb{R}}(3)\times\mathrm{SO}_{\mathbb{R}}(3)$
	for every $w\in F_{2}\backslash\left\{ 1\right\}$. 
	
	Now, the group $\mathrm{SO}_{\mathbb{C}}(3)$, which contains $\mathrm{SO}_{\mathbb{R}}(3)$,
	is a complex algebraic group with a Zariski topology, and it is not
	hard to see that $\mathrm{SO}_{\mathbb{R}}(3)$ is dense in $\mathrm{SO}_{\mathbb{C}}(3)$
	in this topology (one says that $\mathrm{SO}_{\mathbb{R}}(3)$ is
	``Zariski dense'' in $\mathrm{SO}_{\mathbb{C}}(3)$). Moreover,
	we may consider $w$ also as a map from $\mathrm{SO}_{\mathbb{C}}(3)\times\mathrm{SO}_{\mathbb{C}}(3)$
	to $\mathrm{SO}_{\mathbb{C}}(3)$ in the exact same way, and set $\tilde{\Omega}_{w}=w^{-1}(1)$.
	Observe that $\tilde{\Omega}_{w}$ contains $\Omega_{w}$. Therefore
	(since $\mathrm{SO}_{\mathbb{R}}(3)$ is Zariski dense in $\mathrm{SO}_{\mathbb{C}}(3)$),
	if it were the case that $\Omega_{w}=\mathrm{SO}_{\mathbb{R}}(3)^{2}$,
	it would follow that $\tilde{\Omega}_{w}=\mathrm{SO}_{\mathbb{C}}(3)^{2}$.
	Thus, it is enough for us to show (for every $w\in F_{2}\backslash\left\{ 1\right\}$) that $\tilde{\Omega}_{w}$
	is not all of $\mathrm{SO}_{\mathbb{C}}(3)^{2}$.
	
	Therefore, it is enough for us to show that $\mathrm{SO}_{\mathbb{C}}(3)$
	contains a free group, because, if $x,y\in\mathrm{SO}_{\mathbb{C}}(3)$
	generate a free group, it clearly follows that $(x,y)\notin\tilde{\Omega}_{w}$
	for every $w\in F_{2}\backslash\left\{ 1\right\} $.
	
	Thus, we wish to find two elements in $\mathrm{SO}_{\mathbb{C}}(3)$ generating $F_2$.
	(Of course, Exercise~\ref{f2:so3} gives us such a subgroup, but the whole point of the present proof is to find a different way of proving there exists such a subgroup.)
	
		Consider the following:
	$$\begin{matrix}
		&\mr{SO}_{\mb C}\l(3\r) &= &\l\{A\in M_3\l(\mb C\r):A^t\matI A=\matI,\;\det A=1\r\}\\
		
		\cong \;&\mr{SO}_{\mb C}\l(2,1\r) &= &\l\{A\in M_3\l(\mb C\r):A^t\mat1 A=\mat1 \! ,\;\det A=1\r\}\\
		
		\supseteq\; &\mr{SO}^\circ_{\mb R}\l(2,1\r) &= &\l\{A\in M_3\l(\mb R\r):A^t\mat1 A=\mat1 \! ,\;\det A=1\r\}^\circ\\
		
		\cong &\psl &=  &\mr{SL}_2\l(\mb R\r)/\l\{\pm 1\r\}
	\end{matrix}$$
	
	\begin{exe}
		Prove that $\mr{SO}_{\mb C}\l(3\r)\cong \mr{SO}_{\mb C}\l(2,1\r)$.
	\end{exe}
	\begin{hint}
		Over $\mb C$, all non-degenerate quadratic forms are isomorphic.
	\end{hint}
	\begin{exe}
		Prove that $\mr{SO}^\circ_{\mb R} \l(2,1\r)$ (the connected component of the identity element of $\mr{SO}_{\mb R} \l(2,1\r)$) is isomorphic to $\psl$.
	\end{exe}
	\begin{hint}
		\def \g {\mf{g}}
		Consider
		\begin{align*}
			G&\coloneqq \mr{SL}_2\l(\mb R\r)=\l\{A\in M_2\l(\mb R\r):\det A=1\r\},\\
			\g &\coloneqq \mathfrak{sl}_2\l(\mb R\r)=\l\{X\in M_2\l(\mb R\r):\mr{tr}X=0\r\}.
		\end{align*}		
		The map $\Ad:G\to\ \mr{GL}\l(\g\r)$ defined by $A\mapsto \Ad A$, where $\Ad A\l(X\r)=AXA^{-1}$, gives us an action of $G$ on $\g$ (since $\mr {tr}\l(AXA^{-1}\r)=\mr {tr}\l(AA^{-1}X\r)=0$). The kernel $\ker \l(\Ad\r)=\l\{\pm1\r\}$ is exactly the centre. Show that $\Ad\l(G\r)=\mr{SO}^\circ_{\mb R}\l(2,1\r)$. On $\g$, take the determinant form (which is a quadratic form of signature $2,1$ which is preserved by $\Ad\l(G\r)$).
		\def \g {\gamma}
	\end{hint}

Therefore, it is enough for us to find a subgroup isomorphic to $F_2$ inside $\psl$ (since $\psl$ is isomorphic to $\mr{SO}_{\mb R}\l(2,1\r)$, which is contained in $\mr{SO}^\circ_{\mb C}\l(2,1\r) $, which is isomorphic to $\mr{SO}_{\mb C}\l(3\r)$).

This will be accomplished in the next section, using the ping-pong lemma, which is the most common tool for finding free subgroups.\end{proof}

\begin{rem}
	Observe that we have proved the following: if a connected Lie group contains a free subgroup, then almost every pair of elements ``satisfies no relations" (and, in particular, generate a free subgroup). Then we used the fact that this property goes down to Zariski dense subgroups in order to deduce that, if $\mr{SO}_{\mb C}(3)$ contains a free subgroup, so does $\mr{SO}_{\mb R}(3)$. 
\end{rem}

\section{The Ping-Pong Lemma}

A useful tool for proving $F_2$ is contained in certain groups is the ping-pong lemma.
\subsection{First Formulation}
\begin{lem}[The Ping-Pong Lemma]
Let $\Gamma$ act on $X$. Suppose we have elements $a,b$ in $\Gamma$ and disjoint subsets $A^+,A^-,B^+,B^-$ in $X$ such that:
\begin{align*}
a.\l(A^+\cupdot B^-\cupdot B^+\r)\subseteq A^+,\\
a^{-1}.\l(A^-\cupdot B^-\cupdot B^+\r)\subseteq A^-,\\
b.\l(B^+\cupdot A^-\cupdot A^+\r)\subseteq B^+,\\
b^{-1}.\l(B^-\cupdot A^-\cupdot A^+\r)\subseteq B^-.\\
\end{align*}
Then $\l<a,b\r>=F_2$.
\end{lem}
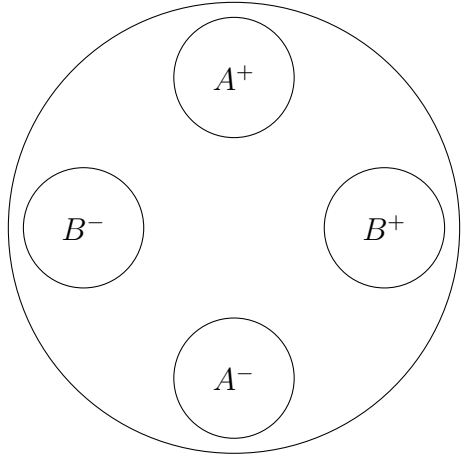
\begin{center}
\begin{tikzpicture}
\draw (0,0) circle (3cm);
\draw (0,2) node {$A^+$} circle (0.8cm);
\draw (0,-2) node {$A^-$} circle (0.8cm);
\draw (2,0) node {$B^+$} circle (0.8cm);
\draw (-2,0) node {$B^-$} circle (0.8cm);
\end{tikzpicture}
\end{center}
\begin{proof}
We need to show that every nonempty word in $a,b$ is a nontrivial element, so let $u$ be such a word. We can assume without loss of generality that $u_1=a$, where $u_1$ is the first letter of $u$. Take any $x\in B^+,y\in B^-$. Then $u_1$ sends both $x$ and $y$ to $A^+$. The next letter in $u$, call it $u_2$, can be anything except $u_1^{-1}=a^{-1}$, so it follows from our equations that it will necessarily send $A^+$ to a certain $\Omega_2\in \{A^-,B^+,B^-\}$ (and, in particular, it will send both $x$ and $y$ to $\Omega_2$). Observe that any $\xi \in \l\{a,a^{-1},b,b^{-1}\r\}$ sends the elements in all but one of the sets $\l\{A^+,A^-,B^+,B^-\r\}$ to a certain $\Omega_\xi\in\l\{A^+,A^-,B^+,B^-\r\}$, and observe that any $\eta\neq\xi^{-1}$ ($\eta\in\l\{a,a^{-1},b,b^{-1}\r\}$) sends $\Omega_\xi$ to some $\Omega_\eta\in\l\{A^+,A^-,B^+,B^-\r\}$. Thus, we can repeat this observation any finite number of times, so we can deduce that $x,y$ live together happily every after, or at least that there is some $\Omega_u\in\l\{A^+,A^-,B^+,B^-\r\}$ such that $u.x,u.y\in\Omega_u$. Since $x\in B^+$ and $y\in B^-$, it follows that $u$ is a nontrivial element.
\end{proof}
\begin{exa}
\label{lalala}
Consider the action of $\mr{SL}_2\l(\mb{R}\r)$ on the projective line $\mb{RP}^1$ (the space of lines in $\mb{R}^2$ that go through the origin) defined as follows: $a\cdot\l[x\r]=\l[ax\r]$, where $[x]$ denotes the  line $\l\{\lambda x:\lambda \in\mb R\r\}$. It's easy to see that this is indeed an action.

Consider $a=\mathrm{diag}\l(n,n^{-1}\r)$ for large $n$. A line which isn't very close to the $y$-axis becomes a line which is very close to the $x$-axis. Conjugate by a rotation of angle $\pi/4$ and obtain $b$; $b$ has similar dynamics to $a$, only around a different axis. (See Figure \ref{Figure} below.)

We have: 
\begin{align*}
a\l(\mb{RP}^{1}\setminus A^-\r)&\subseteq A^+,\\
b\l(\mb{RP}^{1}\setminus B^-\r)&\subseteq B^+,\\
a^{-1}\l(\mb{RP}^{1}\setminus A^+\r)&\subseteq A^-,\\
b^{-1}\l(\mb{RP}^{1}\setminus B^+\r)&\subseteq B^-,
\end{align*}
and hence $\l<a,b\r>\cong F_2$ by the ping-pong lemma.
\end{exa}
\begin{figure}
\centering
\caption{}
\label{Figure}

\begin{tikzpicture}
\draw (4,0) circle (6cm);
\draw[fill, color=black] (4,0) circle (0.05);
\path (4,6) coordinate (up);
\path (4,-6) coordinate (down);
\path (-2,0) coordinate (left);
\path (10,0) coordinate (right);
\path (8.2426,4.2426) coordinate (up-right);
\path (8.2426,-4.2426) coordinate (down-right);
\path (-0.2426,4.2426) coordinate (up-left);
\path (-0.2426,-4.2426) coordinate (down-left);
\draw[color=gray] (up) -- (down);
\draw[color=gray] (left) -- (right);
\draw[color=gray] (down-left) -- (up-right);
\draw[color=gray] (down-right) -- (up-left);

\path (5.04189, 5.90885) coordinate (a2);
\path (2.95811, -5.90885) coordinate (a3);
\draw[color=red] (a2)  -- (a3);

\path (2.95811, 5.90885) coordinate (a1);
\path (5.04189, -5.90885) coordinate (a4);
\draw[color=red] (a1) node[below right] {$A^-$} -- (a4)node[above left] {$A^-$};

\path (-1.90885, 1.04189) coordinate (a5);
\path (9.90885, -1.04189) coordinate (a8);
\draw[color=red] (a5) -- (a8);

\path (-1.90885, -1.04189) coordinate (a5);
\path (9.90885, 1.04189) coordinate (a8);
\draw[color=red] (a5)node[above right] {$A^+$} -- (a8)node[below left] {$A^+$};

\path ({7.44146, 4.91491}) coordinate (b1);
\path ({0.558541, -4.91491}) coordinate (b4);
\draw[color=blue] (b1)   -- (b4);

\path ({8.91491, 3.44146}) coordinate (b2);
\path ({-0.914912, -3.44146}) coordinate (b3);
\draw[color=blue] (b2) node[left] {$B^+$} -- (b3)node[right] {$B^+$};

\path ({0.558541, 4.91491}) coordinate (b5);
\path ({7.44146, -4.91491}) coordinate (b8);
\draw[color=blue] (b5) -- (b8);

\path ({-0.914912, 3.44146}) coordinate (b6);
\path ({8.91491, -3.44146}) coordinate (b7);
\draw[color=blue] (b6)node[right] {$B^-$} -- (b7) node[below left] {$B^-$};

\end{tikzpicture}

\end{figure}
\begin{exa}
Let $\mr{SL}_2\l(\mb C\r)$ act on $\hat {\mb C}$ by M\"obius transformations. $$\pmtb a &b\\ c &d\pmte z=\frac{az+b}{cz+d}.$$Recall M\"obius transformations send generalised circles to generalised circles.

Consider the sets $A^+, A^-,B^+,B^-$ and the maps $f,g$ defined by:

\

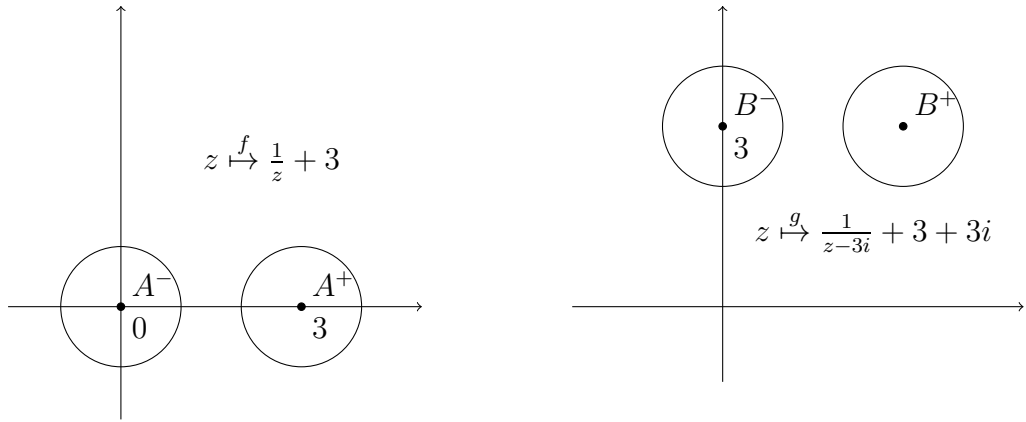
\begin{tikzpicture}
\pa (0,0) \at (O);
\pa (-1.5,0) \at (l);
\pa (0,-1.5) \at (d);
\pa (4,0) \at (r);
\pa (0,4) \at (u);
\pa (2.4,0) \at (a);

\dr[->] (l) -- (r);
\dr[->] (d) -- (u);
\dr[fill] (O)node[below right] {$0$} circle (0.05);
\dr (O)node[above right] {$A^-$} circle (0.8);
\dr[fill] (a)node[below right] {$3$} circle (0.05);
\dr (a)node[above right] {$A^+$} circle (0.8);
\dr (2,2)node{$z\overset{f}{\mapsto}\frac{1}{z}+3$};

\pa (8,0) \at (O);
\pa (6,0) \at (l);
\pa (8,-1) \at (d);
\pa (12,0) \at (r);
\pa (8,4) \at (u);
\pa (8,2.4) \at (bm);
\pa (10.4,2.4) \at (bp);

\dr[->] (l) -- (r);
\dr[->] (d) -- (u);
\dr[fill] (bm)node[below right] {$3$} circle (0.05);
\dr (bm)node[above right] {$B^-$} circle (0.8);
\dr[fill] (bp) circle (0.05);
\dr (bp)node[above right] {$B^+$} circle (0.8);
\dr (10,1)node{$z\overset{g}{\mapsto}\frac{1}{z-3i}+3+3i$};
\end{tikzpicture}

The maps $f,g$ satisfy the ping-pong lemma:
\begin{align*}
f(\hat {\mb C}\setminus A^-)&\subseteq A^+\\
f^{-1}(\hat {\mb C}\setminus A^+)&\subseteq A^-\\
g(\hat {\mb C}\setminus B^-)&\subseteq B^+\\
g^{-1}(\hat {\mb C}\setminus B^+)&\subseteq B^-
\end{align*}
And hence $\l<f,g\r>\cong F_2$.
\end{exa}

\subsection{Second Formulation}
There is a second formulation to the ping-pong lemma, which is more general:

\begin{lem}[Second Formulation of the Ping-Pong Lemma]
Let a group $\Gamma$ act on a set $X$. If there exist $G,H\le\Gamma$, $\l| G\r|\ge 3$, and disjoint nonempty $A,B\subseteq X$ such that $gB\subseteq A$ and $hA\subseteq B$ for all nontrivial $g\in G$, $h\in H$  respectively, then $\l<G,H\r>\cong G*H$\ (the free product of $G$ and $H$).
\end{lem}
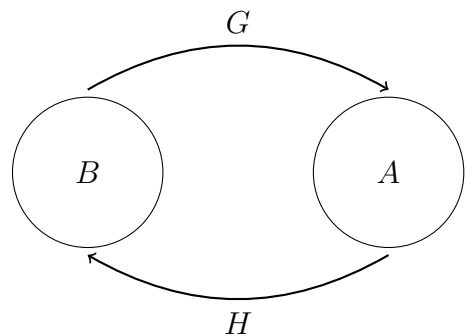
\begin{center}
\begin{tikzpicture}

\pa (0,0);
\pa (4,0) \at (A);
\pa (0,0) \at (B);

\dr (A)node {$A$} circle (1);
\dr (B)node {$B$} circle (1);

\pa[thick, ->] (0,1.1) edge [bend left] (4,1.1);
\pa[thick, ->] (4,-1.1) edge [bend left] (0,-1.1);

\dr (2,2) node {$G$};
\dr (2,-2) node {$H$};

\end{tikzpicture}
\end{center}
\begin{proof}
Every nontrivial word in $\l<G,H\r>$ is conjugated to a word starting and ending with elements in $G$:
\begin{align*}
h_1g_1\cdots h_ng_n&\mapsto \t g^{-1}h_1g_1\cdots h_n g_n \t g,\quad\tilde g\ne g_n^{-1}\\
g_1h_1\cdots g_n h_n&\mapsto \t g^{-1}g_1h_1\cdots g_n h_n \t g,\quad \t g\ne g_1^{-1}\\
h_1g_1\cdots g_{n-1}h_n&\mapsto \t g^{-1} h_1 g_1\cdots g_{n-1} h_n \t g.
\end{align*}
(Observe we're using the fact $\l|G\r|\ge 3$). A word starting and ending with elements in $G$ must be nontrivial because it sends $B$ into $A$, so our proof is done.
\end{proof}
\begin{rem}
The requirement $\l|G\r|\ge 3$ is necessary; take $X=\mb Z/{2\mb Z}=\l\{\bar 0,\bar 1\r\}$, $A=\l\{\bar 0\r\}$, $B=\l\{\bar 1\r\}$, $\Gamma=G=H=\mb Z/{2\mb Z}$.
\end{rem}
\begin{exe}
Deduce the first formulation of the ping-pong lemma from this one. (Check $\l<a\r>\cong\mb Z\cong \l<b\r>$, and observe $\mb Z*\mb Z=F_2$).
\end{exe}
\begin{exa}
Let $X=\Gamma=F_2$, acting on itself by left translations.

Denote by $a,b$ the generators of $F_2$, and by $F_x$ the set of words in $F_2$ starting with $x$. Define
\begin{align*}
A^+&=F_a,\\
A^-&=F_{a^{-1}},\\
B^+&=F_b,\\
B^-&=F_{b^{-1}},
\end{align*}
and $A=A^+\cup A^-$, $B=B^+\cup B^-$. So, by the second formulation of the ping-pong lemma we get: $F_2=\l<a,b\r>\cong\l<a\r>*\l<b\r>\cong F_2$.
\end{exa}
\begin{exe}
$\psl$ acts on the projective line $\mb P^1$; if $L_{\l(x,y\r)}=\l\{\l(tx,ty\r\}:t\in\mb R\r\}\in\mb P^1$ and $\l(\bsmt a &b\\ c &d\esmt\r)\in\psl$ then $\l(\bsmt a &b\\ c &d\esmt\r)L_{\l(x,y\r)}=L_{\l(\bsmt a &b\\ c &d\esmt\r)\l(\bsmt x\\ y\esmt\r)}$.

Show that $\l<\l(\bsmt 1 &2\\0 &1\esmt\r),\l(\bsmt1 &0\\2 &1\esmt\r)\r>=F_2$.
\end{exe}
Consider $\mr{SL}_2\l(\mb Z\r)=\l\{A\in M_2\l(\mb Z\r):\det A=1\r\}$.
\begin{exe}
Check $\mr{SL}_2\l(\mb Z\r)$ is a group.
\end{exe}
By the solution to the previous exercise, $F_2$ is isomorphic to a subgroup of $\mr{SL}_2\l(\mb Z\r)$.
\begin{exe}[Important]
Show $W=\l(\bsmt0 &1\\-1 &0 \esmt\r)$, $U=\l(\bsmt1 &1\\0 &1\esmt\r)$ generate $\mr{SL}_2\l(\mb Z\r)$.
\end{exe}
\begin{hint}
Show that if $x=\l(a,b\r)$ is a primitive vector (i.e., $\gcd\l(a,b\r)=1$)\footnote{$\gcd$ stands for the \emph{\textbf{greatest common divisor}}.} then there exist a sequence $\l\{r_i\r\}_{i=1}^n\subseteq\l\{U,W\r\}$ such that $xr_1r_2\cdots r_n=\l(1,0\r)$. In order to do this you can apply the Euclidean algorithm.

Show $\l\{g\in\ \mr{SL}_2\l(\mb Z\r):g\l(\bsmt1\\0\esmt\r)=\l(\bsmt1\\0\esmt\r)\r\}=\l<U\r>$.
\end{hint}
\begin{exe}
Let $T=WU=\l(\bsmt0 &1\\-1 &-1\esmt\r)$. Compute the order of $T$\ and the order of $W$ (as elements in $\mr{SL}_2\l(\mb Z\r)$), and prove that $\l<W,T\r>=\mr{SL}_2\l(\mb Z\r)$.
\end{exe}
\begin{exe}
Denote by $\bar W, \bar T$ the corresponding matrices in $\pslz=\mr{SL}_2\l(\mb Z\r)/\l\{\pm1\r\}$. Then $\bar T$ is of order $3$ and $\bar W$ is of order $2$. Show that $\l<\bar T\r>$, $\l<\bar W\r>$ satisfy the second formulation of the ping-pong lemma.
\end{exe}
\begin{hint}Consider the following diagram:

\
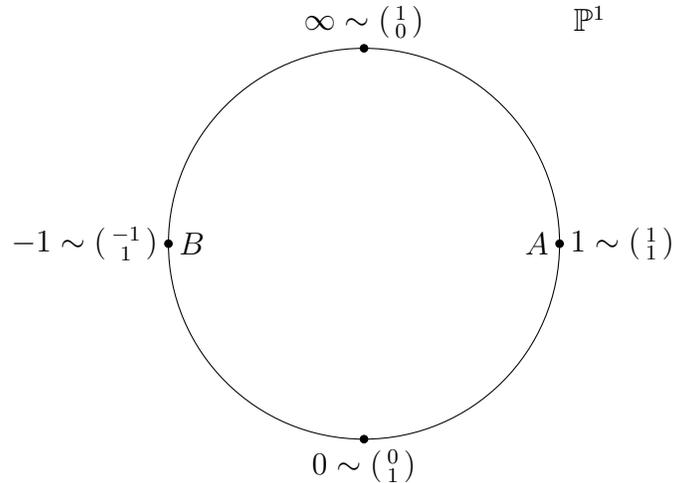
\begin{center}
\begin{tikzpicture}
\pa (0,0) \at (O);
\pa (0,2.6) \at (u);
\pa (0,-2.6) \at (d);
\pa (-2.6,0) \at (l);
\pa (2.6,0) \at (r);

\dr (O) circle (2.6);

\dr[fill] (u)node[above]{$\infty\sim\l(\bsmt1\\0\esmt\r)$} circle (0.05);
\dr[fill] (d)node[below]{$0\sim\l(\bsmt0\\1\esmt\r)$} circle (0.05);
\dr[fill] (l)node[left]{$-1\sim\l(\bsmt-1\\1\esmt\r)$} node[right]{$B$} circle (0.05);
\dr[fill] (r)node[right]{$1\sim\l(\bsmt1\\1\esmt\r)$} node[left]{$A$} circle (0.05);

\dr (3,3)node{$\mb P^1$};

\end{tikzpicture}
\end{center}
$A$ is the arc from $0$ to $\infty$ going anticlockwise, $B$ the arc going clockwise. Then $W$ is a rotation by $180$ degrees, which switches them.

Compute what $T$ does and deduce, using the second formulation of the ping-pong lemma, that $\pslz\cong C_3*C_2$. If you know what an amalgamated product is, you may also show that $\mr{SL}_2\l(\mb Z\r)$ is isomorphic to an amalgamated product of two finite groups.
\end{hint}

\section*{Further Reading}To explore the original papers on this paradox by Hausdorff and by Banach and Tarski, the readers are referred to \cite{hausdorff1914bemerkung} and \cite{banach1924decomposition}. Those interested in delving deeper into this paradox and its variations (e.g., in the hyperbolic space) may refer to \cite{tomkowicz2016banach}.

\chapter{Amenability}

\section{Definition}
In Chapter 1, we showed that there is no finitely additive probability measure defined on all subsets of $\mb S^n$ that is invariant under rotations, if $n$ is at least $2$. But what if $n=1$?
\begin{question}
Is there a finitely additive probability measure defined on all subsets of $\mb S^1$ that is invariant under rotations?
\end{question}
\begin{thm}[Banach]
Yes.
\end{thm}
The property of admitting such a measure is called \textit{amenability}, and this chapter is mostly devoted to investigating different, equivalent definitions of it. Amenability has countless applications, as we will see.
\begin{defn}
A group $\Gamma$ is called \emph{\textbf{amenable}} if it admits a finitely additive left-invariant (i.e., invariant under left translations) probability measure defined on all its subsets.
\end{defn}
\begin{exa}
Every finite group is amenable (give each singleton the measure $\l|G\r|^{-1}$).
\end{exa}

\begin{rem}
There is nothing special about ``left", of course. If $m$ is a left-invariant, finitely additive measure, then $$v\l(A\r)\coloneqq m\l(A^{-1}\r),$$ where $A^{-1}=\l\{a^{-1}\mi a\in A\r\}$, is a right invariant, finitely additive measure.
\end{rem}

\section{Means}
In this section we define \textit{means} and show the existence of a left-invariant mean is equivalent to the existence of a left-invariant, finitely additive probability measure.

Consider $\lig=\l\{f:\Gamma\to\mb R\mi\sup_{\gamma\in\Gamma}\l|f\l(\gamma\r) \r|<\infty\r\}$. It is a linear space, and it has a natural norm, $\l\|f\r\|=\sup_{\gamma\in\Gamma}\l|f\l(\gamma\r)\r|$. In fact, this is a Banach space (i.e., it is complete with respect to this norm).

If $B$ is a normed vector space over $\mb F$, where $\mb F=\mb C$ or $\mb F=\mb R$, we denote by $B^*$ the space of all continuous linear functionals $\varphi:B\to\mb F$ (with respect to the norm above). It is a linear space, and it has a norm $\l\|\varphi\r\|=\sup_{f\in B,\l\|f\r\|=1}\l|\varphi\l(f\r)\r|$. In fact, if $B$ is a Banach space, then $B^*$ is also a Banach space with respect to this norm.
\begin{rem}
We can replace $\infty$ with $p$ for $p\ge 1$, and get $$\lpg=\l\{f:\G\to\mb R\mi (\sum_{\g\in\G}\l|f\l(\g\r)\r|^p)^{1/p}<\infty\r\}$$ with the norm $\l\|f\r\|=\l(\sum_{\g\in\G}\l|f\l(\g\r)\r|^p\r)^{1/p}$.\footnote{Note that if $\sum\l|f\l(\g\r)\r|^p<\infty$ then $\{\g|f\l(\g\r)\ne 0\}$ is countable.}  This is also a Banach space.
\end{rem}
\begin{defn}
Let $\G$ be a group. A linear functional $\varphi\in\lig^*$ is a \emph{\textbf{mean}} on $\Gamma$ if $\l\|\varphi\r\|=1$ and $\varphi$  is positive (i.e., $\varphi\l(f\r)\ge0$ for all positive functions $f\in\lig$).
\end{defn}
\begin{exe}
Show $\varphi\in\lig^*$ is a mean if, and only if, $\varphi$ is positive and $\varphi\l(1\r)=1$, where $1\in\lig$ is the constant function $f\equiv 1$.
\end{exe}
$\Gamma$ acts on $\lig$ by left translation: $$L_\gamma(f)\l(\alpha\r)=f\l(\gamma^{-1}.\alpha\r)$$ for all $\gamma\in\Gamma$, $f\in\lig$, $\alpha\in\Gamma$. This is called the \emph{\textbf{left regular action}}. A linear functional $\varphi\in\lig^*$ is \emph{\textbf{left-invariant}} if $\varphi\l(L_\gamma f\r)=\varphi\l(f\r)$  for all $f\in\lig$ and for all $\gamma\in\Gamma$.
\begin{prop}
$\Gamma$ is amenable if and only if there is a left invariant mean on $\ell_\infty(\G)$.
\end{prop}

\begin{proof}
Given a left-invariant mean we define a left-invariant measure, and vice versa.

If $\varphi\in\lig^*$ is a left-invariant mean, define $m:\mathcal P\l(\Gamma\r)\to\l[0,1\r]$ by $A\mapsto \varphi\l(1_A\r)$.

 If $m$ is a left-invariant, finitely additive probability measure on $\Gamma$, define $\varphi$ as follows: set $\varphi\l(1_A\r)=m\l(A\r)$, extend linearly to linear combinations of indicator functions, and then for general $f\in\lig$  write $f$ as a limit of a finite-valued function sequence $\l(f_n\r)_{n=1}^\infty$ and set $\varphi\l(f\r)=\lim_{n\to\infty}\varphi\l(f_n\r)$.
\end{proof}
\begin{exe}
Show $\varphi$ is well defined, and finish the details of the proof (it's straightforward).
\end{exe}

Of course, it's easy to see in the same way that amenability is equivalent to the existence of a right-invariant mean. In fact, a much stronger statement holds:
\begin{exe}
Show that, if a group $\G$ is amenable, then it has a 2-sided invariant mean.
\end{exe}
\subsection{Banach--Alaoglu Theorem}
If $B$ is a normed vector space over $\mb F=\mb R$ or $\mb F=\mb C$, we define the \emph{\textbf{weak topology on $B$}} to be the weakest topology such that all elements in $B^*$ are continuous. (Recall that $B^*$ is the space of continuous linear functionals $f: B\to \mb F$, with respect to the topology on $B$ induced by its norm.)

In the real case, a base for the weak topology is $$\l\{\varphi^{-1}\l(I\r)\mi\varphi\in B^*,I\subseteq \mb R\text{ is an open interval}\r\}.$$

We say a sequence $\l(x_n\r)_{n=1}^\infty$ \emph{\textbf{converges weakly to $x$}}, and denote $x_n\overset{w}{\too }x$, if $x_n$ converges to $x$ in the weak topology. This is equivalent to the condition that $\varphi\l
(x_n\r)\too\varphi\l(x\r)$ for every $\varphi\in B^*$.
\begin{exe}
Consider $\ltz$ and $\delta_n\in \ltz$ given by $\delta_n=1_{\l\{n\r\}}$. Show that $\delta_n\wkly 0$.
\end{exe}

We define the \emph{\textbf{weak-* topology}} on $B^*$ to be the weakest topology such that all evaluation functions $\delta_x$ for $x\in B$, defined by $\varphi\mapsto \varphi\l(x\r)$,  are continuous. In the real case, a base for this topology is $$\left\{\delta_x^{-1}\l(I\r) \middle| x\in B, I\subseteq \mb R\text{ is an open interval}\right\}.$$
\begin{rem}
Note that $B^*$ has both the weak topology (the weakest topology such that all the elements in $B^{**}$ are continuous) and the weak-* topology (the weakest topology such that all the evaluation functions are continuous). Clearly, the weak topology is stronger than the weak-* topology (since the evaluation functions are also elements in $B^{**}$).
\end{rem}
\begin{rem}
Note that a sequence $\l(\varphi_n\r)_{n=1}^\infty$ in $B^*$ converges to $\varphi$ in the weak-* topology  if and only if $\varphi_n\l(x\r)\too\varphi\l(x\r)$ for all $x\in B$.
\end{rem}
\begin{thm}[Banach–Alaoglu Theorem]
If $X$ is a normed vector space over $\mb R$ or $\mb C$, the unit ball of $X^*$ is weak-* compact (i.e., compact with respect to the weak-* topology).
\end{thm}
\begin{proof}[Sketch of proof]
Embed $B\l(X^*\r)$ (the unit ball in $X^*$) in $I^{B\l(X\r)}$ (the set of functions from the unit ball in $X$ to the interval $I=\l[-1,1\r]$), by\ $$\varphi\mapsto \l(\varphi\l(x\r)\r)_{x\in B\l(X\r)}.$$

We need to show:
\begin{enumerate}
\i This function is injective.
\i It is a homeomorphism onto its image.
\i The image is closed.
\end{enumerate}
And then the theorem follows from Tychonoff's theorem.
\end{proof}
\begin{exe}
Write a detailed proof of the theorem.
\end{exe}
\begin{exa}
        The space of means on a group $\G$, as well as the space of left-invariant means, is weak-* compact (as a subset of $\ell_\infty(\G)^*$).
\end{exa}
\begin{proof}
The set of means on $\lig$ is a subset of the unit ball of $\lig^*$. This set is weak-* closed, because $m\in \lig ^*$ is a mean if and only if $m\ge 0$\ (i.e., $m\l(f\r)\ge 0$ for all $f\ge 0$) and $m\l(1\r)=1$, and these conditions are weak-* closed (because if $m_n\overset{w\text{-}*}{\too}m$ then $m\l(f\r)=\lim m_n\l(f\r)$). Hence, the space of means is weak-* compact by the Banach–Alaoglu theorem.

We want to prove the subspace of left-invariant means is also weak-* closed and hence compact. First, we go over the definitions again:

$\G$ acts on $\lig$ by left translations. That is, for all $f\in \lig$, $x\in \G$, $$L_\g f\l(x\r)\coloneqq f\l(\g^{-1}x\r).$$$\G$ acts on $\lig^*$ too, by $$\gamma.m\l(f\r)=m\l(L_{\g^{-1}}\l(f\r)\r).$$ So $m$ is left-invariant $\iff$ $\g. m=m$ $\forall\g\in\G$ $\iff$ $m\l(f\r)\l(x\r)=m\l(L_{\g^{-1}}f\r)\l(x\r)$ $\forall f\in\lig,\g\in\G,x\in \G$. This is a weak-* closed condition, so the subspace of left-invariant means is closed and hence compact.
\end{proof}

In general, a Banach space $B$ can be embedded in $B^{**}$ by evaluation maps ($\delta_x\l(\varphi\r)=\varphi\l(x\r)$ for all $x\in B$, $\varphi \in B^*$). Consider $\leg=\l\{f:\G\to\mb R\mi \sum_{\g\in\G}\l|f\l(\g\r)\r|<\infty\r\}$.

\begin{exe}
\label{l1:means}
By the remarks above, $\leg$ can be embedded in $\lig^*$. Which functions correspond to means?
\end{exe}
\section{Asymptotically Invariant Nets}
We know $\G$ is amenable if and only if the space of left invariant means is nonempty. We want to find more such conditions. One such condition is the existence of \textit{asymptotically invariant nets}. First, we need to introduce the notion of nets.
\begin{rem}
Usually, we will only deal with countable groups, where it is enough to consider asymptotically invariant sequences (and later, F\o lner sequences), rather than nets. Therefore, if you  don't want to be bothered with nets, you will not lose much (from a practical point of view) by simply assuming all nets are sequences (and all subnets are subsequences).\footnote{Later, in the last chapter of this part as well as the second part of the book, we will deal mostly with \textit{topological} groups, which will usually not be countable. However, we will then deal mostly with second-countable groups, where it is again enough to consider sequences (rather than general nets).}
\end{rem}
\subsection{Nets}
A net is a generalisation of the notion of a sequence. 

Recall that a sequence in a space $X$ is a function $f:\mb N\to X$. If $X$ is a topological space, a sequence $f:\mb N\to X$ converges to a point $x\in X$ if, for every open neighbourhood $U$ of $x$, there is some $N\in \mb N$ such   $f\l(n\r)\in U$ for all $n\ge N$.
In metric spaces, or more generally first-countable topological spaces, sequences offer equivalent definitions both for continuity of functions and for closedness and compactness of sets: a function $f$ is  continuous at $x$ if and only if, for every $\l(x_n\r)_{n\in\mb N}$ such that $x_n\too x$, the sequence $f\l(x_n\r)$ converges to $f\l(x\r)$; a set $A$ is closed if and only if, for every $\l(x_n\r)_{n\in \mb N}$ in $A$  such that $x_n\too x$, we have $x\in A$; a set $K$ is compact if and only if every sequence $\l(x_n\r)_{n\in \mb N}$ in $K$ admits a convergent subsequence.

Unfortunately (or excitingly), none of this holds in general topological spaces. This is the reason to use nets; if you replace the word `sequence' with the word `net` in each of these statements, it will be correct even in general topological spaces.

The way this is done is the following: we replace the set $\mb N$ with any directed set, where
\begin{defn}
A \emph{\textbf{directed set}} is a partially ordered set $J$ such that, for every $\alpha,\beta$ in $J$, there exists some $\gamma$ in $J$ for which $\alpha,\beta\le\gamma$.
\end{defn}

With this terminology, we have
\begin{defn}
A \emph{\textbf{net}} in a space $X$ is a function $f: J\to X$ for some directed set $J$. We frequently denote a net as $\l(x_j\r)_{j\in J}$ instead of as a function $x: J\to X$.
\end{defn}

For the statements we made earlier to make sense, we need to define what a \emph{sub}net is, as well. 
\begin{defn}
If $\l(x_j\r)_{j\in J}$ is a net, $K$ is a directed set, and $f:K\to J$ is a function that satisfies 
\begin{enumerate}[i.]
\i $k_1\le k_2\then f\l(k_1\r)\le f\l(k_2\r)$
\i $f\l(K\r)$ is cofinal in $J$: that is, every $a\in J$ admits some $b\in K$ such that $a\le f(b)$,
\end{enumerate} 
then $x\circ f=\l(x_{f\l(k\r)}\r)_{k\in K}$ is called a \emph{\textbf{subnet}} of $\l(x_j\r)_{j\in J}$.
\end{defn}
\begin{rem} 
Observe that, if we consider a sequence as a net, it has subnets which aren't sequences.
\end{rem}
The notion of convergence in nets is formalised as follows:
\begin{defn}
If $X$ is a topological space, a net $\l(x_j\r)_{j\in J}$ \emph{converges} to a point $x\in X$ if, for every open neighbourhood $U$ of $x$, there exists some $\alpha\in J$ such that $x_j\in U$ for every $j\ge a$.
\end{defn}

\subsection{Asymptotically Invariant Nets}
We want conditions equivalent to amenability. We already know that the existence of invariant means is one such; we can now define the notion of \textit{asymptotically invariant nets}, and show it has an invariant mean as a limit of a subnet.
\begin{defn}[Important]
Let $\G$ be a group. A net $\l(f_j\r)_{j\in J}$ in $\lpg$ of functions with norm $1$ is said to be \emph{\textbf{asymptotically invariant}}, or \emph{\textbf{almost invariant}}, if $\l\|L_\g\l(f_j\r)-f_j\r\|_p\too 0$ for all $\g\in\G$.

An asymptotically invariant net whose domain is $\mb N$ (with the usual order) is called an \emph{\textbf{asymptotically invariant sequence}}. 
\end{defn}
At times, we will abbreviate these terms as ``a.i.n." and ``a.i.s.", respectively.
\begin{exe}
Show there exists an a.i.n.\ in $\leg$ if and only if there exists an a.i.n. in $\ltg$. (In fact, this is equivalent to the existence of an a.i.n. in $\lpg$ for $p<\infty$).
\end{exe}
\begin{hint}
Use the Mazur map: $$M_{p,q} :\lpg\to\lqg,$$ for $1\le p,q<\infty$, defined by $$f\mapsto\mr{Sign}\l(f\r)\cdot \l|f\r|^{p/q}.$$ The Mazur map maps the unit sphere of $\lpg$ to the unit sphere of $\lqg$ uniformly continuously.
\end{hint}
\begin{prop}
A group $\G$ is amenable if and only if there exist an almost invariant net in $\ell_1(\G)$.
\end{prop}
\begin{proof}
First, assume there is an almost invariant net in $\ell_1(\G)$. We want to find a subnet (up to normalisation) that converges to an invariant mean. First, given an a.i.n.\ $\l(f_j\r)$, we want to construct a new a.i.n.\ of positive functions. The way to do this is to first replace it with $(\l\|f_j\r\|)$ (where $\l\|f\r\|(\g)=|f(\g)|$ for any $f\in\ell_1(\G)$), and then normalise it.

By the solution to Exercise~\ref{l1:means}, if $f\ge 0$\ and $\l\|f\r\|_1=1$ then the evaluation function $\delta_f$ is a mean. Hence, if $\l(f_j\r)$ is a positive, norm $1$ a.i.n., then $\l(\delta_{f_j}\r)$ is a net of means. We know $$(L_\g\l(f_j\r)-f_j)\too 0, \quad \forall\g\in\G,$$ so in particular $$\varphi\l(L_\g\l(f_j\r)-f_j\r)\too0,\quad\forall\varphi\in\lig,\g\in\G,$$ which means $$\l(\g^{-1}. \delta_{f_j}-\delta_{f_j}\r)\l(\varphi\r)=\varphi\l(L_\g\l(f_j\r)\r)-\varphi\l(f_j\r)\too 0.$$ 
Every $\delta_{f_j}$ lies in the unit ball of $\lig^*$, so by the Banach--Alaoglu theorem there is some subnet $\delta_{f_{j_k}}$ converging to a mean $m$. We also have $\l(\g^{-1}.m-m\r)\l(\varphi\r)=0$ for all $\varphi\in\lig$ by continuity.\footnote{We are using here the statements we made earlier about nets in general topological spaces, which we haven't proved.}

So $\g^{-1}.m=m$ for all $\g\in\G$, and thus $m$ is an invariant mean.

The other direction will follow from the discussion below about F\o lner
nets.
\end{proof}

The last proposition is a useful tool for proving amenability. We can now demonstrate an example of an amenable group that isn't finite.
\begin{exa}
The group $(\mb Z,+)$ is amenable, because $\fr{1}{2n+1}\cdot1_{\l[-n,n\r]}$ is an asymptotically invariant sequence in $\loz$.
\end{exa}
The proof this sequence is asymptotically invariant relies on the fact large intervals in $\mb Z$ are ``almost invariant sets'', which is formalised by the notion of \textit{F\o lner nets}, to which the next section is devoted.

\section{F\o lner Nets}

\begin{defn}
A net $(F_j)_{j\in J}$ of finite subsets of a group $\G$ is called \emph{\textbf{a F\o lner net}} if $$\fr{\l|\g F_j\triangle F_j\r|}{\l|F_j\r|}\too 0,\quad\forall\g\in\G.$$ A F\o lner net with $\mb N$ as its domain is called a \emph{\textbf{F\o lner sequence}}.
\end{defn}
\begin{exa}
$F_n=\l[-n,n\r]\subseteq \mb Z$ is a F\o lner sequence, because $\fr{\l|\l(k+\l[-n,n\r]\r)\triangle \l[-n,n\r]\r|}{2n+1}\le \fr{2k}{2n+1}$.
\end{exa}

If $F_j$ is a F\o lner net, then $\fr{1}{\l|F_j\r|}1_{F_j}$ is an a.i.n. in $\leg$. So, if $\G$ admits a F\o lner net, it is amenable. This result is clear, but important, so we state it as its own lemma:
\begin{lem}
\label{f:a}
If a group $\G$ admits a F\o lner net, then it's amenable.
\end{lem}

\begin{exa}
The squares $\l[-n,n\r]^2\subseteq \mb Z^2$  form a F\o lner sequence. (This is because the ratio of the boundary over the area approaches zero).
\end{exa}
\begin{cor}
$\mb Z^2$ is amenable.
\end{cor}
\begin{exe}
By definition, a net $\l(F_j\r)$ of finite subsets of a group $\G$ is F\o lner if $\fr{1}{\l|F_j\r|}\l|\g F_j\triangle F_j\r|\too 0$ for all $\g\in\G$. Show this is equivalent to the condition $\fr{1}{\l|F_j\r|}\l|\g F_j\sm F_j\r|\too 0$ for all $\g\in\G$, and that this is equivalent to the condition $\fr{1}{\l|F_j\r|}\l|S\cdot F_j\sm F_j\r|\too 0$ for all finite $S\subseteq \G$. 
\end{exe}
Another way to go about this is to use the ``F\o lner condition".
\begin{defn}
Let $\G$ be a group, $S\subseteq \G$ finite and $\ep>0$. A subset $F$ is called \emph{\textbf{$\l(S,\ep\r)$-F\o lner}} if $$\mr{max} _{s\in S}\fr{|FS\triangle F|}{|F|}<\ep.$$
\end{defn}
\begin{lem}
\label{F-cond}
A group $\G$ admits a F\o lner net if and only if, for every finite $S\subseteq \G$ and $\ep>0$, there exists a finite subset $F\subseteq \G$ that is $\l(S,\ep\r)$-F\o lner.
\end{lem}
\begin{proof}
Take $J$ to be $\mc S\x \mb N$ with the dictionary order, where $\mc S$ is the set of all finite subsets of $\G$ ordered by inclusion, and $F_{S,n}$ to be an $\l(S,\fr{1}{n}\r)$-F\o lner subset.
\end{proof}
In fact, the opposite of Lemma \ref{f:a} holds as well.
\begin{thm}
\label{th1}
A group is amenable if and only if it admits a F\o lner net.
\end{thm}

When proving this theorem, we will prove another theorem, interesting on its own. It is easy to see that, if a group is paradoxical, then it is not amenable. It turns out that the converse is true as well:

\begin{thm}
\label{na:pd}
A group is amenable if and only if it does not admit a paradoxical decomposition.
\end{thm}

We will prove these two theorems in the next section.

\begin{exe}
Prove there are (at least) two different invariant means on $\mb Z$. (In fact, there are uncountably many).
\end{exe}
\begin{hint}
Consider $F_n=\l[2^n,2^n+2^{n-1}\r)$, 
$\phi_n=\l[2^n-2^{n-2},2^n\r)$
as F\o lner sequences, and the respective functionals $m_n=\fr{1}{\l|F_n\r|}1_{F_n}$, $v_n=\fr{1}{\l|\phi_n\r|}1_{\phi_n}$ in $\liz^*$ (where we consider $\loz$ as a subset of $\liz^*$). They both have limit points $m$, $v$ in the weak-* topology which are different (there is a set $A$ such that $m\l(1_A\r)=1$, $v\l(1_A\r)=0$).

\end{hint}
\begin{exe}
\label{a:p}
Show that if $\G$ has a paradoxical decomposition, and it acts freely on a set $X$, then $X$\ is paradoxical with respect to $\G$.
\end{exe}

\section{Back to Amenability}
We have three conditions on a group $\G$:
\begin{enumerate}
\i $\G$ admits a F\o lner net.
\i $\G$ is amenable.
\i $\G$ has no paradoxical decompositions.
\end{enumerate}
We want to prove these three conditions on $\G$ are equivalent. We already know that if $\G$ has a F\o lner net then it's amenable (Lemma~\ref{f:a}), and that if $\G$ is amenable then it has no paradoxical decompositions (since that rules out the possibility of a finitely additive probability measure). All there is left to show is that, if $\G$ has no paradoxical decompositions, then it admits a F\o lner net.
This is, of course, equivalent to proving that, if $\G$ doesn't admit a F\o lner net, then it does admit a paradoxical decomposition.

The proof we present is due to \cite{deuber1995note}. The fact that amenability implies the existence of F\o lner sets was known before, but the proofs were more complicated; the fact that the existence of F\o lner sets implies amenability was first shown there.

For the proof, we will need a variant of the Marriage theorem:
\begin{thm}
Let $G\l(X,Y\r)$ be a bipartite graph, and suppose $\l|N_Y\l(A\r)\r|\ge2\l|A\r|$ for all finite $A\subseteq X$, where $N_Y\l(A\r)$ is the set of all neighbours of elements in $A$. Then there are two matches from $X$ into $Y$ with disjoint images in $Y$.
\end{thm}

\begin{proof}[Proof of Theorems \ref{th1} and \ref{na:pd}]
Let $\G$ be a group that doesn't admit a F\o lner net. 
By Lemma \ref{F-cond}, there is some finite subset $S\ss \G$ and some $\ep>0$ such that there is no finite subset $A\ss \G$ that is $(S,\ep)$-F\o lner. Assume without loss of generality that $S$ contains the identity element. Then $AS\supseteq A$ for every $A\ss \G$, so that $|AS\triangle A|=|AS|-|A|$. Thus, for all finite $A\subseteq \G$, we have $\l|A\cdot S\r|\ge \l(1+\ep \r)\l| A\r|$. Take $n$ such that $\l(1+\ep\r)^n\ge 2$; then $\l|A\cdot S^n\r |\ge 2\l|A\r|$.

We want to use the Marriage theorem. Consider the bipartite graph $G$ which has two copies of $\G$ as its vertices, $\G_1=\l\{\l(\g,1\r)\mi \g\in \G\r\}$ and $\G_2=\l\{\l(\g,2\r)\mi \g \in G\r\}$, and the pairs $\l\{\l\{\l(x,1\r),\l(y,2\r)\r\}\mi y=xs\text{ for }s\in S^n\r\}$ as its edges. 

This is clearly a bipartite graph, which satisfies the conditions for the Marriage theorem. Hence, there are two matches from $\G_1$ to $\G_2$ such that their images in $\G_2$ are disjoint, which means $2\cdot \G\prec\G$ by decomposing $\G$ to $\l|S^n\r|$  pieces (the $i$th piece consists of all the elements in $\G$ sent by the $i$th element in $S^n$).

Thus, $\G$ admits a paradoxical decomposition, as desired.
\end{proof}

\begin{cor}\label{Amen-fg}
A group $\G$ is amenable if and only if every finitely generated  subgroup of $\G$ is amenable.
\end{cor}
\begin{proof}
A group is amenable if and only if for every finite $S\subseteq \G$, $\ep>0$ there exists an $\l(S,\ep\r)$-F\o lner set. Choose $S$ to be the generating subset.
\end{proof}
\begin{exe}
Let $\G$ be a countable group.
\begin{enumerate}
    \item Show that $\G$ is amenable if and only if it admits a F\o lner sequence. 
    \item Show that $\G$ is amenable if and only if it admits an almost invariant sequence.
\end{enumerate}
\end{exe}

\section{Abelian Groups}
Our original aim was to prove $\mb S^1$ is amenable. We prove this by proving all abelian groups are amenable. By Corollary \ref{Amen-fg}, it is enough to prove this for finitely generated abelian groups.

Recall that every finitely generated abelian group is a finite direct sum of cyclic groups, i.e.\ of the form $$G\oplus\mb Z\oplus\cdots\oplus\mb Z=G\oplus \mb Z^k,$$ where $G$ is finite.
\begin{exe}
Construct F\o lner sets in $G\oplus \mb Z^k$, and deduce all finitely generated abelian groups are amenable.
\end{exe}
\begin{cor}
Every abelian group is amenable, and in particular $\mb S^1$ is.
\end{cor}

\section{Cayley Graphs}

\begin{defn}
If $\G$ is a group and $S\subseteq \Gamma$ is a finite generating subset, then the \emph{\textbf{Cayley graph}} $\mr{Cay}\l(\G,S\r)$ is the graph whose vertices are elements of $\G$ and whose edges are $\l\{\l(\g,\g s\r)\mi \g\in\G,s\in S\r\}$.
\end{defn}

In general, if $\Delta$ is a graph and $A\subseteq \Delta$ is finite (identifying $\Delta$ with its set of vertices), we denote $$\bd A=\l\{x\in\Delta\!\sm\!A\mi x\text{ is a neighbour of an element in }A\r\}.$$ 
In the case of a Cayley graph $\mr{Cay}\l(\G,S\r)$, $\bd A$  is exactly $A\!\cdot\! S\!\sm\!A$.

There is a notion of F\o lner nets in graphs as well:
\begin{defn}
A \emph{\textbf{F\o lner net}} in a connected graph $\Delta$ is a net of finite subsets $\l(A_j\r)_{j\in J}$ of $\Delta$ such that $$\fr{\l|\bd A_j\r|}{\l|A_j\r|}\too 0.$$ 

A graph is called \emph{\textbf{amenable}} if it has a F\o lner net.

A \emph{\textbf{F\o lner sequence}} in a graph is a F\o lner net whose domain is $\mb N$.
\end{defn}

\begin{exe}
\label{F:g-g}
Let $\G$ be an amenable group and let $(F_j)_{j\in J}$ be a F\o lner net of $\G$. Show that, for every finite generating subset $S\subseteq \G$, the net $(F_j)$ is also a F\o lner net in $\mr{Cay}\l(\G,S\r)$, when we consider the $F_j$'s as subsets of vertices of this graph. In particular, $\mr{Cay}\l(\G,S\r)$ is an amenable graph.
\end{exe}
\begin{exe}
Let $\G$ be a group and $S\ss \G$ a finite generating subset. Show that, if $\mr{Cay}(\G,S)$ is amenable, then $\G$ is amenable.
\end{exe}
\begin{hint}
Suppose $F_n$ is a F\o lner sequence. We know $$\fr{\l|\sigma F_n\triangle F_n\r|}{\l|F_n\r|}\too 0$$ for all $\sigma\in\Sigma$. Prove this implies $$\fr{\l|\g F_n\triangle F_n\r|}{\l|F_n\r|}\too 0$$ for all $\g\in\G$ by induction on the length of $\g$ (when represented as a multiplication of elements in $\Sigma$). If $\g=\sigma_1\sigma_2$, then $$\fr{\l|\sigma_1\sigma_2 F_n\triangle F_n\r|}{\l|F_n\r|}=\fr{\l|\sigma_2 F_n\triangle \sigma_1^{-1}F_n\r|}{\l|F_n\r|}\le\fr{\l|\sigma_2 F_n\triangle F_n\r|+\l|\sigma_1^{-1} F_n\triangle F_n\r|}{\l|F_n\r|}.$$
\end{hint}

Accordingly, we will at times abuse notations a bit, and confuse between a group and its Cayley graph(s), as well as the amenability of the two.

Another useful notion is that of the Cheeger constant.
\begin{defn}
The \emph{\textbf{Cheeger constant}} of an infinite graph $\Delta$ is $$C\l(\Delta\r)=\inf_{A\subseteq\Delta,\l|A\r|<\infty}\fr{\l|\bd A\r|}{\l|A\r|}.$$
\end{defn}
\begin{rem}
If $\Delta$ is finite, one takes the infimum over the sets $A$ such that $\l|A\r|\le \fr{1}{2}\l|\Delta\r|$. 
\end{rem}
\begin{exe}
\label{f:ch}
Show $C\l(\Delta\r)=0$ if and only if there exists a F\o lner net in $\Delta$.
\end{exe}

\section{Summary}
We've shown a lot of equivalent definitions for amenability:
\begin{enumerate}[(1)]
\i $\G$ admits a finitely additive, left-invariant probability measure.
\i $\G$ admits a left invariant mean on $\lig$\ (i.e., an element $\varphi\in\lig^*$ that is left invariant, positive and of norm $1$).
\i $\G$ admits an almost invariant net in $\leg$ (equivalently, $\G$ admits an almost invariant net in $\lpg$ for $1\le p <\infty$).
\i $\G$ admits a F\o lner net (equivalently, $C\l(\mr{Cay}\l(\G,S\r)\r)=0$ for all/some finite generating $S\subseteq \G$).
\i $\G$ doesn't admit a paradoxical decomposition.
\end{enumerate}

Examples of amenable groups include all finite groups and all abelian groups (e.g., $\mb Z^n$ and $\mb S^1$). Examples of non-amenable groups include all groups admitting non-abelian free subgroups (e.g., $\mr{SO}\l(3\r)$).

\section*{Further Reading}
For a more thorough development of the theory of amenable groups, readers may refer to \cite{paterson1988amenability}, \cite{runde2004lectures} and
\cite{Kate}.

\chapter{Actions of Amenable Groups}

\section{Actions on Compact-Convex Sets}
In this section we introduce two important mathematical concepts: compact-convex sets and affine maps. We explore the notion of affine group actions on such sets, and provide a characterization of amenability in terms of the fixed point property of such actions. Specifically, we prove that a group is amenable if and only if it has a fixed point under any continuous affine action on a compact-convex set.

Recall that a subset $C$ of a vector space $V$ is \textit{convex} if, for every $x,y\in C$ and $t\in [0,1]$, necessarily $x+t(y-x)$ belongs to $C$ as well.
\begin{defn}
A \emph{\textbf{compact-convex set}} is a compact convex subset of some real topological vector space in which the continuous functionals separate points (i.e., for every $x,y$, there is some continuous functional $f$ such that $f\l(x\r)\ne f\l(y\r)$).
\end{defn}
The requirement that continuous functionals separate points might seem a bit artificial, but in fact `most' topological real vector spaces satisfy it. 
\begin{thm}[Hahn--Banach]
If $V$ is a topological real vector space that is locally convex (i.e., for every $p\in V$ and a neighbourhood $U\subseteq V$ of $p$, there is another neighbourhood $W\subseteq U$ of $p$ that is convex), then its continuous functionals separate points.
\end{thm}
In particular, Banach spaces with the topology induced by the norm, the weak topology, and (if it's a dual space) the weak-* topology all satisfy the condition that continuous functionals separate points.
\begin{defn}
An \emph{\textbf{affine map}} of a compact-convex set $K$ is a continuous map $f:K\to K$ such that $$f\l(\alpha x+\l(1-\alpha\r)y\r)=\alpha f\l(x\r)+\l(1-\alpha\r)f\l(y\r)$$ for all $x,y\in K$, $0\le\alpha\le1$.
\end{defn}
\begin{exe}
Show this is equivalent to the condition $$f\l(\tfrac{1}{2}\l(x+y\r)\r)=\tfrac{1}{2}\l(f\l(x\r)+f\l(y\r)\r).$$
\end{exe}
\begin{hint}
Use induction.
\end{hint}
\begin{rem}
    In fact, one can show affine maps must also satisfy the same condition for the average of $n$ weighted points (for every $n$).
\end{rem}
\begin{thm}
\label{am:af}
A group $\G$ is amenable if and only if for every action of $\G$ on a compact-convex set $K$ by affine, continuous maps, there is a global fixed point (i.e., a point $k\in K$ such that $\g\l(k\r)=k$ for every $\g\in\G$).
\end{thm}
\begin{proof}
Let an amenable group $\G$ act on a compact-convex set $A$ by continuous affine maps. We want to prove it admits a global fixed point.

Let $\l(F_n\r)_{n\in J}$ be some F\o lner net in $\G$; then for all $\g\in \G$ we have $$\fr{\l|\g F_n\triangle F_n\r|}{\l|F_n\r|}\longrightarrow0.$$Take some $p\in A$, and define $$p_n=\fr{1}{\l|F_n\r|}\sum_{f\in F_n}f.p\in A.$$
(The point $p_n$ is in $A$ because it's a convex combination of points in $A$).

$A$ is compact, so $\l(p_n\r)$ has some convergent subnet converging to a point $p_0\in A$. For simplicity, denote this subnet again by $\l(p_n\r)$. We claim this $p_0$ is a global fixed point, as desired.

Let $\g\in\G$. We need to prove $\g.p_0=p_0$. In order to do this, we prove $\varphi\l(\g.p_0\r)=\varphi\l(p_0\r)$ for all continuous functionals $\varphi:V\to \mb R$ (where $V$ is the vector space containing $A$); then it follows $\g.p_0=p_0$ because the continuous functionals on $V$ separate points (i.e., if $\g.p_0\ne p_0$ then there exists a continuous functional $\varphi$ such that $\varphi\l(\g.p_0\r)\ne\f\l(p_0\r)$, a contradiction).

Let $\varphi: V\to \mb R$ be a continuous functional. Note that it's bounded, since $A$ is compact. Then
\begin{align*}
\f\l(\g.p_0\r)&=\f\l(\lim_{n\to \infty}\g.p_n\r)=\lim_{n\to \infty} \f\l(\g.p_n\r)=\lim_{n\to \infty}\l(\fr{1}{\l|F_n\r|}\sum_{f\in F_n}\f\l(\g.f.p\r)\r)\\
&=\lim_{n\to\infty}\fr{1}{\l|F_n\r|}\sum_{f\in \g.F_n}\f\l(f.p\r)
\end{align*}
\begin{align*}
&=\lim_{n\to\infty}\l(\fr{1}{\l|F_n\r|}\sum_{f\in \g F_n\cap F_n}\f\l(f.p\r)+\underset{\too0}{\underbrace{\fr{1}{\l|F_n\r|}\sum_{f\in\g F_n\sm F_n}\f\l(f.p\r)}}\r)\\
&=\lim_{n\to\infty}\l(\fr{1}{\l|F_n\r|}\sum_{f\in \g F_n\cap F_n}\f\l(f.p\r)\r)\\
\end{align*}
\begin{align*}
&=\lim_{n\to\infty}\l(\fr{1}{\l|F_n\r|}\sum_{f\in \g F_n\cap F_n}\f\l(f.p\r)\r)+\underset{=0}{\underbrace{\lim_{n\to\infty}\l(\fr{1}{\l|F_n\r|}\sum_{f\in F_n\sm \g F_n}\f\l(f.p\r)\r)}}\\
&=\lim_{n\to\infty}\l(\fr{1}{\l|F_n\r|}\sum_{f\in F_n}\f\l(f.p\r)\r)=\f\l(\lim_{n\to\infty}\fr{1}{\l|F_n\r|}\sum_{f\in F_n}\l(f.p\r)\r)=\f\l(p_0\r),
\end{align*}
as desired.
\paragraph{}
Conversely, assume that whenever a group $\G$ acts by continuous affine maps on convex-compact sets, it admits a global fixed point. We prove $\G$ is amenable.

Choose $A$ to be the space of means on $\G$, which is contained in $\lig^*$ with the weak-* topology, in which the continuous functionals separate points by the Hahn--Banach theorem. We proved (using the Banach--Alaoglu theorem) $A$ is compact, and clearly it's convex.

Recall that $\G$\ acts on $A$; for all $\g\in\G, m\in A, h\in \lig$, 
$$\l(\g.m\r)\l(h\r)=m\l(\g^{-1}h\r),$$ where for all $g\in G$ $$\l(\g^{-1}h\r)\l(g\r)=h\l(\g g\r).$$
For each $\g\in\G$, $\g:A\to A$ is affine and continuous, so by our assumption it has a global fixed point. This fixed point is an invariant mean, so $\G$ is amenable.
\end{proof}

\begin{rem}
If $\G$ is countable, the space of means on $\G$ is second-countable. Therefore, if a countable group $\G$ admits a global fixed point whenever it acts on a \textit{second-countable} compact-convex subset (by affine, continuous transformations), this is already enough to deduce it is amenable. 
\end{rem}
\section{Actions on Compact Spaces}
The purpose of this section is to prove the following theorem:
\begin{thm}
\label{am:ch}
A group $\G$ is amenable if and only if whenever it acts continuously on a compact Hausdorff space $K$, $K$ admits a $\G$-invariant, Borel probability measure.\footnote{A Borel measure is a measure defined on the Borel $\s$-algebra, which is by definition the smallest $\s$-algebra containing all the open subsets.}
\end{thm}
\begin{rem}
If $X$ is a topological space and $\g:X\to X$ is continuous for every $\g\in\G$, then it's actually a homeomorphism, since $\g^{-1}:K\to K$ is an inverse map which is also continuous.
\end{rem}

For the proof, we need a theorem from functional analysis, which identifies two linear spaces related to compact Hausdorff space $K$. Consider $C\l(K\r)$, the space of continuous functions $f:K\to \mb R$ with the sup norm. This is a Banach space, so it has a dual space $C\l(K\r)^*$; this is the first space. The second space is $M\l(K\r)$, the space of all bounded signed regular Borel measures, where
\begin{enumerate}

\i A bounded signed measure is a generalisation of the concept of finite measures, that is exactly the same except its value can be any real number. It can't be $\infty$ or $-\infty$ (it is `bounded'), and therefore there is no issue of convergence in the requirement of $\s$-additivity.
\i A regular Borel measure is a Borel measure that satisfies 
\begin{align*}m\l(A\r)&=\sup\l\{m\l(C\r)\mi C\subseteq A\text{ compact}\r\}\\
m\l(B\r)&=\inf\l\{m\l(U\r)\mi U\supseteq B\text{ open}\r\}\end{align*}
for all open subsets $A$ and Borel subsets $B$. A signed Borel measure is regular if it's a difference of two positive regular Borel measures.
\end{enumerate}
There is a natural map from $M(K)$ to $C(K)^*$, sending a measure $\mu$ to the functional defined by $f\mapsto \int f\mr{d}\mu$. The theorem we need is:
\begin{thm}[The Riesz Representation theorem]
Let $K$ be a compact Hausdorff space. The map $M(K)\to C(K)^*$ defined by $f\mapsto 
\int f\mr d\mu$ is an isomorphism.
\end{thm}
Now we can prove the theorem.
\begin{proof}[Proof of Theorem \ref{am:ch}]
Assume $\G$ is amenable. Let $K$ be a compact Hausdorff space on which $\G$ acts continuously. We will prove that it admits a $\G$-invariant, Borel probability measure. 

Consider the $\G$-action on $\mr{Prob}\l(K\r)$, the set of (positive) regular Borel probability measures on $K$. $\mr{Prob}\l(K\r)$ is a compact-convex subset of $M\l(K\r)=C\l(K\r)^*$ with the weak-* topology: it is compact because it is a closed subset of the unit ball (which is compact by the Banach--Alaoglu theorem), and it is convex because a convex combination of positive probability measures is also a positive probability measure. Therefore, by the previous section, the action admits a global fixed point. This fixed point is a $\G$-invariant Borel probability measure, as desired.

\paragraph{}Conversely, assume that whenever $\G$  acts continuously on a compact Hausdorff space $K$, $K$ admits a $\G$-invariant, Borel probability measure. We prove $\G$ is amenable, by proving it has a global fixed point whenever it acts by continuous affine maps on a compact-convex set (using the previous subsection).

Thus, let $\G$ act by affine, continuous maps on some Hausdorff compact-convex set. By our assumption, there is a Borel probability measure $m$ on $A$ which is $\G$-invariant. We define the notion of the \textit{barycentre} of a measure on a compact-convex set, and then we prove this barycentre is a point fixed by $\G$.

\begin{defn}
        Let $m$ be a (Borel) probability measure on a compact-convex set $A$. A barycentre of $m$ is a point $C\l(m\r)\in A$ such that $$\f\l(C\l(m\r)\r)=\int_{A}\f\; \mr dm$$for all continuous functionals $\f$ on $V$, where $V$ is the vector space containing $A$.
\end{defn}

\begin{lem}
        Let $m$ be a probability measure on a compact-convex set $A$. Then it has a unique barycentre $C\l(m\r)$.
\end{lem}
\begin{proof}
        Clearly, if such a point exists it is unique, because if $x_1,x_2\in A$ are barycentres then $\f\l(x_1\r)=\f\l(x_2\r)$ for all continuous functionals $\f$ and hence $x_1=x_2$.\footnote{Recall that by definition the continuous functionals on $V$ separate points.} We prove such a point does indeed exist.
        
For any finite set $F$ of continuous functionals on $V$, define 
        $$C\l(F\r)=\l\{c\mi\f\l(c\r)=\int_{A}\!\f\;\mr dm\;\forall\f\in F\r\}.$$
        $C\l(F\r)$ is closed because of continuity. We claim it's nonempty. Define the map
$$G:A\to \mb R^n\quad x\overset{G}{\mapsto}\l(\f_1\l(x\r),\dots,\f_n\l(x\r)\r)$$
for $F=\l\{\f_1,\dots,\f_n\r\}$, and push the measure on $A$ to $\mb R^n$ with it, using
\begin{defn}
Let $\l(X_1,\S_1,\mu\r)$ be a measure space, $\l(X_2,\S_2\r)$ be a measurable space (without a measure), and $T:X_1\to X_2$ some measurable function. The \emph{\textbf{pushforward}} measure $T_*\mu$ of $\mu$  is defined by setting $$T_*\mu\l(A\r)=\mu\l(T^{-1}A\r)$$
for every measurable set $A\ss \Sigma_2$.
\end{defn}
Now we know the barycentre of $G_*m$ exists, because it's the integral
$$\int_{G\l(A\r)}\vec v\;\mr{d}G_*m\l(\vec v\r),$$

and its preimage is in $C\l(F\r)$, and hence $C\l(F\r)$ is nonempty.
        
        $A$ is compact and the collection $\l\{C\l(F\r)\r\}$ satisfies the finite intersection property, so $\bigcap C\l(F\r)$ is also nonempty, and its elements are barycentres by definition. 
\end{proof}

\begin{rem}
The barycentre may be viewed as the integral of the identity function $\mr{id}$ over $A$, $$\int_A\!\mr{id}\;\mr dm,$$where we define integration of functions $f:A\to W$ ($W$ being a real topological vector space in which the continuous functionals separate points) by setting $$\int_A\!f\;\mr dm$$ to be the unique point $c\in W$ satisfying $$\f\l(c\r)=\int_A\!\f\ci f\;\mr dm,$$for all continuous functionals $\f$ on $W$. The existence and uniqueness of such a point can be proved by a similar argument.
\end{rem}

Our proof is almost done. We only need two more things:
\begin{exe}
Let $A\subseteq V$ be a compact-convex set with a probability measure $m$, and  $W$ some other real topological vector space where the continuous functionals separate points.

Prove that for every linear transformation $T:V\to W,$ $$T\l(C\l(m\r)\r)=T\l(\int_A\!\mr{id}\;\mr dm\r)=\int_{T(A)}\!\mr{id}\;\mr dT_*m=C\l(T_*m\r).$$
\end{exe}

Now we just need to show any $\g:A\to A$ can be viewed as a linear transformation.
\begin{exe}
Let $A\subseteq V$ be a compact-convex set. Assume 
that $\mr{Sp}\l(A\r)=V$, and let $T:A\to A$ be an affine map.

Show you can embed $V$ in $V\oplus\mb R$ such that $A$ will go to $V+\l(0,1\r)=\{(v,1)|v\in V\}$ and such that there exists a linear map $\t T$ that extends $T$.
\end{exe}

Now we can prove $C\l(m\r)$ is fixed by $\g$, because it follows from the last two exercises that 
$$\g \l(C\l(m\r)\r)=\int_{\g(A)}\!\mr{id}\;\mr d\g_*m=\int_A\!\mr{id}\;\mr dm=C\l(m\r),$$
because $m$ is $\G$-invariant.
\end{proof}

\section{Operations Preserving Amenability}

So far, we've only proved finite groups and abelian groups are amenable. In this section, we prove that a few operations on groups `preserve' amenability.

\begin{thm}
\label{ele}
        Let $\G$ be an amenable group. Then:
        \begin{enumerate}
                \i \label{7:1}All subgroups of $\G$ are amenable.
                \i \label{7:2}All quotients of $\G$ are amenable. (Equivalently, homomorphic images of $\G$ are amenable).
        \suspend{enumerate}
        Moreover:
        \resume{enumerate}
                \i \label{7:3} All extensions of amenable groups by amenable groups are amenable. That is, if $N \nrml G$ is an amenable normal subgroup of some group $G$, and $G/N$ is also amenable, then $G$ is amenable.
                \i \label{7:4}A direct limit of amenable groups is amenable. In particular, if $\G$ is some group and $\l(\G_n\r)_{n\in J}$ is a net of amenable subgroups of $\G$ such that $n<m\then \G_n\subseteq\G_m$, and $\G=\bigcup_{n\in J}\G_n$, then $\G$ is amenable.\footnote{If you don't know what a direct limit is in general, don't worry about it: we will only use (and prove) this special case.}
        \end{enumerate}
\end{thm}

\begin{exe}
        Prove  Theorem \ref{ele} (\ref{7:1}). That is, show that if $\G$ is amenable and $H\le \G$, then $H$ is amenable.
\end{exe}
\begin{hint}
        Either prove this by contradiction, showing that if $H$ has a paradoxical decomposition then so does $\G$,
        or prove directly that $H$ has an invariant mean (since $\G$ does). 
\end{hint}
\begin{proof}[Proof of Theorem \ref{ele} (\ref{7:2})]
        One way to prove this is using the fixed point property on convex-compact spaces. Suppose $\G$ is amenable and that $N\nrml \G$, and let $\G/N$ act on some convex-compact space $A$ (by affine transformations). An action of $\G/N$ on a space always induces an action of $\G$ on that space  (by letting each element act as its equivalence class acts).
        Since $\G$ is amenable, this induced action (which is still by affine transformations) admits a fixed point, and this point is clearly fixed by the action of $\G/N$ as well. Thus, $\G/N$ is amenable as well.
\end{proof}

\begin{proof}[Proof of Theorem \ref{ele} (\ref{7:3})]
        Let $N$\ be a normal subgroup of $\G$ such that $N,\G/N$ are both amenable. Let $A$ be some compact-convex set on which $\G$ acts affinely and continuously; we want to prove $\G$\ admits a global fixed point (and hence, it's amenable as well).

        Define $\mr{Fix}\l(N\r)\subseteq A$ to be the points fixed by the elements in $N$. $N$ is amenable, so $\mr{Fix}\l(N\r)\ne \varnothing$. It's convex because the action is affine, and it's closed because it's continuous. Therefore, $\mr{Fix}(N)$ is moreover compact (because it lives inside a compact space), hence a compact-convex set.  It's $\G$-invariant: let $\g\in\G$, $p\in\mr{Fix}\l(N\r)$ and $n\in N$; then $$n.\g p=\g\g^{-1}n\g p=\g\l(\g^{-1}n\g\r)p=\g p,$$ because $\g^{-1}n\g\in N$, as $N$ is normal. 

        Thus, $\mr{Fix}\l(N\r)$ is a compact-convex set on which $\G$ acts affinely and continuously. Clearly, $N$ acts trivially on $\mr{Fix}\l(N\r)$, so every element in a given coset $\g N$ acts in the same way on $\mr{Fix}\l(N\r)$, and hence $\G/N$ acts on $\mr{Fix}\l(N\r)$ affinely and continuously as well.

        We assumed $\G/N$ is amenable, so there is a point $p\in\mr{Fix}\l(N\r)$ which is fixed by all of $\G/N$, and therefore by all of $\G$. Thus, $p\in A$ is a point fixed by all of $\G$, as desired.
\end{proof}
\begin{cor}
        Solvable groups are amenable.
\end{cor}
\begin{proof}
        Let $\G$ be a solvable group. Then for some $n$,
        $$\l\{1\r\}=\G^{\l(n\r)}\nrml\cdots\nrml\G^{\pr\pr}\nrml\G^\pr\nrml\G,$$
        and since $\G/\G^\pr$ is always amenable, we get by inductively using Theorem \ref{ele} (\ref{7:3}) that $\G$ is amenable.
\end{proof}

\begin{exe}
        Show that the class of solvable groups is the class of groups obtained by repeatedly using operation (\ref{7:3}) of Theorem \ref{ele} on abelian groups.
\end{exe}

\begin{proof}[Proof of Theorem \ref{ele} (\ref{7:4})]
        Suppose $\l(\G_n\r)_{n\in J}$ is a net of amenable groups such that $n<m$ implies $\G_n\subseteq\G_m$, and that $\G=\bigcup_{n\in J}\G_n$. We can prove $\G$ is amenable in more than one way.

        \paragraph{Proof using the F\o lner condition}
        For every finite $F\subseteq \G$ and $\ep>0$, take $n\in J$ such that $F\subseteq \G_n$ (which exists by induction on the size of $F$). $\G_n$ is amenable, so there is an $\l(F,\ep\r)$-F\o lner set in $\G_n$, which is clearly  $\l(F,\ep\r)$-F\o lner in $\G$. Thus, for all finite $F\subseteq \G$ and $\ep>0$ there exists an $\l(F,\ep\r)$-F\o lner set in $\G$, and therefore it's amenable.
        
        \paragraph{Proof using the fixed point property}
        Each $\G_n$ has a set of fixed points, $\mr{Fix}\l(\G_n\r)$, which is nonempty and compact. Moreover, the collection $\{\mr{Fix}(\G_n)\}_{n\in J}$ is ordered by inclusion (i.e., $\mr{Fix}(\G_m)\ss \mr{Fix}(\G_n)$ whenever $m\ge n$), and $J$ is directed, so this collection satisfies the finite intersection property. Since $A$ is compact, this means $$\bigcap_{n\in J}\mr{Fix}\l(\G_n\r)$$ is nonempty, and hence $\G$ admits a fixed point as well, and therefore it is amenable.
\end{proof}

\subsection{Elementary Amenable Groups}

\begin{defn}
        A group is called \emph{\textbf{elementary amenable}} if one of the following conditions holds:
        \begin{enumerate}
                \i It's finite. 
                \i It's abelian.
                \i It's obtained from some other elementary amenable groups using the operations described in Theorem \ref{ele}.
        \end{enumerate}
\end{defn}
 Theorem \ref{ele} implies that elementary amenable groups are indeed amenable.

\chapter{Some Applications of Amenability}
\section{Banach Limits}
Banach wanted to define a notion of a generalised limit, $\mr{LIM}$, that will be defined for all bounded sequences $\l(x_n\r)_{n\in \mb N}$ in $\mb R$, and that will satisfy the following:
\begin{enumerate}
\i If $\lim x_n$ exists, then $\lim x_n=\L x_n$.
\i It is shift invariant, i.e.\ $\L_{n\to \infty} x_n=\L_{n\to \infty}x_{n+k}$ for every $k$.
\i It is linear, i.e.\ $\L\l(ax_n+by_n\r)=a\L x_n+b\L y_n$.
\i It satisfies $\l|\L x_n\r|\le \sup \l|x_n\r|$,
\end{enumerate}
\begin{thm}
Such $\L$ exists.
\end{thm}
$\L$ is actually a functional of norm $1$ on $\lin$, the space of bounded sequences. This space, like $\lig$, is a Banach space, so its space of functionals $\lin^*$ is also a Banach space. The only difference is that $\mb N$ is only a semigroup, so for the proof we'll need the notion of amenability in semigroups. 

Recall we can identify $\lon$ as a subspace of $\lin^*$ by identifying the elements of $\lon $ with their evaluation maps in $\lin$. That is, if $x\in \lon$, we can think of it as a functional on $\lin$, by setting $$x\l(y\r)=\sum_{i=1}^\infty x\l(i\r)\cdot y\l(i\r)$$ for all $y\in \lin$.

With this in mind, the proof is very easy.
\begin{proof}
Consider the functions $\fr{1}{n}1_{\l[1,n\r]}\in\lon\subseteq\lin^*$. Clearly they all lie in the unit ball of $\lin^*$, so by the Banach--Alaoglu theorem they have some limit point, $\L\in \lin^*$. This $\L$ satisfies our conditions.
\end{proof}
\begin{exe}
Check this indeed satisfies all our requirements. For shift invariance, note that the sets $\l[1,n\r]$ form a F\o lner sequence and hence the functions $\fr{1}{n}1_{\l[1,n\r]}$ form an almost invariant sequence (this is where the notion of amenability in semigroups is used).
\end{exe}

\section{Von Neumann Ergodic Theorem}
Let $\l(X,B,m\r)$ be a probability measure space. We want to define the notion of \textit{ergodicity}.
\begin{rem}
\label{m0}
When dealing with measure spaces, we will constantly ignore sets of measure zero. This means that assertions that fail in sets of measure zero may still be regarded as true; for example, a function $T:X\to X$ is called invertible if there is some $ U:X\to X$ such that $T\l(U\l(x\r)\r)=U\l(T\l(x\r)\r)=x$ \emph{almost} everywhere. \end{rem}
\begin{defn}
Let $\l(X,\Sigma,m\r)$ be a probability space. A measurable function $T:X\to X$ is \emph{\textbf{measure-preserving}} if $m\!\l(T^{-1}\l(A\r)\r)=m\!\l(A\r)$ for all $A\in \Sigma$. An \emph{\textbf{automorphism}} on a probability space is an invertible, measure-preserving function $T:X\to X$.
\end{defn}

If $T$ is an automorphism, its inverse is also an automorphism, and we have
\begin{align*}
m\!\l(TA\r)&=m\l(A\r),\quad\forall A\in B\\
m\!\l(X\sm TX\r)&=0.
\end{align*}
\begin{defn}
Let $\l(X,\Sigma ,m\r)$ be a probability space, $T:X\to X$ a measure preserving map. $T$ is \emph{\textbf{ergodic}} if $T^{-1}\l(A\r)=A$ implies $m\!\l(A\r)=0$ or $m\!\l(A\r)=1$.
\end{defn}
\begin{rem}
This is equivalent to saying that, if $T\l(A\r)=A$ almost everywhere, then  $m\l(A\r)\in\l\{0,1\r\}$ (in accordance with Remark \ref{m0}).
\end{rem}

Consider $$L_2\l(X,m\r)=\l\{f:X\to \mb R\text{ measurable}\mi\l(\int_X\l|f\r|^2\mr dm\r)^\fr{1}{2}<\infty\r\}.$$ This is a Hilbert space with $$\l<f,g\r>=\int_X f\cdot g\;\mr dm.$$
\begin{thm}[Von Neumann ergodic theorem]
If $T:X\to X$ is ergodic, then, for all $f\in L_2\l(X,m\r)$, the functions $$g_n=\fr{1}{n}\sum_{i=1}^nf\circ T^i\in L_2\l(X,m\r)$$ converge (in $L_2$) to the constant function $$g\equiv\int_X f\;\mr dm.$$
\end{thm}
\begin{proof}

We want to use the language of operators.

Define the operator $U_T:L_2\l(X,m\r)\to L_2\l(X,m\r)$ by setting $U_T\l(f\r)=f\circ T$. Note that $U_T$ is a unitary operator.

We want to prove the sequence $V_i\coloneqq\fr{1}{n}\sum_{i=1}^nU_T^i$ converges pointwise to the projection $\pi$ onto $\mc C\coloneqq\mr{Sp}\l(1\r)$ (the space of constant functions, which is the span of the function that is constantly $1$). This is because $$1\cdot \l<f,1\r> =\int_Xf\;\mr dm$$ for all $f\in L_2\l(X,m\r)$.

In order to do this, we consider $L_2$ as a direct sum of $\mc C$,  the space of constant functions, with $\mc C^\perp$,  the space of functions with whose integral is $0$. Then, we show $V_i$ converges pointwise in each of these subspaces, and deduce by linearity it converges in all of $L_2$.

\begin{exe}
Prove rigorously $U_T\l(f\r)=f$ if and only if $f$ is constant. (Recall we assume $T$ is ergodic).
\end{exe}

Clearly, this implies $V_i\too \pi$ in the subspace $\mc C$.

More interestingly, consider the operator $U_T-1$, $$\l(U_T-1\r)\l(f\r)\l(x\r)=f\l(T\l(x\r)\r)-f(x).$$ The exercise shows that $\ker \l(U_T-1\r)=\mc C$. 

It is true that for any operator $A$, the equality $\ker A=(\mr{Im}A^*)^\perp$ holds, and that in any Hilbert space $A^{\perp\perp}=\bar{A}$. Thus, $$\mc C^\perp=\l(\ker\l(U_T-1\r)\r)^\perp=\l(\mr{Im}\l(U_T-1\r)^*\r)^{\perp\perp}=\overline{\mr{Im}\l(U_T-1\r)^*}.$$ But $\l(U_T-1\r)^*=U_T^*-1^*=U_T^{-1}-1$, so $$\mc C^\perp=\overline{\mr{Im}\l(U_T^{-1}-1\r)}.$$ Recall we want to prove $V_i\to \pi$ in $\mc C^\perp$ (pointwise). First we do it for $f\in\mr{Im}\l(U_T^{-1}-1\r)$, then for general $f\in\overline{\mr{Im}\l(U_T^{-1}-1\r)}$.

Let $f\in\mr{Im}\l(U_T^{-1}-1\r)$. Then there exists some $g\in L_2$ such that $f=U_T^{-1}\l(g\r)-g$, and then $$V_i\l(f\r)=\fr{1}{n}\sum_{i=1}^nU_T^i\l(f\r)=\fr{1}{n}\sum_{i=1}^nU_T^i\l(U_T^{-1}\l(g\r)-g\r)=\fr{1}{n}\l(g-U_T^n\l(g\r)\r),$$ and this converges to zero (which is $\pi\l(f\r)=1\cdot \l<f,1\r>$) because $U_T^n$ is unitary (and hence $\l\|g-U_T^n\l(g\r)\r\|\le2\cdot\l\|g\r\|$).

Given $f\in\overline{\mr{Im}\l(U_T^{-1}-1\r)}$, take $g\in L_2$ such that $$\l\|f-\l(U_T^{-1}-1\r)g\r\|<\fr{\ep}{2}$$ and take $N$ such that $n>N$ implies 
$$\l\|\fr{1}{n}\sum_{i=1}^nU_T^i\l(U_T^{-1}(g)-g\r)\r\|< \fr{\ep}{2}$$ 
and deduce that $n>N$ implies $$\l\|\fr{1}{n}\sum_{i=1}^nU_T^i\l(f\r)\r\|<\ep.$$
\end{proof}

You might wonder where we used amenability here. Once again, we used the amenability of the semigroup $\mb N$, by considering its action on $X$ by measure-preserving transformations; each $n$ acts on $X$ by $T^n$. We can replace $\mb N$ with any amenable group or semigroup, and then the theorem becomes
\begin{thm}
If $\G$ acts ergodically on a probability space $\l(X,\Sigma,m\r)$ by measure-preserving transformations, and $\l(F_j\r)_{j\in J}$  is a F\o lner net (or sequence), then $$\fr{1}{\l|F_j\r|}\sum_{\g\in F_j}L_{\g}f\longrightarrow\int_X f\;\mr dm,$$ where $L_\g f\l(x\r)=f\l(\g^{-1}x\r)$.
\end{thm}
\begin{rem}
    If $\G$ does not act ergodically, we have the following convergence:
    $$\fr{1}{\l|F_j\r|}\sum_{\g\in F_j}L_{\g}f\longrightarrow\int_X f\;\mr dm,$$ where $\pi:L_2\to L_2^{\G}$ is the projection onto the space of invariant functions.
\end{rem}

\begin{rem}
The Birkhoff pointwise ergodic theorem \cite{Birkhoff}
says that, if $T$ is an ergodic measure preserving transformation of a probability space $(X,m)$,
then for every $f\in L_1(X,m)$ the ergodic average (or ``time average") converges to the integral (or ``space average") almost everywhere, i.e. 
$$
 \fr{1}{n}\sum_{i=1}^nf( T^{-i}x)\too \int f \mr{d}m
$$
for all $x$ outside a set of measure $0$.

The proof of this theorem is more involved than that of the Von Neumann ergodic theorem and does not extend to every F\o lner sequence in every amenable group. 
However, E. Lindenstrauss \cite{Lindenstrauss} proved that every amenable group admits a F\o lner sequence with respect to which the pointwise ergodic theorem holds.
 
\end{rem}

\chapter{Growth of Groups}

\section{Definition}

Given two increasing, non-zero functions $f,g:\mb N\to\mb R$, we denote \emph{\textbf{$f\prec g$}} if there is some $c>0$ such that $f\l(n\r)<cg\l(cn\r)$ for every $n\in\mb N$. We denote \emph{\textbf{$f\sim g$}} if both $f\prec g$ and $g\prec f$. This is an equivalence relation, and two functions are said to have the same \emph{\textbf{growth type}} if $f\sim g$.
\begin{defn}
Given a finitely generated group $\G$, $\G=\l<\Sigma\r>$ (where $\Sigma$ is finite, symmetric and contains the identity element $1$), the \emph{\textbf{growth type}} of $\G$ is the growth type of the function $$f_\Sigma\l(n\r)\coloneqq\l|\Sigma^n\r|,$$
where $\Sigma^n=\{\g_1\cdots\g_n|\g_i\in \Sigma\}$, the ball around the origin of radius $n$ in $\mr{Cay}\l(\G,\S\r)$.
\end{defn}
\begin{exe}[Important]
The growth type of $f_\Sigma$ is independent of $\Sigma$, i.e.\ the growth of $\G$ is well defined.
\end{exe}
A finitely generated group $\G=\l<\Sigma\r>$, $\l|\Sigma\r|<\infty$, has \emph{\textbf{polynomial growth}} if there is some $k\in\mb N$ such that $\l|\S^n\r|\le n^k$ for every $n\in \mb N$. $\G$ has \emph{\textbf{exponential growth}} if there is some $c>1$ such that $\l|\S^n\r|\ge c^n$ for every $n\in\mb N$. $\G$ has \emph{\textbf{sub-exponential growth}} if, for every $c>1$, there is some $N\in\mb N$ such that $n>N$ implies $\l|\S^n\r|< c^n$.

Observe that $\S^{n+m}=\S^n\S^m$ and hence $\l|\S^{n+m}\r|\le\l|\S^n\r|\!\cdot\! \l|\S^m\r|$, i.e.\ $f_\S$ is sub-multiplicative. This is equivalent to saying $\log\circ f_\S$ is sub-additive.
\begin{exe}
Assume $g:\mb N\to \l[0,\infty\r)$ is sub-additive. Show the limit $\lim\fr{g\l(n\r)}{n}$ exists and equals $\inf\fr{g\l(n\r)}{n}$.
\end{exe}
This means $\lim \fr{1}{n}\log\l|\S^n\r|$ exists, and hence $\lim\l|\S^n\r|^{1/n}=e^{\lim\fr{1}{n}\log\l|\S^n\r|}$ exists and equals $d\coloneqq\inf\l|\S^n\r|^{1/n}$, and we can define this number to be the \emph{\textbf{exponential growth rate of $\S$}}. Thus, $\G$ has exponential growth if and only if $d>1$ (just pick $c=d$) and it has sub-exponential growth if and only if $d=1$ (because, for every $c>1$, you can pick $\ep=\fr{c-d}{2}$, and then there is some $N\in\mb N$ such that $\l|\S^n\r|^{1/n}<d+\ep<c$ for any $n\ge N$). 
\section{Growth of Nilpotent and Solvable Groups}

Let $\G$ be a group. We define the \emph{\textbf{lower central series}} $\l(\G_n\r)_{n\in\mb N}$  by $$\G_0=\G,\quad\G_{n+1}=[\G,\G_n],$$ where $\l[a,b\r]=aba^{-1}b^{-1}$ and $\l[M,N\r]=\l<\l\{\l[m,n\r]\mi m\in M,n \in N\r\}\r>$.

Recall that a subgroup $N$ of $\G$ is normal if $g^{-1}Ng=N$ for all $g\in\G$. We say a subgroup $N$ of $\G$ is \emph{\textbf{characteristic}}  if $\varphi\l(N\r)=N$ for all automorphisms $\varphi$ of $\G$. Clearly, this is a stronger requirement, since $x\mapsto g^{-1}xg$ is an automorphism for all $g\in \G$.
\begin{exe}
Let $\G$ be a group, $N,M\le \G$. Show the following:
\begin{enumerate}
\i If $N\nrml \G$ then $\l[N,M\r]\le N$.
\i If $N,M\nrml \G$ then $\l[N,M\r]\nrml \G$.
\end{enumerate}
\end{exe}

By a very similar solution to the last exercise, one can show that if $N,M\le \G$ are characteristic, then so is $\l[N,M\r]$, and hence $\G_n$ are characteristic subgroups of $\G$.

\begin{defn}
A group $\G$ is \emph{\textbf{nilpotent}} if $\G_n=\l\{1\r\}$ for some $n\in\mb N$.
\end{defn}

\begin{exa}
The group of upper triangle matrices with ones in the diagonal, $$\G=\l\{\begin{pmatrix}1 & &*\\ &\ddots & \\ 0 & &1 \end{pmatrix}\r\},$$ is nilpotent.
\end{exa}

In group theory, given some property $P$, it's common to say a group is \emph{\textbf{virtually $P$}} if it has a subgroup with finite index that is $P$. For example, a group $\G$ is virtually nilpotent if there exists some finite index subgroup $H\le \G$ that is nilpotent.
\begin{thm}[Gromov]
Let $\G$ be a finitely generated group. Then it has polynomial growth if and only if it's virtually nilpotent.
\end{thm}
It is pretty straightforward to show that a virtually nilpotent group has polynomial growth. The other direction is much deeper; we recommend reading the original proof \cite{Gromov-poly} as well as later proofs \cite{kleiner2010new,Ozawa}.

Another series of characteristic subgroups of $\G$ is the \emph{\textbf{derived series}}, defined by $$\G^{\l(0\r)}=\G,\quad\G^{\l(n+1\r)}=\G^{\l(n\r)\pr}=[\G^{\l(n\r)},\G^{\l(n\r)}].$$
Recall that a group $\G$ is \emph{\textbf{solvable}} if $\G^{\l(n\r)}=\l\{1\r\}$ for some $n\in\mb N$.
Obviously we have $\G^{\l(n\r)}\subseteq\G_n$, so nilpotent  groups are always solvable. The converse does not hold.
\begin{exa}
The group $$\G=\l\{\begin{pmatrix}* & &*\\ &\ddots &\\ 0 & &* \end{pmatrix}\r\}$$ of upper triangle matrices is a solvable, non-nilpotent group.
\end{exa}

\begin{defn}The \emph{\textbf{free semigroup}} on a set $S$ is the semigroup of finite sequences of elements in $S$, with the operation of concatenation. More formally, $$\G=\bigcup_{n\in\mb N}S^n=\bigcup_{n\in\mb N}\l\{f\mi f:[n]\to S\r\},$$ where $[n]=\l\{1,\dots,n\r\}$, and for all $f_1:[r_1]\to S, f_2:[r_2]\to S\in \G$ we define $f_1\cdot f_2$ by
$$f_1\cdot f_2\l(n\r)=
\begin{cases}
f_1\l(n\r) &n\in [r_1]=\l\{1,\dots,r_1\r\},\\
f_2\l(n-r_1\r) &n\in 
\l\{r_1+1,\dots,r_2\r\}.
\end{cases}$$
A semigroup is said to be \emph{\textbf{free}} if it's isomorphic to some free semigroup. In this case, it has some  generating set $S$ such that all words 
in $S$ are different.
\end{defn}
\begin{rem}
Note that the free semigroup on a set $S$ is abelian if and only if $\l|S\r|=1$.
\end{rem}
\begin{exe}
Show that if a finitely generated group $\G$ contains a non-abelian free subsemigroup, then it has exponential growth.
\end{exe}

\begin{thm}
A solvable group which isn't virtually nilpotent admits a non-abelian free subsemigroup. \end{thm}
We will not prove this theorem either, but an immediate corollary of it is:
\begin{cor}
A finitely generated solvable non-virtually nilpotent group has exponential growth. Hence, a finitely generated solvable group either has polynomial growth or exponential growth.
\end{cor}
\section{Growth of Non-Amenable Groups}
\begin{prop}
\label{ne:a}
If $\G$ has sub-exponential growth then it's amenable, and there is a F\o lner sequence consisting of balls $\S^{n_k}$ around the identity.
\end{prop}
\begin{proof}
Assume the contrary, i.e.\ $\l(\S^n\r)_{n\in\mb N}$ has no F\o lner subsequence, which means there is some $\ep>0$ such that $$\fr{\l|\partial\Sigma ^n\r|}{\l|\S^n\r|}\ge\ep$$ in the graph $\mr{Cay}\l(\G,\S\r)$, for every $n$. But $\partial \S^n=\S^{n+1}\sm\Sigma^n$, so this implies $$\l|\S^{n+1}\r|\ge\l|\S^n\r|\cdot\l(1+\ep\r)$$ (for $n$ large enough), and hence $$\l|\S^n\r|\ge\l(1+\ep\r)^n,$$ and $\G$ has exponential growth, a contradiction.
\end{proof}
\begin{cor}Non-amenable groups have exponential growth.
\end{cor}
\begin{rem}
We have seen that solvable groups have either exponential or polynomial growth, and that non-amenable groups always have exponential growth. There are, however, amenable (necessarily non-solvable) groups with intermediate growth (slower than exponential, but faster than polynomial). The first such example was given by Rostislav Grigorchuk; see \cite{grigorchuk2008groups}, as well as \cite{grigorchuk14burnside,grigorchuk1984degrees} for the original proof.
\end{rem}
\begin{exe}
Let $a\in \mb R$ be a transcendental number. Consider the group of affine transformations of $\mb R$ generated by multiplication by $a\in \mb R$ and addition by $1$. That is,$$\G\coloneqq\l<\begin{pmatrix}a &0\\0 &1\end{pmatrix},\begin{pmatrix}1 &1\\0 &1\end{pmatrix}\r>\curvearrowright\l\{\begin{pmatrix}r\\1\end{pmatrix}\mi r\in\mb R\r\}.$$
The group $\G$ is amenable (because it's solvable). Prove it has exponential growth and deduce the converse of Proposition \ref{ne:a} doesn't hold.

Prove that, in fact, $$\begin{pmatrix}a &0\\0 &1\end{pmatrix},\begin{pmatrix}1 &1\\0 &1\end{pmatrix}$$ generate a free semigroup.
\end{exe}
\begin{rem}
In an amenable group with exponential growth, you can't find a F\o lner sequence of balls around the identity in $\G$.
\end{rem}
\subsection{The Ping Lemma}
Back when we wanted to find free subgroups, we used the ping-pong lemma. When looking for free subsemigroups, there is a similar tool:
\begin{lem}[The Ping Lemma]
Let $\G\acts X$. If there exist elements $a,b\in\G$ and disjoint subsets $A,B\subseteq X$ such that 
\begin{align*}
a(A\cup B)&\subseteq A,\\
b(A\cup B)&\subseteq B,
\end{align*}
Then the subsemigroup generated by $a$ and $b$ is a free semigroup.

\begin{center}
\begin{tikzpicture}

\pa (0,0);
\pa (3,0) \at (A);
\pa (0,0) \at (B);

\dr (A)node {$A$} circle (0.8);
\dr (B)node {$B$} circle (0.8);

\pa[->] (0,0.9) edge [bend left] (3,0.9);
\pa[->] (3,-0.9) edge [bend left] (0,-0.9);

\dr (1.5,1.6) node {$a$};
\dr (1.5,-1.6) node {$b$};
\end{tikzpicture}
\end{center}
\end{lem}
\begin{exe}
Prove the Ping Lemma.
\end{exe}
\begin{exa}
Consider $$\mr{Aff}(\mb R)=\l\{\begin{pmatrix}a &b\\0 &1\end{pmatrix}\mi a,b\in\mb R\r\}\acts\l\{\begin{pmatrix}v\\1\end{pmatrix}\mi v\in\mb R\r\},$$ where $$\begin{pmatrix}a &b\\0 &1\end{pmatrix}\cdot\begin{pmatrix}v\\1\end{pmatrix}=\begin{pmatrix}av+b\\1\end{pmatrix}.$$ Take $g$ to be the map defined by $x\mapsto \fr{1}{3}x$, i.e.\ $$g=\begin{pmatrix}\fr{1}{3} &0\\0 &1\end{pmatrix},$$ and $h$ to be the map defined by $x\mapsto \fr{1}{3}(x-1)+1$, i.e.\ $$h=\begin{pmatrix}\fr{1}{3} &\fr{2}{3}\\0 &1\end{pmatrix}.$$ Then, for $A=\l(-\fr{1}{2},\fr{1}2{}\r)$, $B=\l(\fr{1}{2},\fr{3}{2}\r)$, the Ping Lemma holds.
\end{exa}

\begin{rem}
For some time there was a lack of additional interesting examples of amenable groups besides elementary amenable groups and groups of sub-exponential  growth. However, the groundbreaking paper of Juschenko and Monod  
\cite{juschenko2013cantor} which produced the
first example of a (finitely generated, infinite) simple amenable group revealed a new source for plenty of examples.

\end{rem}

\chapter{Topological Groups}

When a group admits a topology which behaves ``nicely enough'' with
its group structure, it is called a \emph{topological group}. A lot
of the groups we have been discussing admit natural underlying topologies,
making them into topological groups. 

In the first section of this chapter, we define topological groups. 
Throughout the rest of the book, unless otherwise stated, all groups
are topological groups; when we want to consider groups without topology,
we will call them ``abstract groups", or ``discrete groups" (since abstract groups are more 
or less the same as topological groups with the discrete topology).
Moreover, we will always assume (unless otherwise stated) that 
locally compact groups are second-countable.

\section{Definition and Examples}
\begin{defn}
	A \textbf{\emph{topological group}} is a group with a Hausdorff topology
	such that the map
	\begin{align*}
		G\times G & \to G\\
		(x,y) & \mapsto x\cdot y
	\end{align*}
	and the map 
	\begin{align*}
		G & \to G\\
		x & \mapsto x^{-1}
	\end{align*}
	are continuous (where $G\times G$ is endowed with the product topology).
\end{defn}

\begin{rem}
	To be more formal, we should have defined a topological group to be
	a triple $(G,\cdot,\tau)$, such that $(G,\cdot)$ is a group and
	$\tau$ is a topology on $G$, such that multiplication and inversion
	are continuous. But, as a standard abuse of notation, we will think
	of $G$ as being both a group and a topological space simultaneously.
\end{rem}

\begin{rem}
	It was enough for us to require the topology of the group to be $T_{1}$.
	This is because a space is Hausdorff if and only if its diagonal is
	closed, and the diagonal is the preimage of $\{1\}$ under $(x,y)\mapsto x\cdot y^{-1}$,
	which is continuous by Exercise \ref{tg:c}. In fact, it is not hard
	to show that it was even enough to require the topology to be $T_{0}$,
	and that this implies not only that it is Hausdorff, but even $T_{3.5}$.
\end{rem}

\begin{rem}
	\label{tg:c} Given a group $G$ endowed with some topology, show that the map $(x,y)\mapsto x\cdot y^{-1}$ is continuous if and only if both of the maps $(x,y)\mapsto x\cdot y$ and  $x\mapsto x^{-1}$
	are continuous.
\end{rem}

\begin{exa}
	Any group with a discrete topology is a topological group. Some groups
	do not admit any other topologies making them into topological groups. 
\end{exa}

\begin{exa}
	If $G$ is a topological group, then any subgroup of it (with the
	induced topology) is a topological group. If $G_{1},G_{2}$ are topological
	groups, then $G_{1}\times G_{2}$ is also a topological groups (with
	the product topology).
\end{exa}

\begin{exa}
	The group $\mathrm{GL}(n,\mathbb{R})$ of $n\times n$ invertible
	matrices over $\mathbb{R}$ is a topological group (with the subspace
	topology of $\mathbb{R}^{n^{2}}$). 
\end{exa}

\begin{exa}
	The group $\mathrm{SO}(n,\mathbb{R})$ is a topological group (with
	the subspace topology of $\mathbb{R}^{n^{2}}$).
\end{exa}

\begin{exa}
    \label{autT}
	If $T$ is a locally finite graph, then the group $\mathrm{Aut}(T)$
	of its automorphisms has a topology that makes it into a topological
	group. This topology may be defined by the following basis of neighbourhoods
	around the identity element: $\left\{ \mathrm{Fix}(K)\middle|K\subseteq T\text{ is finite}\right\} $,
	where $\mathrm{Fix}(K)$ is the subgroup of elements fixing every
	vertex in the finite set $K$. One says that an automorphism is ``close
	to $1$'' if it fixes (pointwise) a large subset of vertices, and two elements
	are close to each other if they act similarly on a large subset of vertices. This
	group does not admit a structure of a Lie group.
\end{exa}

\begin{exe}
    Let $G$ be a topological group. Show that it is locally compact if
    and only if it admits a compact neighbourhood of the identity element.
\end{exe}
    
\begin{notation}
	From now on, when we say that a topological group is \textbf{\emph{locally compact}},
	we will in fact mean that it is locally compact and second-countable.
	
\end{notation}

\begin{exa}
        The group $\mathrm{Aut}(T)$ (where $T$ is a connected locally
	finite graph) is locally compact. In fact, the basis of identity neighbourhoods
	we gave when defining its topology consists of open, compact subgroups.
\end{exa}

\begin{exa}
	If $k$ is a topological field\footnote{That is, a field with a topology such that addition, multiplication and the inversion of both are continuous.}  (such as $\mathbb{R},\mathbb{C}$
	or $\mathbb{Q}_{p}$), then $\mathrm{GL}_{n}(k)$ is a topological
	group, which is locally compact if and only if $k$ is.
\end{exa}

\section{Lie Groups}

Possibly the most important class of topological groups is the class
of Lie groups. We assume in this section that the reader is at least
vaguely familiar with what a smooth manifold is.
\begin{defn}
	A \textbf{\emph{Lie group}} is a group which is also a smooth manifold,
	such that the maps of multiplication and inversion are smooth.
\end{defn}
\begin{exe}
    Show that the group $\mr{GL}(n,\mb R)$, with the  smooth structure induced from $\mb R^{n^2}$, is a Lie group.
\end{exe}

\begin{thm}[Cartan's theorem]
    If $G$ is a Lie group and $H\leqslant G$ is a closed subgroup (that is, a subgroup which is closed as a subset of $G$), then $H$ is also a Lie group with the induced topology and smooth structure.
\end{thm}
We immediately get that $\mr{SO}(n,\mb R)$ and $\mr{SL}(n,\mb R)$ are Lie groups as well (although it is not hard to show this directly, either). A Lie group $G$ which admits an injective homomorphism $f:G\to \mr{GL}(n,\mb R)$ (for some $n$) is called \textit{linear}. One can show that every linear Lie group $G$ admits an injective homomorphism $f:G\to \mr{GL}(n,\mb R)$ which has a closed image. Not all Lie groups are linear, but every Lie group $G$ admits a homomorphism $f:G\to \mr{GL}(n,\mb R)$ whose kernel is discrete. 

It is not very difficult to show that, if a group is an \emph{analytic
}manifold and the maps of multiplication and inversion are
smooth, then they are necessarily also analytic. For a long time,
people didn't know whether it was enough for them to be continuous.
This was essentially Hilbert's fifth problem. It was resolved in the 1950s,
by Montgomery, Zippin and Gleason:
\begin{thm}[Gleason--Montgomery--Zippin]
	Let $G$ be a topological group which is locally Euclidean (that
	is, admits an identity neighbourhood homeomorphic to $\mathbb{R}^{n}$
	for some $n$). Then $G$ admits a unique structure of an analytic
	manifold for which the maps of multiplication and inversion are analytic.
\end{thm}

\begin{rem}
	Since the smooth structure in this case is unique, we will often abuse
	notations a little, and say that a topological group ``is'' a Lie
	group, rather than say that it admits a structure of a Lie group.
\end{rem}
\begin{rem}
	It is moreover true (and much easier to show) that, if $f:G\to H$
	is a homomorphism between two Lie groups which is continuous, then
	it is necessarily analytic. Therefore, one can think of Lie groups
	as a class of topological groups, which admit some extra structure.
	In the language of category theory, one says that the category of
	Lie groups is a ``full subcategory'' of the category of topological
	groups. 
\end{rem}

The research on Hilbert's fifth problem also led to the following
remarkable results.
\begin{defn}
	A topological group is said to have \emph{no small subgroups }(or
	\emph{NSS}), if it admits an identity neighbourhood which doesn't
	contain any nontrivial subgroups.
\end{defn}

\begin{thm}
	Let $G$ be a topological group. Then $G$ is a Lie group if and only
	if it is locally compact and has NSS.
\end{thm}

\begin{thm}
	Let $G$ be a connected locally compact group. Then every identity
	neighbourhood of $G$ contains a compact normal subgroup $K$ such that $G/K$
	is a Lie group.
\end{thm}

\section{Totally Disconnected Locally Compact Groups}

Another very important class of topological groups is the class of
totally disconnected locally compact (t.d.l.c.) groups. A good motivation
for studying this class of groups is the following:
\begin{notation}
	When $G$ is a topological group, we denote the connected component
	of its identity element by $G^{\circ}$. This is a normal subgroup.
\end{notation}

\begin{prop}
	Let $G$ be a locally compact group. Then $G/G^{\circ}$ is locally
	compact totally disconnected.
\end{prop}

The following is a very important and fundamental tool in the research
of t.d.l.c. groups:
\begin{thm}[The Van Dantzig Lemma]
	Let $G$ be a totally disconnected locally compact group. Then it
	admits a basis of identity neighbourhoods consisting of compact open
	subgroups.
\end{thm}

\section{Inverse Limits}

An important construction of topological groups is given by \emph{inverse
	limits}.\footnote{Inverse limits are used to construct a lot of objects in
	mathematics, and not just groups.}
\begin{defn}
	Let $I$ be a partially ordered set, and suppose that we have a topological
	group $G_{i}$ for every $i\in I$, as well as a continuous epimorphism
	\[
	f_{ij}:G_{i}\to G_{j}
	\]
	for very $i>j$ in $I$, such that 
	\[
	f_{jk}\circ f_{ij}=f_{ik}
	\]
	for every $i>j>k$ in $I$. Then the group 
	\[
	\lim_{\longleftarrow}G_{i}\coloneqq\left\{ (g_{i})\in\prod_{i\in I}G_{i}\middle|f_{ij}(g_{i})=g_{j}\quad\forall i>j\right\} 
	\]
	is called the \emph{inverse limit }of the \emph{inverse system }$(G_{i},f_{ij})_{ij}$.
\end{defn}

\begin{exe}
	Show that the inverse limit of an inverse system of topological groups
	is always a topological group. Hint: it's not hard.
\end{exe}

Observe that, for every $i\in I$, the group $G_{i}$ is a quotient
of the group $\underset{\longleftarrow}{\lim}G_{i}$. 

The next theorem demonstrates the strength of this construction.
\begin{thm}
	Every connected locally compact group is (isomorphic to) an inverse
	limit of Lie groups.
\end{thm}

\begin{cor}
	\label{invlim}
	Every compact group is (isomorphic to) an inverse limit of (compact)
	Lie groups.
\end{cor}

The reason we put ``compact'' in parentheses is that, if a compact
group is an inverse limit of an inverse system of groups, then these
groups must be compact as well (as continuous images of it). Alternately,
the following is not hard to prove:
\begin{exe}
	Show that an inverse limit of compact groups is compact. Show that
	an inverse limit of totally disconnected groups is totally disconnected.
\end{exe}

In particular, an inverse limit of finite groups must be compact and
totally disconnected. In fact, it turns out that the converse is true
as well!
\begin{defn}
	A topological group is \textbf{\emph{profinite }}if it is (isomorphic
	to) an inverse limit of finite groups
\end{defn}

\begin{exe}
	A topological group is profinite if and only if it is compact and
	totally disconnected.
\end{exe}
\begin{hint}
    First show that the inverse limit of compact groups is compact, and that the inverse limit of totally disconnected groups is totally disconnected. It then follows that the inverse limit of finite groups is both compact and totally disconnected.

    Conversely, if $G$ is a compact, totally disconnected group, then it admits a basis of identity neighbourhoods consisting of open subgroups. If $H\leqslant G$ is an open subgroup, then its cosets form a disjoint covering of $G$ by open sets; hence, there are only finitely many of them. Show that this implies that $H$ admits a finite index subgroup that is normal inside $G$. (In general, if $A$ is a finite index subgroup of a group $B$, then the intersection of all the conjugates of $A$ inside $B$ is a normal, finite index subgroup of $B$.) Deduce that $G$ admits a basis consisting of \textit{normal} open subgroups.

    Now, show that one may form an inverse limit from the quotients of $G$ by the subgroups in this basis, and find a natural continuous map from $G$ to this inverse limit. Show that it is closed and injective. Deduce that $G$ is profinite.
\end{hint}

\begin{rem}
	Another way to prove this fact is to use Corollary \ref{invlim}, as well as the fact (which is nontrivial) that quotients of compact totally disconnected groups are necessarily
	totally disconnected as well.
\end{rem}

\begin{exa}
	Fix some prime $p$, and, for every $n\in\mathbb{N}$, set $G_{n}=\mathbb{Z}/p^{n}\mathbb{Z}$.
	For $n>m$, let $f_{n,m}:G_{n}\to G_{m}$ be the natural projection.
	Then $\underset{\longleftarrow}{\lim}G_{n}$ is called the group of
	$p$-adic integers, denoted by $\mathbb{Z}_{p}$. Thinking of $G_{n}$
	as a ring and not just as a group, we obtain a ring structure on $\mathbb{Z}_{p}$
	as well, making it into an integral domain. Its fraction field is
	denoted by $\mathbb{Q}_{p}$, the field of $p$-adic numbers. It is locally compact and totally disconnected.
\end{exa}

Here is an example of an inverse limit of connected Lie groups.
\begin{exa}
	For $n\in\mathbb{N}$, set $G_{n}=\mathbb{R}/\mathbb{Z}$. We order
	$\mathbb{N}$ by the relation of divisibility. If $n|m$, we set $f_{m,n}(x+\mathbb{Z})=\frac{m}{n}x+\mathbb{Z}$.
	The group $\underset{\longleftarrow}{\lim}G_{n}$ is called \emph{the
		Solenoid}. As an abstract group, it is isomorphic to $(\mathbb{R},+)$,
	but its topology is very different from the standard topology on $\mathbb{R}$.
	It is connected and compact, but not path connected. In fact, as an abstract group, this group is also isomorphic to $(\mb Q_p,+)$, $(\mb C,+)$, and the additive group of any field of characteristic zero and cardinality $2^{\aleph_0}$.
\end{exa}

\section{Haar Measure}

Every locally compact group admits a natural (and very useful) measure,
called its \emph{Haar measure}. Before we define it, we need some
preliminaries. We recall some definitions:
\begin{defn}
	If $X$ is a topological space, then the $\sigma$-algebra generated
	by its open sets is called the \emph{Borel $\sigma$-algebra. }Subsets
	belonging to this $\sigma$-algebra are called \emph{Borel. }A measure
	defined on this $\sigma$-algebra is called a \emph{Borel measure}. 
	We assume all measures are nontrivial, in the sense that they neither 
	assign the measure $0$ to all subsets, nor assign the measure
	$\infty$ to all nonempty subsets.
\end{defn}
Recall the following definition:
\begin{defn}
	Given a Borel measure $\mu$, a Borel subset $A$ is called \emph{inner
		regular }if $$\mu(A)=\sup\left\{ \mu(F)\middle|F\subseteq A,\text{ \ensuremath{F} is compact}\right\}.$$
	It is called \emph{outer regular }if $$\mu(A)=\inf\left\{ \mu(V)\middle|V\supseteq A,\text{ \ensuremath{V} is open}\right\}.$$
	The measure $\mu$ is called \textit{regular} if every open subset is inner
	regular and every measurable subset is outer regular.
\end{defn}

\begin{defn}
	Let $G$ be a topological group. A measure $\mu$ on $G$ is called
	\textbf{\emph{left-invariant }}if $\mu(gA)=\mu(A)$ for every $g\in G$
	and every measurable $A\subseteq G$ (where $gA=\left\{ ga\middle|a\in A\right\} $).
	A left-invariant regular Borel measure is called a \textbf{\emph{Haar
			measure}}.
\end{defn}

\begin{thm}[Haar]
	Every locally compact group admits a Haar measure. This measure is
	unique up to a multiplicative constant.
\end{thm}

It is a theorem of Weil that \emph{only }locally compact groups admit a
Haar measure.
\begin{rem}
	Sometimes we will also consider the completion of the Haar measure,
	which is defined on a larger $\sigma$-algebra than the Borel $\sigma$-algebra.
	This completion is often still called a Haar measure.
\end{rem}

\begin{exe}
	Let $G$ be a locally compact group with a Haar measure $\mu$. Show that
	$\mu(V)>0$ for every open $V\subseteq G$ and that $\mu(K)<\infty$
	for every compact $K\subseteq G$.
\end{exe}

\begin{exe}
	Let $G$ be a locally compact group, and let $\mu$ be a Borel measure
	on $G$ (not necessarily a Haar measure). Prove that the following
	are equivalent:
	\begin{enumerate}
		\item If $K\subseteq G$ is compact, then $\mu(K)<\infty$.
		\item The measure $\mu$ is \emph{locally finite}, in the sense that, for
		every $g\in G$, there is a neighbourhood $V\ni g$ for which $\mu(V)<\infty$.
	\end{enumerate}
\end{exe}

\begin{rem}
	In locally compact second-countable spaces, every locally finite Borel measure is regular. Therefore, since we always assume our locally compact groups are second-countable, we won't need to worry about regularity too much.
\end{rem}

\begin{rem}
        We demanded left invariance in the definition, but if $\mu$ is a left-invariant Haar measure, then the measure $J$ defined by $$J\l(A\r)=\mu\l(A^{-1}\r)$$ is right-invariant, and satisfies all other conditions of  Haar measures. Such a measure is called a \emph{\textbf{right Haar measure}}. A locally compact group is called \emph{\textbf{unimodular}} if its left Haar measure is right-invariant as well.
\end{rem}
\begin{exa}
        For discrete groups, the counting measure is a Haar measure.
\end{exa}
\begin{exa}
        The Lebesgue measure is the Haar measure of $\l(\mb R^n,+\r)$.
\end{exa}
\begin{exe}
        \label{fm:c}
        Let $G$ be a locally compact group, $\mu$ its Haar measure. Show $G$ is compact if and only if $\mu\l(G\r)<\infty$.
\end{exe}

\section{Amenability of Topological Groups}
There is a notion of amenability of topological groups. This notion is different (and weaker) then the previously defined notion of amenability. In fact, a group is amenable in the previous sense if and only if it's amenable in the new sense when given the discrete topology.

\subsection{\texorpdfstring{$L_p$}{L\_p} spaces}
Let $\l(X,m\r)$ be a measure space, and let $1\le p<\infty$. Consider the non-negative function $\l\|\cd\r\|_p$ defined on all measurable functions from $X$ to $\mb R$, defined by
$$\l\|f\r\|_p=\l(\int_X\!\l|f\r|^p\r)^{1/p}.$$
The set of all functions $f$ such that $\l\|f\r\|_p<\infty$,
$$\l\{f:X\to \mb R\mi\int_X\!\l|f\r|^p\;\mr dm<\infty\r\},$$
forms a vector space with the natural operations. $\l\|\cd\r\|_p$ is not a norm in this space, because $\l\|f\r\|_p=0$ doesn't imply $f=0$. However, if we take the quotient of that space over the kernel of $\l\|\cd\r\|_p$, i.e.\ if we identify two functions $f,g$ if they agree almost everywhere,\footnote{That is, if $\{x\in X|f(x)\ne g(x)\}$ is of measure zero.} we get a new vector space, in which $\l\|\cd\r\|_p$ \emph{is} a norm.

This space, along with the norm $\l\|\cd\r\|_p$, is denoted by \emph{\textbf{$L_p\l(X,m\r)$}}, and is a Banach space.

\paragraph{}For $p=\infty$, we take $\l\|\cd\r\|_\infty$ to be the essential supremum,
$$\l\|f\r\|_\infty=\inf\l\{s\mi f\l(x\r)\le s\text{ almost everywhere}\r\}=\inf\l\{s\mi m\l(\l\{x\mi f\l(x\r)>s\r\}\r)=0\r\}.$$
Then we take the set of functions $f$ such that $\l\|f\r\|_\infty<\infty$ (the set of essentially bounded functions) quotient by the kernel of $\l\|\cd\r\|_\infty$ (which again means we identify functions which agree almost everywhere). This set, along with the natural operations and $\l\|\cd\r\|_\infty$, is denoted by \emph{\textbf{$L_\infty\l(X,m\r)$}} and is a Banach space as well.

\begin{rem}
        For all $1\le p\le \infty$, $\ell_p\l(X\r)=L_p\l(X,m\r)$ where $m$ is the counting measure (the measure assigning each set its size).
\end{rem}

\begin{rem}
        If $G$ is a locally compact group, we denote $L_p\l(G,\mu\r)$ (where $\mu$ is the Haar measure of $G$) by $L_p\l(G\r)$.
\end{rem}

\subsection{Amenability}
\begin{defn}
        Let $G$ be a locally compact group. A \emph{\textbf{mean}} on $G$ is a positive, norm $1$ functional on $L_\infty\l(G\r)$.
\end{defn}
The \emph{\textbf{left regular action}} of a group $G$ on $L_p\l(G\r)$ is defined by $$L_g\l(f\r)\l(x\r)=f\l(g^{-1}x\r)$$ for all $g\in G$, $f\in L_p\l(G\r)$, $x\in G$. A functional $\f\in L_p\l(G\r)^*$ is \emph{\textbf{left-invariant}} if $\f\l(L_g\l(f\r)\r)=\f\l(f\r)$ for all $g\in G$, $f\in L_p\l(G\r)$.
\begin{defn}
        A locally compact group $G$ is \emph{\textbf{amenable}} (as a topological group) if there exists a left-invariant mean on $G$.
\end{defn}
We had a lot of equivalent definitions for the original notion of amenability. We want to define analogous notions in topological groups.

\begin{defn}
        An action of a topological group $G$ on a topological space $X$ is \emph{\textbf{continuous}} if the two-variable map $\l(g,x\r)\mapsto gx$ from $G\x X$ to $X$ is continuous.
\end{defn}
\begin{rem}
        Previously when a group $G$ acted on a topological space $X$ we demanded each $g:X\to X$ for $g\in G$ to be continuous. Clearly, the definition above is stronger: it requires $\l(g,x\r)\mapsto gx$ to be continuous as a two-variable map, and not just continuous in $x$.
\end{rem}
\begin{defn}
        A topological group $G$ has \emph{\textbf{almost invariant functions}} (shortened by a.i.f.) if, for every compact $K\subseteq G$ and for every $\ep>0$, there exists $f\in L_p\l(G\r)\setminus\{0\}$ such that
        $$\l\|L_g\l(f\r)-f\r\|_p\le \ep\l\|f\r\|_p$$
         for every $g\in K$.
\end{defn}
The following theorem is left as an exercise.
\begin{thm}
        \label{top:am}
        Let $G$ be a locally compact group. The following are equivalent:
        \begin{enumerate}[1)]
                \i \label{ta:1} $G$ is amenable (that is, admits a left-invariant mean).
                \i \label{ta:2} $G$ admits a finitely additive, left-invariant measure defined on all Haar measurable sets.\footnote{Recall the Haar measure is a completion of a Borel measure.}
                \i \label{ta:3} $G$ admits almost invariant functions in $L_2\l(G\r)$.
                \i \label{ta:4} $G$ admits almost invariant functions in $L_1\l(G\r)$.
                \i \label{ta:5} Whenever $G$ acts continuously on a compact Hausdorff space $K$, $K$ admits a $G$-invariant Borel probability measure.
                \i \label{ta:6} Whenever $G$ acts continuously by affine maps on a compact-convex set, it admits a global fixed point.
        \end{enumerate}
\end{thm}
\begin{exe}
        Prove Theorem \ref{top:am}. 
\end{exe}
\begin{hint}
        Some of the implication are very similar to the proofs we had for discrete groups: (\ref{ta:1})$\iff$(\ref{ta:2}); (\ref{ta:3})$\iff$(\ref{ta:4}); (\ref{ta:5})$\iff$(\ref{ta:6}); (\ref{ta:4})$\then$(\ref{ta:2}); (\ref{ta:6})$\then$(\ref{ta:1}). The tough part is (\ref{ta:2})$\then$(\ref{ta:4}).
\end{hint}
\section{Lattices}
\begin{defn}
A \emph{\textbf{lattice}} in a locally compact group $G$ is a discrete subgroup $\G$ such that the topological space $G/\G$ admits a $G$-invariant Borel probability measure.
\end{defn}

\begin{exe}
        Show that $\mb Z^n$ is a lattice inside $\mb R^n$.
\end{exe}
The next example discusses the Hyperbolic plane $\mb H^2$. If you don't know what this is, you may skip it without loss of continuity. A subgroup $H\le G$ is called \textit{\textbf{cocompact}} if the topological quotient space $G/H$ is compact.
\begin{exe}
        $\mr{SL}_2\l(\mb R\r)\acts \mb H^2$ by isometries. There is a hyperbolic structure on $\S$, the surface of genus $2$, such that the fundamental group $\pi_1\l(\S\r)$ is embedded in $\mr{SL}_2\l(\mb R\r)$ as a co-compact lattice.
\end{exe} 
\begin{exa}
        The subgroup $\mr{SL}_n\l(\mb Z\r)$ of $\mr{SL}_n\l(\mb R\r)$ is a lattice, which is not cocompact.
\end{exa}
\begin{exa}
        Take $G=\mr{Aut}(T_3)$, where $T_3$ is the $3$-regular tree. Colour the edges by three colours in a `legal' way (that is, such that the three edges coming out of every vertex have different colours), and let $\G\le G$ be the group of all automorphisms that preserve this colouring. $G/\G$ is the space of legal colourings, and $\G$ is a lattice.
\end{exa}
\begin{thm}\
        \begin{enumerate}[1)]
                
                \i A closed subgroup of an amenable topological group is amenable.
                \i If $\G\le G$ is a lattice, then $\G$ is amenable if and only if $G$ is amenable.
        \end{enumerate}
\end{thm}
\begin{proof}\
        \begin{enumerate}[1)]
                \i We're going to need the following fact, which we will not prove.
                \begin{fact}
                        If $H$ is a closed subgroup of the locally compact group $G$, and $G/H$ is the set of right cosets, then there exists a set $A$ of representatives, which contains precisely one representative of each coset, such that
                        $$\psi:A\x H\to G\quad \l(a,h\r)\overset{\psi}{\mapsto} ah$$ is a measurable isomorphism.
                \end{fact}
                Now, suppose $G$ is amenable, $H\le G$ closed. Let $m$ be a right invariant mean on $G$; we will use it to construct an invariant mean on $H$.
                
                For each $f\in L_\infty\l(H\r)$ define $\t f\in L_\infty\l(G\r)$ by setting, for all $g\in G$,
                $$\t f\l(g\r)=f\l(h\r),$$
                where $\psi^{-1}\l(g\r)=\l(a,h\r)$ for some $a\in A$. That is, $\t f=f\circ p_2\circ \psi^{-1}$, where $p_2:A\x H\to H$ is the projection onto the second coordinate.
                
                Now, define $\nu$ by $\nu\l(f\r)=m(\t f)$, and $\nu$ is a right invariant mean on $L_\infty\l(H\r)$.

                \i If $G$ is amenable and $\G\le G$ is a lattice, then $\G$ is amenable as a closed subgroup, using
                \begin{exe}
                        Let $G$ be a group, $H\le G$ a discrete subgroup. Prove $H$ is closed (i.e., has no accumulation points).
                \end{exe}
                \begin{hint}
                        If $\g_n\too x$ then $\g_n^{-1}\g_m\too 1$.
                \end{hint}
                Let $\G\le G$ be an amenable lattice. We prove $G$ is amenable.
                
                Let $G$ act continuously by affine maps on some compact-convex set $A$, and let $p\in A$ be a $\G$-fixed point.
                
                Push the probability measure on $G/\G$ to $A$ via the orbit map of $p$ defined by $g\mapsto gp$ (which is well defined since $p$ is fixed by $\G$). We get a probability measure on $A$ which is $G$-invariant. Its barycentre is a $G$-fixed point.\qedhere
        \end{enumerate}
\end{proof}
\begin{rem}
        Notice that in the discrete topology, all subgroups are closed, so the last theorem implies all subgroups of discrete amenable groups are amenable, which we already stated in Theorem \ref{ele}.
\end{rem}
\begin{rem}
Every group admitting an amenable lattice is amenable, but not every amenable group admits a lattice. There are solvable groups which are not unimodular, and hence don't admit any lattices (see Corollary \ref{latuni}).
\end{rem}
\begin{cor}
        \label{C:A}
        Compact groups are amenable.
\end{cor}
\begin{proof}
        The trivial group is a lattice, because $G/\l\{e\r\}\cong G$ is compact and its Haar measure is $G$-(left-)invariant and finite.
\end{proof}
\begin{exa}
        $\rot$ is compact, and hence it's amenable (as a topological group).
        Recall $\rot$ is \emph{not} amenable when given the discrete topology, because it contains $F_2$.
\end{exa}

\begin{rem}
	A very useful property of Lie groups is that every connected Lie group $G$ admits a maximal normal solvable subgroup. That is, a normal subgroup that is solvable and that contains every other normal solvable subgroup. This subgroup is called the \textit{radical} of $G$, and is sometimes denoted by $\mr R(G)$. The quotient $G/\mr R(G)$ has a trivial radical; in other words, it admits no normal solvable subgroup except the trivial subgroup. A Lie group with a trivial radical is called \textit{semisimple}. The centre of a semisimple Lie group is discrete, and the quotient of a semisimple Lie group by its centre is always a product of simple (centre-free) simple Lie groups. Thus, questions about Lie groups are often separated to the solvable case and to the simple case. 
    
    In general, topological (or abstract) groups do not admit a maximal normal solvable subgroup. However, one can use a similar notion: it is not hard to show that every topological group admits a maximal normal \textit{amenable} subgroup. This subgroup is called the \textit{\textbf{amenable radical}} of $G$. If $G$ is a locally compact group and $A$ is its amenable radical, then the amenable radical of $G/A$ is trivial. Moreover, $G/A$ admits a finite index subgroup $\hat{G}$ such that $\hat{G}\cong S\times H$, where $S$ is a connected semisimple Lie group without compact factors, and $H$ is a totally disconnected locally compact group. 
\end{rem}

\section*{Further Reading}
The theory of topological groups is vast, with many possible directions for pursuing. In order to explore Lie groups further, readers may consult \cite{Varadarajan,hall2013lie}, \cite{knapp1996lie} and \cite{onishchik2012lie} (and many more books). For a treatment delving into the theory of algebraic groups as well, see \cite{onishchik1993lie}.
For the theory of locally compact groups and Haar measures, the reader is referred to \cite{folland2016course}. 
For Hilbert's fifth problem, see \cite{Kaplansky,MoZi}.
For more general topological groups (and even more general topological structures), the reader may refer to \cite{arhangel2008topological}.
The recommended source about lattices in Lie groups is Raghunathan's book \cite{Rag}. For a crash course on lattices see \cite{Gel-lattices}.

\setcounter{section}{0}
\chapter{Property (T)}

\section{Definition and First Results}

Given a Hilbert space $\mathcal{H}$, we denote by $\mathcal{B}(\mathcal{H})$
the algebra of bounded operators on $\mathcal{H}$. The \textbf{\emph{strong
operator topology}} on $\mathcal{B}\left(\mathcal{H}\right)$ is the
weakest topology such that, for every $v\in\mathcal{H}$, the evaluation
map 
\begin{align*}
\delta_{v}:\mathcal{B}\left(\mathcal{H}\right) & \to\mathcal{H}\\
A & \mapsto A(v)
\end{align*}
is continuous. In this topology, a sequence $\left(U_{n}\right)_{n\in\mathbb{N}}$
of operators converges to $U$ if and only if it converges pointwise.
That is, $U_{n}\longrightarrow U$ in $\mathcal{B}(\mathcal{H})$
if and only if $U_{n}(v)\longrightarrow U(v)$ in $\mathcal{H}$ for
every $v\in\mathcal{H}$. \footnote{If we endow $\mathcal{H}^{\mathcal{H}}$ (the set of all maps from
$\mathcal{H}$ to itself) with the product topology, and consider
$\mathcal{B}(\mathcal{H})$ as a subset of $\mathcal{H}^{\mathcal{H}}$,
then the subspace topology on $\mathcal{B}(\mathcal{H})$ is the same
as the strong operator topology.}

We also endow $\mc U(\mc H)$, the group of unitary linear transformations from $\mc H$ to itself (which is a subset of $\mc B(\mc H)$), with the subspace topology. This makes $\mc U(\mc H)$ into a topological group, which is second-countable but not locally compact. 
\begin{defn}
A \textbf{\emph{unitary representation}} of a group $G$ is a continuous
homomorphism $\rho:G\to U\left(\mathcal{H}\right)$ (for some Hilbert space $\mc H$). 
\end{defn}

This homomorphism is continuous if and only if the action $G\acts\mathcal{H}$
defined by $\left(g,v\right)\mapsto\rho\left(g\right)\left(v\right)$
is continuous. At times, in order to simplify notations, we will omit
$\rho$, and simply write (for example) $g.v$ instead of $\rho(g)(v)$.
\begin{defn}
Let $G$ be a group acting on a Hilbert space $\mathcal{H}$. Let
$K\subseteq G$ be a subset and $\varepsilon>0$. A nonzero vector
$v\in\mathcal{H}\backslash\left\{ 0\right\} $ is called \textbf{\emph{$(K,\varepsilon)$-invariant
}}if, for every $k\in K$, we have
\[
\frac{\left\Vert k.v-v\right\Vert }{\left\Vert v\right\Vert }<\varepsilon.
\]
It is called \textbf{\emph{$K$-invariant}} if it is fixed by all
$k\in K$ (i.e., $k.v=v$), and simply \textbf{\emph{invariant}} if
it is $G$-invariant.
\end{defn}

\begin{defn}
A unitary representation of $G$ has \textbf{\emph{almost invariant
vectors}} (shortened by a.i.v.) if, for every compact $K\subseteq G$
and for every $\ep>0$, there exists a $(K,\ep)$-invariant vector.
\end{defn}

\begin{defn}
A locally compact group $G$ is \textbf{\emph{Kazhdan}} (or, alternately, is said
to have \textbf{\emph{property (T)}}), if every unitary representation
of $G$ with almost invariant vectors admits an invariant vector.
\end{defn}

\begin{prop}
\label{C:T} Compact groups have property (T).
\end{prop}

\begin{proof}
Let $G$ be a compact group, $\left(\pi,\mc H\right)$ a unitary representation
with almost invariant vectors. Let $\mu$ be the Haar measure such that
$\mu\left(G\right)=1$.

Take some $v\in\mc H$, $\l\|v\r\|=1$, which is $\left(G,\sqrt{2}\right)$-invariant.
Define 
\[
\bar{v}=\int_{G}\!g.v\ \mathrm{d}\mu\left(g\right).
\]
We will show that $\bar v$ is a non-zero invariant vector.
\begin{rem}
We are taking an integral here of a map from $G$ to $\mathcal{H}$.
The reader may take it as an exercise to define its meaning precisely,
and to show that it behaves as well as integrals of functions into
$\mathbb{C}$ or $\mathbb{R}$.
\end{rem}
The vector $\bar{v}$ is the `averaging' of $v$ with respect to the
action of $G$, or, if you prefer such language, the expected value
of $g.v$ (with respect to the probability measure $\mu$). Whenever
one defines a vector in such a way, it is going to be $G$-invariant;
the problem is usually to show that it is nonzero. In this case, this
is easy to prove, since $g.v$ always belongs to the same halfplane
(since it is always at distance at most $\sqrt{2}$ from $v$). To
be more detailed:

To see that $\bar{v}$ is invariant, let $h\in G$, and consider the
following computation: 
\begin{align*}
h.\bar{v} & =h.\left(\int_{G}\!g.v\;\mathrm{d}\mu\left(g\right)\right)=\int_{G}\!\left(h\cd g\right).v\;\mathrm{d}\mu\left(g\right)\\
 & =\int_{G}\!g.v\;\mathrm{d}\mu\left(h^{-1}g\right)\\
 & =\int_{G}\!g.v\;\mathrm{d}\mu\left(g\right)=\bar{v},
\end{align*}
where the penultimate equality follows from the left-invariance of the
Haar measure.

Now, in order to see that $\bar{v}\neq0$, observe that $\l\|g.v-v\r\|<\sqrt{2}$
implies that $\l<g.v,v\r>>0$ (for every $g\in G$), so 
\[
\l<\bar{v},v\r>=\int_{G}\!\l<g.v,v\r>\;\mr{d}\mu\left(g\right)>0.
\]
Thus, $\bar{v}\ne0$. 

\end{proof} \begin{cor} A locally compact group is compact if and
only if it is both amenable and Kazhdan. \end{cor} \begin{proof}
One direction follows from Corollary \ref{C:A} and Proposition \ref{C:T}.
Conversely, suppose $G$ is both amenable and Kazhdan, and consider
its action on $L_{2}\left(G\right)$. By amenability, it admits almost
invariant vectors, so, by property (T), it admits an invariant vector.
This vector is a constant function, because the action of $G$ on
itself is transitive. It is also in $L_{2}$, so its integral is finite
and equals (up to multiplication by a constant) to the measure of
$G$. Thus the (Haar) measure of $G$ is finite, so $G$ is compact
by Exercise \ref{fm:c}. \end{proof} 
Thus, property (T) is somewhat of an `opposite' to amenability: a non-compact (locally compact) amenable group cannot have property (T).

\begin{prop} A locally compact Kazhdan group is compactly generated.
In particular, a discrete Kazhdan group is finitely generated. \end{prop}

For the proof, we will use direct sums of Hilbert spaces. If $I$ is a set and $\mc H_\alpha$ is a Hilbert space for each $\alpha\in I$,
                the \emph{\textbf{direct sum}} $\bigoplus_{\alpha\in I}\mc H_\alpha$ of $\l(\mc H_\alpha\r)_{\alpha\in I}$ is the set 
                $$\l\{\l(v_\alpha\r)_{\alpha\in I}\mi \sum_{\alpha\in I}\l\|v_\alpha\r\|^2<\infty\r\}$$
                with the inner product
                $$\l<\l(v_\alpha\r),\l(u_\alpha\r)\r> = \sum_{\alpha\in I} \l<v_\alpha,u_\alpha\r>.$$
                The requirement $\sum\l\|v_\alpha\r\|^2<\infty$ requires in particular that $v_\alpha\ne 0$ for only countably many values of $\alpha$.

\begin{proof} Let $G$ be a locally compact Kazhdan group. We prove
it's generated by some compact set.

Let $O\subseteq G$ be a compact set with non-empty interior and that
contains the identity element $1$ (the existence of which is guaranteed
by local compactness), and for any compact $K$, define $G_{K}\coloneqq\l<KO\r>$.
$KO$ has nonempty interior (because $O$ has nonempty interior and
$KO=\bigcup_{k\in K}kO$), and it is compact (because $K\x O\subseteq G\x G$
is compact and multiplication is continuous). \begin{exe} \label{Tos}
Let $H\le G$ be a subgroup with nonempty interior. Prove it's open.
\end{exe}

Thus, $G_{K}$ is open and $G/G_{K}$ is discrete. This means that $L_2(G/G_K)$ and $\ell_2(G/G_K)$ are the same, since the measure on $G/G_K$ is equivalent to the counting measure. Consider 
 $$\mc H\coloneqq
 \bigoplus_{K\text{ compact}} \ell_2\l(G/G_K\r).$$

        $\mc H$ is a huge Hilbert space; each vector in $\mc H$ is of the form $\l(v_K\r)_{K\text{ compact}}$ where $v_K\in\ell_2\l(G/G_K\r)$. The left regular action of $G$ on $\mc H$, defined by $$L_\g\cd \l(v_K\r)_{K\text{ compact}}=\l(L_\g\cd v_K\r)_{K\text{ compact}},$$ is a unitary representation. It has almost invariant vectors -- in fact, for each compact $C$ there exists a $C$-invariant vector: take $\l(v_K\r)$ such that $v_K=0$ for $K\ne C$ and $v_C\in \ell_2\l(G/G_C\r)$ is the Dirac function of the trivial coset (i.e., $v_C\l(gG_C\r)=0$ for $g\notin G_C$ and $v_C\l(G_C\r)=1$). By property (T) there exists a globally invariant function $\l(u_K\r)$; each $u_K\in \ell_2\l(G/G_K\r)$ is a constant function because the action of $G$ on $G/G_K$ is transitive. By definition $\l(u_K\r)\ne 0$, so there is some $K$ such that $u_K\ne 0$. Thus, this $u_K$ is constant and in $\ell_2\l(G/G_K\r)$, so the set $G/G_K$ must be finite (otherwise the sum wouldn't converge).

$G$ is generated by $G_K$ and the representatives in $G/G_K$. $G_K$ is generated by $KO$ which is compact, and we just showed $G/G_K$ is finite so we can take a finite number of representatives from it. Hence, $G$ is generated by the union of $KO$ with some finite set, which is also a compact set.
\end{proof}

We could have defined property (T) using the notion of \textit{Kazhdan pairs}.
\begin{defn}
Let $G$ be a locally compact group, $K\subseteq G$ and $\ep>0$. We say $\l(K,\ep\r)$ is a \emph{\textbf{Kazhdan pair}} if for all unitary representations $\rho:G\to U\l(\mc H\r)$, the existence of a $\l(K,\ep\r)$-invariant unit vector implies existence of an invariant vector.

If $\l(K,\ep\r)$ is a Kazhdan pair, we say $K$ is a \emph{\textbf{Kazhdan set}} and $\ep$ is its \emph{\textbf{Kazhdan constant}}.
\end{defn}
Using this terminology, we could have defined property (T) differently.
\begin{thm}[Alternative definition for property (T)]\

A locally compact group $G$ is Kazhdan if and only if there exists a compact $K\subseteq G$ and $\ep>0$ such that $\l(K,\ep\r)$ is a Kazhdan pair.
\end{thm}
\begin{exe}
Prove this theorem.
\end{exe}
In fact, we have the following:
\begin{prop}\
\label{T:gen}
\begin{enumerate}[(1)]
\i \label{gen:1} If $G$ has property (T) and $G=\l<K\r>$ then $K$ is a Kazhdan set.
\i \label{gen:2} If $K\ss$ G is a Kazhdan set with nonempty interior then $G=\l<K\r>$.
\end{enumerate}
\end{prop}
\begin{proof}[Proof of (\ref{gen:1})]
We may assume $K$ is symmetric (i.e.\ $K=K^{-1}$), because in general $$\l\|\rho\l(g\r)v-v\r\|=\l\|\rho\l(g^{-1}\r)v-v\r\|$$ for all unitary representations $\rho$. We may also assume it contains the identity, $1\in K$. In fact, we may even assume $K$ is closed.
\begin{exe}
\label{cl-K}
Let $G$ be a locally compact group, $K\ss G$. Prove that if $\overline{K}$ is a Kazhdan set, then $K$ is a Kazhdan set.
\end{exe}
\begin{exe}
\label{cov-K}
Show that, if $\l(K^n,\ep\r)$ is a Kazhdan pair, then $\l(K,\ep/n\r)$ is also a Kazhdan pair.
\end{exe}
Let $M$ be some compact Kazhdan set in $G$. Observe it's enough for us to prove there exists some $n\in\mb N$ such that $K^n\supseteq M$, since it will follow that $K^n$ is Kazhdan, and hence (by Exercise \ref{cov-K}) $K$ is Kazhdan.

We know $\bigcup_{n\in\mb N} K^n$ covers $G$ and hence $M$. If it were an open cover, we could take a finite subcover and be done. In fact, it would have been enough if $K$ had nonempty interior: then $\bigcup_{n\in\mb N}(\mr{Int}K)K^n$ is an open cover (because whenever $U$ is an open subset and $S$ is any subset, $US=\bigcup_{s\in S}Us$ is an open subset as well), hence has a subcover $\bigcup_{n=1}^{n_0}(\mr{Int}K)K^n$, and hence $\bigcup_{n=2}^{n_0+1}K^n=K^{n_0+1}$ covers $M$ as well (because $\mr{Int}K\subseteq K$ and hence $(\mr{Int}K)K^i\subseteq K^{i+1}$).

We don't know $K$ has nonempty interior, but we do know there exists a $k$ such that $\overline{K^k}$ has nonempty interior, because otherwise $\bigcup_{n\in\mb N}\overline{K^n}$ would have been a cover of $G$\ (which is locally compact and Hausdorff, and hence a Baire space) of closed sets with nonempty interior, a contradiction. Then we can take $\t K\coloneqq \overline{K^k}$, deduce by the above that $\t K$ is Kazhdan, and hence (by Exercise \ref{cl-K}) that $K^k$ is Kazhdan, and hence (by Exercise \ref{cov-K}) that $K$ is Kazhdan.\footnote{In fact, more is true, although we didn't need it for the proof: by the Steinhaus theorem, there is some $k\in \mb N$ such that $K^k$ has nonempty interior, and hence there is some $N\in\mb N$ such that $K^N$ contains $M$. For this, one does not even need to assume $K$ is closed, only that it is measurable (and symmetric, and contains the identity).}

\end{proof}
\begin{proof}[Proof of (\ref{gen:2})]
Let $H=\l<K\r>$. By Exercise \ref{Tos}, $H$ is open. By the following exercise, it's closed.
\begin{exe}
Prove open subgroups are always closed.
\end{exe}
Since $H$ is open, $G/H$ is discrete. Consider $L_2\l(G/H\r)=\ell_2\l(G/H\r)$ and the left regular action of $G$ on $H$. The Dirac function supported on $eH$,
$$f\l(gH\r)=\begin{cases}1 &g\in H\\0 &g\notin H\end{cases},$$is $H$-invariant and hence $\l(K,\ep\r)$-invariant for all $\ep>0$, so (by definition of a Kazhdan set) we deduce there exists an invariant function in  $\ell_2\l(G/H\r)$. Since the action of $G$\ on $G/H$ is transitive, such a function must be constant, and since it's in $\ell_2\l(G/H\r)$ the sum of its values must converge -- and thus, $\l|G/H\r|$ must be finite. We go on to prove it equals $1$.

Consider
$$\ell_2^0\l(G/H\r)=\l\{f\in\ell_2\l(G/H\r)\mi \sum_{gH\in G/H}f\l(gH\r)=0\r\},$$ the orthogonal space to the space of constant functions. $\ell_2^0\l(G/H\r)$ is $G$-invariant (i.e.\ $G\l(\ell_2^0\l(G/H\r)\r)\ss\ell_2^0\l(G/H\r)$).

Assume by contradiction $\l|G/H\r|>1$. Then $\ell_2^0\l(G/H\r)$ is nonzero, and we have a nonzero $H$-invariant function, $$f\l(gH\r)=\begin{cases}1 &g\in H\\-\fr{1}{\l|G/H\r|-1} &g\notin H\end{cases},$$so by property (T) there exists a nonzero invariant function in $\ell_2^0\l(G/H\r)$. Such a function is constant (and nonzero), and satisfies $$\sum_{gH\in G/H} f\l(gH\r)=0,$$which is a contradiction. Hence, $\l|G/H\r|=1$ and $G=H=\l<K\r>$.
\end{proof}
By definition, if $\l(K,\ep\r)$ is a Kazhdan pair and $v$ a unit $\l(K,\ep\r)$-invariant vector, then there exists a globally invariant vector. The next lemma shows we can pick an invariant vector which is close to the original $v$.
\begin{lem}
\label{clo:inv}
Let $\l(K,\ep\r)$ be a Kazhdan pair in a group $G$ and let $\rho:G\to U\l(\mc H\r)$ be a unitary representation. If $\delta\in(0,1)$ and $v\in\mc H$ is a unit vector that is $\l(K,\ep\delta\r)$-invariant, then there exists an invariant vector $v^\pr$ such that $\l\|v-v^\pr\r\|<\delta$.
\end{lem}
\begin{proof}
Write $\mc H=\mc H^\pr\oplus \mc {H^\pr}^\perp$ where $\mc H^\pr$ is the subspace of invariant vectors, $\mc {H^\pr}^\perp$ its orthogonal complement. Both are invariant under $\rho\l(G\r)$.

Let $\pi$ be the projection onto $\mc H^\pr$. We claim $\pi\l(v\r)$ is the desired $v^\pr$.

Let $p$ be the projection onto ${\mc H^\pr}^\perp$, $p=1-\pi$. We need to show $\l\|p\l(v\r)\r\|=\l\|v-\pi\l(v\r)\r\|<\delta$. We know $p\l(v\r)$ is also $\l(K,\delta\ep\r)$-invariant (because $\l\|p\r\|\le 1$), and therefore $\fr{p\l(v\r)}{\l\|p\l(v\r)\r\|}$ is $(K,\fr{\delta\ep}{\l\|p\l(v\r)\r\|})$-invariant. Therefore $\fr{\delta}{\l\|p\l(v\r)\r\|}> 1$, as we wanted (otherwise $\fr{p\l(v\r)}{\l\|p\l(v\r)\r\|}$ would have been $\l(K,\ep\r)$-invariant, and we would have deduced there is an invariant vector in ${\mc H^\pr}^\perp$, a contradiction).
\end{proof}
\begin{exe}
Let $G$ be a group with property (T), and let $\phi:G\to H$ be a continuous homomorphism. Prove that $\overline{\phi(G)}$ has property (T).
\end{exe}
In particular, quotient groups of a group with property (T) have property (T). Since abelian groups are always amenable, and amenable groups with property (T) are compact, We get the following useful corollary:
\begin{cor}
If $G$ is locally compact and has property (T), then $G/\overline{G^\pr}$ is compact.\footnote{Recall that $G'=\l<\{[g,h]|g,h\in G\}\r>$ is the commutator subgroup. Here $\overline{G'}$ is its closure, which is the topological analogue of the commutator subgroup: $G/G'$ is the ``biggest" Hausdorff abelian quotient of $G$.}
\end{cor}
\section{Lattices and Induced Representations}
\subsection{Lattices}

Let $G$ be a locally compact group. A subgroup $H\leqslant G$ is
called \textbf{\emph{cofinite }}if $G/H$ admits a finite Borel (regular)
$G$-invariant measure (with respect to the standard action of $G$
on $G/H$).\footnote{We wrote `regular' in parentheses since, since we always assume $G$
is second-countable, it follows that every finite Borel measure on
$G/H$ (which is also locally compact second-countable) is regular.} Recall that a \textbf{\emph{lattice }}is a cofinite, discrete subgroup. 

If $H\leqslant G$ is any subgroup of a group $G$, a measurable subset $\Omega\subseteq G$
is a \textbf{\emph{(left) fundamental domain }}of $H$ if $G=\bigcupdot_{h\in H}\Omega h$,
a disjoint union. It is a \textbf{\emph{right fundamental domain }}of
$H$ if $G=\bigcupdot_{h\in H}h\Omega$ (still a disjoint union).
Observe that $\Omega$ is a left fundamental domain if and only if
$\Omega^{-1}=\left\{ \omega^{-1}\middle|\omega\in\Omega\right\} $
is a right fundamental domain.
\begin{exe}
Show that a every discrete subgroup admits a fundamental domain. 
\end{exe}

\begin{exe}
Show that a discrete subgroup is cofinite if and only if it admits
a fundamental domain of finite Haar measure.
\end{exe}

\begin{lem}
If $\Gamma\leqslant G$ is a discrete subgroup and $\Omega_{1},\Omega_{2}\subseteq G$
are right fundamental domains of $\Gamma$, then $\mu(\Omega_{1})=\mu(\Omega_{2})$
(for every Haar measure $\mu$).
\end{lem}

\begin{proof}
We have the following equalities:
\begin{align*}
\Omega_{1} & =\Omega_{1}\cap G=\Omega_{1}\cap\left(\bigcupdot_{\gamma\in\Gamma}\gamma\Omega_{2}\right)=\bigcupdot_{\gamma\in\Gamma}\left(\Omega_{1}\cap\gamma\Omega_{2}\right),\\
\Omega_{2} & =G\cap\Omega_{2}=\left(\bigcupdot_{\gamma\in\Gamma}\gamma^{-1}\Omega_{1}\right)\cap\Omega_{2}=\bigcupdot_{\gamma\in\Gamma}\left(\gamma^{-1}\Omega_{1}\cap\Omega_{2}\right).
\end{align*}
By left invariance of $\mu$, we have $\mu\left(\gamma^{-1}\Omega_{1}\cap\Omega_{2}\right)=\mu\left(\Omega_{1}\cap\gamma\Omega_{2}\right)$
for every $\gamma\in\Gamma$. Therefore,
\[
\mu(\Omega_{1})=\sum_{\gamma\in\Gamma}\mu(\Omega_{1}\cap\gamma\Omega_{2})=\sum_{\gamma\in\Gamma}\mu(\gamma^{-1}\Omega_{1}\cap\Omega_{2})=\mu(\Omega_{2}),
\]
as needed.
\end{proof}
\begin{cor}
\label{latuni}
If a group $G$ admits a lattice, it is unimodular.
\end{cor}

\begin{proof}
Let $\mu$ be a Haar measure on $G$; we need to show that it is right
invariant. Let $g\in G$; we need to show that $\mu(Ag)=\mu(A)$ for
every Borel subset $A\subseteq G$. Define a new measure, $\nu_{g}$,
on $G$, by setting
\[
\nu(A)=\mu(Ag)
\]
for every Borel subset $A\subseteq G$. It is easy to show that $\nu$
is a Haar measure on $G$. Therefore, there is some $c_{g}>0$ such
that $\nu_{g}=c_{g}\mu$. We need to show $c_{g}=1$. Therefore, it
is enough to find one Borel subset $A\subseteq G$ of finite, positive
$\mu$-measure for which $\mu(A)=\mu(Ag)=\nu_{g}(A)$. Let $\Omega$
be a right fundamental domain for a lattice $\Gamma$ of $G$. It
is immediate that $\Omega g$ is also a right fundamental domain of
$\Gamma$, since 
\[
G=G\cdot g=\left(\bigcupdot_{\gamma\in\Gamma}\gamma\Omega\right)\cdot g=\bigcupdot_{\gamma\in\Gamma}\gamma(\Omega g).
\]
Since $\mu(\Omega)$ is positive and finite, we are done.
\end{proof}
Let $\Gamma$ be a lattice in a locally compact group $G$, and let
$\Omega$ be a fundamental domain of $\Gamma$. By definition, there
is a finite Borel $G$-invariant measure on $G/\Gamma$. This measure,
just like the Haar measure on $G$, is unique up to scale.
\begin{fact}
Let $G$ be a locally compact group and $\Gamma\leqslant G$ a lattice.
Let $\nu_{1},\nu_{2}$ be two nontrivial finite Borel $G$-invariant
measures on $G/\Gamma$. There is some $c>0$ such that $\nu_{1}=c\cdot\nu_{2}$.
\end{fact}

For every such measure $\nu$ on $G/\Gamma$, there is a (unique)
Haar measure $\mu$ on $G$ such that
\[
\nu\l(A\r)=\mu\l(\pi^{-1}\l(A\r)\cap\Omega\r)
\]
for every Borel $A\subseteq G$, where $\pi$ is the projection $\pi:G\to G/\G$.
Often, we will assume for convenience of notation that $\nu(G/\Gamma)=1$. 
\begin{lem}
\label{fin:com}
If $\G\le G$ is a lattice and is finitely generated, then $G$ is compactly generated.
\end{lem}
\begin{proof}
 Denote by $\pi:G\to G/\G$ the quotient map. Let $\Om$ be a fundamental domain of $\G$. Let $\mu$ be the (unique) Haar measure such that $\mu\l(\Om\r)=1$, and let $\nu$ be the $G$-invariant Borel measure on $G/\G$ such that $\nu(G/\G)=1$, so that $\nu(A)=\mu(\pi^{-1}(A)\cap \Om)$ for every $A\subseteq G/\G$.

 We don't know if $\Om$ is compact, and we want to use compactness, so let $K\ss\Om$ be a compact set which contains most of $\Om$, i.e.\ such that $\mu\l(K\r)>\fr{1}{2}\mu\l(\Om\r)=\frac{1}{2}$. Such $K$ exists by regularity of $\mu$.

We claim $G=\l<K\cup \G\r>$ and hence $G=\l<K\cup S\r>$ for some finite $S$ which generates $\G$. Let $g\in G$; we will show $g\in K \G K^{-1}$. We have $$\nu\l(\pi\l(gK\r)\r)=\nu\l(g\pi\l(K\r)\r)=\nu\l(\pi(K)\r)=\mu\l(K\r)>\frac{1}{2}.$$
Since $\nu\l(G/\G\r)=1$, we get that $\pi\l(gK\r)\cap\pi\l(K\r)\ne \varnothing$ (otherwise we could deduce $\nu\l(G/\G\r)>1$).

Thus, $\l(gK\r)\G \cap K\G\ne\varnothing$ and therefore there exist $\g_1$, $\g_2$ and $k_1$, $k_2$ such that $gk_1\g_1=k_2\g_2$, and therefore $g=k_2\g_2\g_1^{-1}k_1^{-1}$, as we wished. 
\end{proof}
\begin{exe}
Show that the converse is true for uniform lattices. That is, show that, if $G$ is compactly generated and $\G\le G$ is a uniform lattice, then $\G$ is finitely generated.
\end{exe}
\begin{rem}
    This statement is not true for non-uniform lattices.
\end{rem}
\subsection{The Induced Representation}

Let $\Gamma$ be a lattice in a locally compact group $G$, and let
$\Omega$ be a fundamental domain of $\Gamma$.
Let $\rho:\G\to\mc{U}(\mc H)$ be a unitary representation of $\G$.
We seek to define an induced representation of $G$. We have two actions
of $\Gamma$: an action on $G$ by left multiplication, $\gamma.g=\gamma g$,
and an action of $\Gamma$ on $\mathcal{H}$ via $\rho$, $\gamma.v=\rho(\gamma)(v)$.
A map $\varphi:G\to\mathcal{H}$ is called \textbf{\emph{$\Gamma$-equivariant}}
(with respect to these two actions) if 
\[
\gamma.\varphi(g)=\varphi(\gamma g)
\]
for every $\gamma\in\Gamma$, $g\in G$ and $v\in\mathcal{H}$. The
induced unitary representation we will define for $G$ will not be
on $\mathcal{H}$, but on a new Hilbert space $\tilde{\mathcal{H}}$,
consisting of equivariant maps from $G$ to $\mathcal{H}$. Formally, we let
$X$ be the vector space of all measurable, $\Gamma$-equivariant
maps $f$ from $G$ to $\mathcal{H}$ for which 
\[
\int_{\Omega}\left\Vert f(\omega)\right\Vert^2 \mathrm{d}\mu(\omega)<\infty.
\]
We define an inner product on $X$ by setting
\[
\l<f_{1},f_{2}\r>=\int_{\Omega}\!\l<f_{1}\l(\omega\r),f_{2}\l(\omega\r)\r>\;\mr d\mu\l(\omega\r).
\]
This inner product is independent of our choice of the fundamental
domain $\Omega$ (since any two fundamental domains of $\Gamma$ have
the same measure). However, it is not strictly speaking an inner product,
since the norm of some vectors is zero. We therefore say that two
maps $f,g\in X$ are \emph{equivalent}, $f\sim g$, if $\left\Vert f-g\right\Vert =\left\langle f-g,f-g\right\rangle =0$.
We then set 
\[
\mc{\t H}=X/\!\sim\]
In order to keep notations reasonably simple, we will think of elements of $\mc {\t H}$ simply as maps from $G$ to $\mc H$, not as equivalence classes of such maps (just like we did with $L_2(X)$ for  measure spaces $X$). Keeping this in mind, and to summarise the above: 
\[
\mc{\t H}=\l\{f:G\to \mc{H}\text{ measurable}\mi\begin{array}{l}
{\displaystyle \int_{\Om}\l\|f\l(\omega\r)\r\|^2\;\mr d\mu(\omega)<\infty,}\\
f\l(\g g\r)=\rho\l(\g\r)f(g)\quad\forall\g\in\G,g\in G
\end{array}\r\}.
\]
The space $\mc{\t H}$
is a huge Hilbert space, and there is a natural $G$ action on it:
for every $g\in G$, $f\in\mc{\t H}$, $x\in G$, we define $\pi\l(g\r)$
by setting 
\[
g.f\l(x\r)=\l(\pi\l(g\r)\l(f\r)\r)\l(x\r)=f\l(x\cd g\r).
\]
Then $\pi$ is a unitary representation, and is called the \textbf{\emph{induced
representation of $\rho$}}. This is denoted by $\pi=\mr{ind}_{\G}^{G}\l(\rho\r)$.

\begin{exe} Verify $\t{\mc H}$ is a Hilbert space, and that $\pi$
is well-defined. \end{exe}

\subsection{There and Back Again}
\begin{thm}
\label{T:lat}
Let $G$ be a locally compact group, $\G$ a lattice. Then $G$ is Kazhdan if and only if $\G$ is Kazhdan.
\end{thm}

\begin{proof}[Proof of $\implies$]
Assume $G$ has property (T). We want to show $\G$ has property (T). Let $\rho:\G\to U\l(\mc H\r)$ be some unitary representation with almost invariant vectors; we need to show it has an invariant vector.

Let $\pi=\mr{ind}_\G^G\l(\rho\r)$ be the induced representation, $\t{\mc H}$ the corresponding function space. We want to prove $\pi$ admits almost invariant vectors as well. 

Let $\Om$ be a fundamental domain for $G/\G$, and let $\mu$ be the Haar measure such that $\mu\l(\Om\r)=1$. Take some Kazhdan pair $(K,\ep)$, where $K\subseteq G$ is compact and $\ep>0$. First we consider the case $\Om$ is compact and there exists some finite $S\ss \G$ such that $\Om K\ss S\Om$ (observe $\Om K$ is compact).

The representation $\rho:\G\to U\l(\mc H\r)$ admits almost invariant vectors, so there exists an $\l(S,\ep\r)$-invariant unit vector (with respect to $\rho$), call it $\hat v$.
Define $f\in\t{\mc H}$ by setting $f\l(\omega\r)=\hat v$ for every $\omega\in\Om$ and setting $f\l(\g\omega\r)=\rho\l(\g\r)\hat v$ for every $\g\in\G$, $\omega\in\Om$. If $\omega\in\Omega$ and $k\in K$, then $\omega k\in \Omega K\subseteq S\Omega$, so 
there are some $\omega'\in\Omega,s\in S$ such that $\omega k=s\omega'$.
Then $f(\omega k)=f(s\omega')=\rho(s)\hat{v}$ and $f(\omega)=\hat{v}$,
so 
\[
\left\Vert f(\omega k)-f(\omega)\right\Vert <\varepsilon.
\]
Therefore, for every $k\in K$, we have
$$\l\|k\cd f-f\r\|^2=\int_\Om \l\|f\l(\omega k\r)-f\l(\omega\r)\r\|^2\;\mr d\mu\l(\omega \r) < 
\int_\Om\ep^2\;\mr d\mu=\ep^2,$$ so $f\in\t{\mc H}$ is a $\l(K,\ep\r)$-invariant vector. 

In general, we can take $C\ss \Om$ compact and $V\supseteq \Om$ open such that $\mu\l(C\r),\mu(V)$ are as close to ${\mu\l(\Om\r)}=1$ as we wish, and we can take a finite $S\ss \G$ such that $CK\ss SV$, and hence $\mu\l(\Om K\setminus S\Om\r)$ is as small as we wish, so by choosing $\hat v$ which is $\t \ep$ invariant for appropriate $\t \ep>0$ (which depends on $\ep$), we can still construct almost invariant vectors in $\t{\mc H}$.

Therefore, in all cases, we can deduce by property (T) there exists a nonzero vector $\f\in\t{\mc H}$ which is $G$-invariant. As always, this implies that $\f$ is constant, so it is easy to show that its image is $\G$-invariant. But in fact, since $\f$ is really only defined up to null subsets, we need to make an extra step (if we want to be formal). Let $U$ be a subset of $G$ of finite, positive measure, and set $$v_1=\int_U\phi(u)\mr d\mu(u).$$ Then, for every $\g\in\G$, we have (by equivariance): $$\rho(\g)(v_1)=\int_U\rho(\g)(\f(u))\mr d\mu(u)=\int_U\f(\g u)\mr d\mu(u)=\int_{\g U}\f(u)\mr d\mu(u).$$ Since $\f$ is $G$-invariant, and since $\mu(\g U)=\mu(U)$, it is easy to show that $$\int _U\f \mr d\mu=\int_{\g U}\f\mr d\mu=v_1$$ (the reader may take this as an exercise), so $\rho(\g)(v_1)=v_1$ for all $\g\in\G$, as needed.
\end{proof}
\begin{proof}[Proof of $\neht$]
Assume $\G$ has property (T). We want to show $G$ has property (T) as well. Let $\pi:G\to U\l(\mc H\r)$ be a unitary representation with almost invariant vectors; we need to prove it admits an invariant vector.

$\G$ is finitely generated (because it has property (T)), so, by Lemma \ref{fin:com}, $G$ is compactly generated, $G=\l<K\r>$ for some compact $K\ss G$. Assume without loss of generality $K$ admits some $\ep>0$ and a finite subset $S$ such that $(S,\ep)$ is a Kazhdan pair. 

Let $\Om$ be some fundamental domain of $\G$, and let $\mu$ be the Haar measure such that $\mu\l(\Om\r)=1$. As usual, let $\nu$ be the $G$-invariant Borel measure on $G/\G$ such that $\nu(G/\G)=1$, so that $\nu(A)=\mu(p^{-1}(A)\cap \Om)$ for every $A$, where $p:G\to G/\G$ is the projection. By regularity, we can enlarge $K$ so that $\nu\l(p\l(K\r)\r)$ is as close to $1$ as we wish.

Let $\delta \in (0,1)$, and let $\hat v$ be a $\l(K,\ep\delta\r)$-invariant vector. In particular, it's $\l(S,\ep\delta\r)$-invariant, so there is a $\G$-invariant vector $\hat u$ such that $\l\|\hat v-\hat u\r\|<\delta$. Consider the map $G\to \mc H\!: g\mapsto \pi\l(g\r) \hat u$; it's constant on the cosets of $G/\G$, so it defines a map $G/\G\to \mc H\!:g\G\mapsto \pi\l(g\r)\hat u$. Therefore, we can push the $G$-invariant probability measure on $G/\G$ to the unit sphere in $\mc H$.

We know that $\hat v$ is $\l(K,\ep \delta\r)$-invariant and that $\l\|\hat v-\hat u\r\|<\delta$, so $\hat u$ is $\l(K,\t\ep\r)$-invariant for $\t\ep=\ep\delta+2\delta$ (or something like that). Therefore, since $\nu\l(p\l(K\r)\r)$ is as close to $1$ as we wish, the $G$-invariant probability measure on the unit sphere in $\mc H$ is as ``concentrated'' in a disk around $\hat u$ as we wish. Thus, the barycentre of the measure is not zero, and it's clearly fixed by $G$ (because the measure is $G$-invariant), as needed.\end{proof}
\begin{cor}\label{cor:FG}
If $G$ has property (T), then all its lattices are finitely generated  and have finite abelianisation.
\end{cor}

\begin{exa}
If $n$ is at least $3$, then $\mr{SL}_n\l(\mb R\r)$ has property (T) (although we haven't proved this yet), and $\mr{SL}_n\l(\mb Z\r)$ is a lattice in it (which we also haven't proved), so $\mr{SL}_n\l(\mb Z\r)$ has property (T).
\end{exa}

One of the main motivations of Kazhdan to define property $(T)$ was to prove that all lattices in classical groups are finitely generated. For a while there was no alternative proof of that result, but by now there is one which also provides a bound on the number of generators (see \cite{Gelander:rank}).

\section{Expanders}
A graph $\l(V,E\r)$ is called \emph{\textbf{$\ep$-expander}} if, for every finite subset $A\ss V$ with $\l|A\r|<\fr{1}{2}\l|V\r|$, necessarily
$$\fr{\l|\partial A\r|}{\l|A\r|}>\ep.$$ A graph is called an \emph{\textbf{expander}} if it's $\ep$-expander for some $\ep>0$. Every connected finite graph is an expander.

\begin{rem}
A Cayley graph of a group, with respect to some finite generating set, is an expander if and only if the group is non-amenable.
\end{rem}

A natural question to ask in this regard is the following: are there interesting families of connected, finite graphs which are $\ep$-expanders (for the same $\ep$) of bounded degree (i.e.\ with a bounded amount of edges at each vertex)? The first explicit construction of such a family was given by G. A. Margulis.

Denote 
$\G=\mr{SL}_3(\mb Z)$, $\G=\l<S\r>$ for some finite generating subset $S$. By property (T), there exists $\ep>0$ such that $\l(S,\ep\r)$ is a Kazhdan pair.

We have the canonical projection $\pi_n:\mr{SL}_3(\mb Z)\to \mr{SL}_3\l(\mb Z/n\mb Z\r)$. For all $n\in\mb N$, let $\G_n=\l\{A\in \mr{SL}_3(\mb Z)\mi\pi_n\l(A\r)=I_3\r\}=\ker \pi_n$.

Let $\rho_n$ be the regular representation of $\G$ on the finite dimensional Hilbert space $$\ell_2\l(\mr{SL}_3\l(\mb Z/n\mb Z\r)\r).$$ There are invariant vectors in this representation: the constant functions. Define $$\ell_2^o\l(\mr{SL}_3\l(\mb Z/n\mb Z\r)\r)=\l\{1\r\}^\perp=\l\{f\mi\int f=0\r\}.$$ By transitivity of the action of $\G$, only constant functions are invariant, so there are no nonzero invariant functions in $\ell_2^o$.

We claim $\mr{Cay}\l(\mr{SL}_3\l(\mb Z/n\mb Z\r),\pi_n\l(S\r)\r)\eqqcolon X_n$ are $\t \ep$-expanders for some $\t \ep$ (which depends on $\ep$).

Let $A\ss X_n$ such that $\l|A\r|\le \fr{1}{2}\l|X_n\r|$. Take $$v=1_A-\fr{\l|A\r|}{\l|X_n\r|}1_{X_n},\quad\hat v=\fr{v}{\l\|v\r\|}.$$ Since $\ep$ is a Kazhdan constant for $S$, and $v\in\ell_2^o$, there must be some $s\in S$ such that $$\l\|\rho_n\l(s\r)\hat v-\hat v\r\|\ge \ep,$$ and therefore $$\fr{\l|\partial A\r|}{\l|A\r|}>\t \ep$$ for some $\t \ep>0$ (which depends on $\ep$).
\begin{exe}
Find $\t \ep$ explicitly in terms of $\ep$.
\end{exe}

\section*{Further Reading}
The definitive reference for exploring property (T) in depth is \cite{bekka2008kazhdan}. This work provides a comprehensive analysis of the subject, making it an invaluable resource for delving into the intricacies of property (T). For further reading about expanders we refer to \cite{Lub}.

\chapter{Classical Groups}
So far, we haven't given a single example of a non-compact Kazhdan group. In this section, we show that $\mr{SL}_n(\mb R)$ has property (T) whenever $n\geqslant 3$.
First, we establish another property of $\mr{SL}_n(\mb R)$, which is a metric variant of the Howe--Moore property. It is very important and useful in its own right, and it holds for every $n\geqslant 2$.

\section{Metric Howe--Moore}
In this section we prove that in every isometric action of $\mr{SL}_n(\mb R)$, the stabiliser of any point is either compact or all of $\mr{SL}_n(\mb R)$. This is a metric version of the so called Howe--Moore property, and it is an extremely useful property.

\begin{thm}
Let $G=\mathrm{SL}_{n}(\mathbb{R})$ for some $n\geqslant2$. Every continuous action of $G$ on a metric space by isometries that does not have global fixed points is proper.
\end{thm}

An action of a locally compact group $G$ on a topological space is called \emph{proper }if the action map 
\begin{align*}
G\times X & \to X\times X\\
(g,x) & \mapsto(g.x,x)
\end{align*}
is a proper map, which means that the preimage of any compact subset
under it is compact. The following equivalent definition for proper
actions will be useful for the proof of the theorem.
\begin{exe}
Let $G$ be a locally compact group, $X$ a metrisable
space, $G\acts X$ a continuous action. Show that the action is \emph{not
}proper if and only if there exists $x,y\in X$ and a sequence $g_{n}\longrightarrow\infty$
in $G$ converging to infinity such that $g_{n}.x\longrightarrow y$.
\end{exe}

In the proof, we will use the Cartan decomposition of $\mr{SL}_n(\mb R)$, also known as the $KAK$-decomposition. 

We denote $K=\mr{SO}_n(\mb R)$ and \[
A=\left\{ \begin{pmatrix}a_{1} &  & 0\\
 & \ddots\\
0 &  & a_{n}
\end{pmatrix}\middle|a_{1}\geqslant\cdots\geqslant a_{n}>0,\prod_{i=1}^{n}a_{i}=1\right\} .
\]

\begin{lem}
    Every element $g\in \mr{SL}_n(\mb R)$ admits $k_1,k_2\in K$ and $a\in\mr{SL}_n(\mb R)\in A$ such that $g=k_1ak_2$.
\end{lem}
\begin{proof}
Let $g\in\mr{SL}_{n}(\mb R)$. Then $g^{T}g$ is a symmetric matrix,
hence orthogonally diagonalisable. Therefore, there is $k\in\mr{SO}_{n}(\mb R)$
such that $b\coloneqq kg^{T}gk^{-1}$ is diagonal. By multiplying $k$ by elementary orthogonal matrices, we may assume $b\in A$. Since $g^{T}g$
is a positive definite matrix, the elements on the diagonal of $b$
are positive, and we can take $a$ to be the diagonal matrix whose
entries are the square roots of the entries of $b$, so that $a^{2}=b$ and $a\in A$.
Then $p\coloneqq k^{-1}ak$ is a symmetric matrix in $\mr{SL}_{n}(\mb R)$
satisfying $p^{2}=g^{T}g$. Consider $o\coloneqq gp^{-1}$. We want
to show that $o\in\mathrm{SO}_{n}(\mathbb{R})$. We have $(p^{-1})^{T}=p^{-1}$
and $g^{T}gp^{-1}=p$, so 
\begin{align*}
\left\langle ov,ow\right\rangle  & =\left\langle v,(gp^{-1})^{T}gp^{-1}w\right\rangle 
 =\left\langle v,(p^{-1})^{T}g^{T}gp^{-1}w\right\rangle 
 =\left\langle v,w\right\rangle 
\end{align*}
for every $v,w\in\mathbb{R}^{n}$. Therefore $o\in\mathrm{SO}_{n}(\mathbb{R})$. Putting $k'=ok^{-1}$,
we get
$
g=op=k'ak,
$
which is a $KAK$-decomposition, as needed.\end{proof}
\begin{proof}[Proof of the theorem]
Let $G=\mathrm{SL}_{n}(\mathbb{R})$ for some $n\geqslant2$, let
$X$ be a metric space and let $G\acts X$ be a continuous action of
$G$ on $X$ by isometries. Assume by contradiction that the action of
$G$ on $X$ is not proper. 
By the exercise above, we know that there is $x\in X$ and a sequence
$g_{m}\longrightarrow\infty$ in $G$ converging to infinity such
that $g_{m}.x\longrightarrow y$ for some $y\in X$. 

We now use the $KAK$-decomposition. Write
$g_{m}=k_{m}a_{m}c_{m}$ for $a_{m}\in A$, $k_{m},c_{m}\in K$.
Since $K$ is compact, we may assume (up to passing to a subsequence)
that $k_{m}\longrightarrow k$ and $c_{m}\longrightarrow c$
for some $k,c\in K$, and hence $a_{m}\longrightarrow\infty$. 

Observe that $a_{m}c.x\longrightarrow k^{-1}y$. Indeed, by the triangle inequality,
\[
d(a_{m}c.x,k^{-1}y)\leqslant d(a_{m}c.x,a_{m}c_{m}.x)+d(a_{m}c_{m}.x,k_{m}^{-1}.y)+d(k_{m}^{-1}.y,k^{-1}y).
\]
The term equals $d(a_{m}c.x,a_{m}c_{m}.y)=d(c.x,c_{m}.x)$
and as $c_{m}\longrightarrow c$ it does to zero. The second goes to zero because it
is equal to $d(k_{m}a_{m}c_{m}.x,y)=d(g_{m}.x,y)$ and, by assumption,
$g_{m}.x\longrightarrow y$. The third goes to zero because $k_{m}^{-1}\longrightarrow k^{-1}$. 

Set 
$$x'=c.x,$$ 
so that $a_m.x'\to y'=k^{-1}y$. We will show that $x'$ is a global fixed point for $G$.

Towards that aim, we first claim that there exists a sequence $b_{m}\in\left\langle a_{i}\right\rangle _{i=1}^{\infty}\leqslant A$
such that $b_{m}\longrightarrow\infty$ and $b_{m}.x'\longrightarrow x'$.
 To see this, take a subsequence $a_{k_{m}}$ with $k_{m}$ going to infinity fast
enough so that $b_{m}\coloneqq a_{m}^{-1}a_{k_{m}}\longrightarrow\infty$.

Denote
\[
b_{m}=\mathrm{diag}(\beta_{m}^{1},\dots,\beta_{m}^{n}).
\]
Since $b_m\in A$, we have $b_1\geqslant\dots\geqslant b_n> 0$ and $\prod_{i=1}^n\beta_m^i=1$. By further replacing $(b_m)$ with a subsequence, we may assume that, for every $i=1,\dots,n$, either $\beta^i_m/\beta^{i+1}_m$ goes to infinity, or it is bounded.

Now, we define three subgroups: 
\begin{align*}
U_{+} & =\left\{ u\in G\middle|b_{m}^{-1}ub_{m}\longrightarrow1\right\} \leqslant G,\\
U_{-} & =\left\{ u\in G\middle|b_{m}ub_{m}^{-1}\longrightarrow1\right\} \leqslant G,\\
U_{0} & =\left\{ u\in G\middle|\left\{ b_{m}^{-1}ub_{m}\right\} _{m=1}^{\infty},\left\{ b_{m}ub_{m}^{-1}\right\} _{m=1}^{\infty}\text{ are pre-compact}\right\} \leqslant G.
\end{align*}

It is easy to check that these are all subgroups and that $U_{0}$
normalises both $U_{+}$ and $U_{-}$. 

Let us show that $\overline{\left\langle U_{+},U_{-},U_{0}\right\rangle }$
is all of $G$.
For $i\neq j$ and $\alpha\in\mathbb{R}$, we denote by $E_{ij}(\alpha)$ the matrix
with $1$'s on the diagonal, $\alpha$ on the $(i,j)$-entry, and
zero everywhere else. For every
two indices $i>j$, either  $\beta_m^i/\beta_m^j$ goes to infinity or it is bounded. If it goes to infinity, then a direct computation shows that $b_m^{-1}E_{ij}(\alpha)b_m$ goes to the identity element, so that $E_{ij}(\alpha)\in U_+$. If it is bounded, then $E_{ij}(\alpha)\in U_0$. Similarly, if $i<j$, we either have $E_{ij}(\alpha)\in U_-$ or $E_{ij}(\alpha)\in U_0$. It follows that $E_{ij}(\alpha)\in U_{+}\cup U_{-}\cup U_{0}$ for every $i\neq j$ and for every $\alpha$. Since $\left\{ E_{ij}(\alpha)\right\} _{i\neq j,\alpha\in\mathbb{R}}$ generates $G$, it follows that $G$ is generated by $U_{+}\cup U_{-}\cup U_{0}$.

Lastly, let us observe that $U_+$ is infinite: since $b_m\too \infty$, there is some $i_0$ such that $\beta^{i_0}_m/\beta^{i_0+1}_m\too \infty$, so (by the same considerations as above), we get that $E_{i_0,i_0+1}(\alpha)$ belongs to $U_+$. Similarly, $E_{i_0+1,i_0}\in U_-$. 

Since $U_{0}$
normalises both $U_{+}$ and $U_{-}$, $\overline{\left\langle U_{+},U_{-}\right\rangle }$
is normal in $\overline{\left\langle U_{+},U_{-},U_{0}\right\rangle}=G$. But this means that it is everything, since the
only normal subgroups of $G$ are finite.

We now use what is known as the `Mautner phenomenon'. Since $b_{m}.x'\longrightarrow x'$,
we get that $x'$ is fixed by $U_{+}\cup U_{-}$. This is in fact
very easy to see: if, for instance, $u\in U_{+}$, then, on the one
hand, $d(ub_{m}x',b_{m}x')\longrightarrow d(u.x',x')$ (since
$b_{m}x'\longrightarrow x'$). On
the other hand, $d(ub_{m}x',b_{m}x')=d(b_{m}^{-1}ub_{m}x',x')\longrightarrow d(1\cdot x',x')=0$.
Thus $d(u.x',x')=0$ so $x'$ is fixed by $u$.
The same argument works for $u\in U_{-}$, which means that $x'$
is fixed by all of $G$. 
\end{proof}
Whenever a group $G$ acts on a metric space $X$, we may consider the induced action of $G$ on $X\setminus F$, where $F=\l\{x\in X\middle|g.x=x\quad\forall g\in G\r\})\subseteq X$ is the subset of global fixed points. This is a subset of $X$ which is $G$-invariant,  so $G$ acts on it as well. By applying the theorem to this subset, we get the following corollary:
\begin{cor}
    Let $G=\mr{SL}_n(\mathbb{R})$ for some $n\geqslant 2$, and let $G\acts X$ be a continuous action by isometries on a metric space $X$. For every $x\in X$, the stabiliser group $G_x$ is either all of $G$ or compact. 
\end{cor}

\begin{rem}
A connected Lie group is called \textit{semisimple} if it admits no connected normal solvable subgroups apart from the trivial group.
    The theorem and corollary 
    are true for all semisimple Lie groups with finite centre. The proof is similar and not much harder.
\end{rem}
\section{Relative Property (T)}

It is very useful to consider a relative notion of property (T).

\begin{defn}
Let $G$ be a locally compact group, and let $H\leqslant G$ be a
closed subgroup. We say that $H$ has \emph{relative property (T)
in $G$}, or that it is a \emph{relatively Kazhdan subgroup of $G$},
if the following holds: whenever $\pi:G\to\mathcal{U}(\mathcal{H})$
is a unitary representation of $G$ with almost invariant vectors,
there is a nonzero vector $v\in\mathcal{H}\backslash\left\{ 0\right\} $
that is $\pi(H)$-invariant.
\end{defn}

\begin{rem}
Some authors say in such a case that the pair $(G,H)$ has property
(T).
\end{rem}

Our goal is to prove that $\mathbb{R}^{2}$ is relatively Kazhdan
inside $\mathrm{SL}_{2}(\mathbb{R})\ltimes\mathbb{R}^{2}$. For the
proof, we will need some tools from harmonic analysis,
which we presently recall in brief.
\subsection{A Bit of Harmonic Analysis}
If $A$ is a locally compact abelian group, we denote by $\hat{A}$
the \emph{Pontryagin dual} of $A$, which consists of the group of
all homomorphism from $A$ to the circle $\mb S^1$. It is also
a locally compact abelian group, with pointwise addition and the topology
of uniform convergence on compact subsets. 

If $\rho:A\to \mc U(\mc H)$ is a unitary representation of a locally compact abelian group $A$, then there is a \textit{projection valued measure} $E$ defined on the dual group $\hat A$. By `projection valued', we mean that it assigns to each Borel subset
$B$ of $\hat A$, not a number, but a projection operator
on $\mathcal{H}$, denoted by $E(B)$. Accordingly, the integral $\int_{\chi\in\hat A}f(\chi)dE(\chi)$
of a function $f:\hat A\to\mathbb{C}$ is not a
number, but an operator on $\mathcal{H}$. This projection valued
measure $E$ satisfies the following formula, for every $a\in A$:
\begin{equation}
    \label{eqn:int}
    \rho(a)=\int_{\chi\in\hat A}\chi (a)\mathrm{d}E(\chi).
\end{equation}

For every $\xi\in\mathcal{H}$, the projection valued measure $E$
induces an honest-to-god measure on $\hat A$, $E_{\xi}$, defined by
\[
E_{\xi}(B)=\left\langle E(B)(\xi),\xi\right\rangle 
\]
for every Borel $B\subseteq\hat A$. If $\xi$ is
a unit vector, this is a probability measure. We get the following
formula:
\[
\left\langle \rho(a)(\xi),\xi\right\rangle =\int_{\chi\in\hat A}\chi(a)\mathrm{d}E_{\xi}(\chi).
\]

In considering (relative) property (T), we are naturally interested in the following question: When does $\rho: A\to \mc U(\mc H)$ admit a non-zero fixed vector? Let's consider first the simplest example: $A=\mb Z$ and $\mc H$ being finite-dimensional. Then $\rho$ is really just a choice of a unitary operator $U$, and Equation \eqref{eqn:int} becomes
\[
    U=\int_{z\in \mb S^1}z\  \mr{d}E(\{z\}),
\]
where, for any complex number $z\in\mb S^1$, $E(\{z\})$ is the projection onto the eigenspace of $U$ corresponding to $z$. Since there are only finitely many eigenvalues, $E$ is atomic and finitely supported, and the integral is just a sum, so the equation is exactly the classical spectral theorem. In this case, $\rho$ has a non-zero fixed vector if and only if $1$ is an eigenvalue, i.e., if and only if $E$ has an atom at $1$. 

If we take $A=\mb Z^n$, $\rho$ is a choice of $n$ unitary operators which all commute with one another. Since they commute, they preserve each other's eigenspaces and one can decompose the space to a direct sum of spaces, each correspond to a different character. This gives the projection valued measure in that case. The eigenspace of the trivial character is the space of invariant vectors.

The general case is similar. Any representation $\rho:A\to \mc U(\mc H)$ of a locally compact abelian group $A$ has a non-zero fixed vector if and only if $E$ has an atom at $1\in\hat A$, where $1:A\to\mb S^1$ is the trivial homomorphism. We leave it as an exercise to the reader to deduce this from Equation (\ref{eqn:int}).

We will use this theory in the proof that $\mb R^2$ is relatively Kazhdan inside $\mr{SL}_2(\mb R^2)\ltimes \mb R^2$. So, in our case, $A=\mb R^2$. We have an action of $S=\mr{SL}_2(\mb R)$ on $A$ given by matrix multiplication, where we think of the elements of $A$ as column vectors. We have an induced action of $S$ on $\hat A$, given by $g.\xi (a)=\xi(g^{-1}.a)$. One may think of the elements of $\hat A\cong \mb R^2$ as row vectors, so that
\[
    (\lambda,\mu).\begin{pmatrix}x\\
y
\end{pmatrix}=\exp(2\pi i(\lambda x+\mu y)).
\]
Then the action of $S$ on $\hat A$ is given by matrix multiplication as well, $g.(\lambda,\mu)=(\lambda,\mu)\cdot g^{-1}$. We have the following explicit formulation of Equation \eqref{eqn:int}:
\[
    \rho(x,y)=\int_{(\lambda,\mu)\in\widehat{\mathbb{R}^{2}}}\exp\left(2\pi i(\lambda x+\mu y)\right)\mathrm{d}E(\lambda,\mu),
\]
for every $(x,y)\in\mb R^2$. The projection valued measure $E$ is equivariant, in the sense that
\[
E(g.B)=\pi(g^{-1})E(B)\pi(g).
\]
For every $g\in\mathrm{SL}_{2}(\mathbb{R})$ and every unit vector
$\xi\in\mathcal{H}$, the probability measure $E_{g.\xi}$ is exactly
$g^{-1}_{*}(E_{\xi})$, the pushforward of $E_{\xi}$ by $g^{-1}$.

\subsection{The Main Theorem}
\begin{thm}
\label{prop:SL2R-R^2}The normal subgroup $\mathbb{R}^{2}$ of $\mathrm{SL}_{2}(\mathbb{R})\ltimes\mathbb{R}^{2}$
is relatively Kazhdan.
\end{thm}

\begin{proof}
Let $\pi:G\to\mathcal{U}(\mathcal{H})$ be a unitary representation
of $G=\mathrm{SL}_{2}(\mathbb{R})\ltimes\mathbb{R}^{2}$.
Denote $\rho=\pi\restriction_{\mb R^2}$,
and let $E$ be the projection valued measure on $\widehat{\mathbb{R}^{2}}$ satisfying
\[
    \rho(x,y)=\int_{(\lambda,\mu)\in\widehat{\mathbb{R}^{2}}}\exp\left(2\pi i(\lambda x+\mu y)\right)\mathrm{d}E(\lambda,\mu).
\]
Now, let $\xi_{n}\in\mathcal{H}$ be an asymptotically invariant sequence
of unit vectors. The probability measures $E_{\xi_n}$ are asymptotically invariant.
This can be seen by the following computation: for every unit $\xi\in\mathcal{H}$,
every Borel $B\subseteq\widehat{\mathbb{R}^{2}}$ and every $g\in\mathrm{SL}_{2}(\mathbb{R})$,
we get that
\[
\left|E_{\xi}(B)-E_{g.\xi}(B)\right|=\left|\left\langle E(B)(\xi),\xi\right\rangle -\left\langle E(B)\pi(g)(\xi),\pi(g)(\xi)\right\rangle \right|.
\]
By adding and subtracting $\left\langle E(B)(\xi),\pi(g)(\xi)\right\rangle $,
we get that this is equal to
\[
\left|\left\langle E(B)(\xi),\xi-\pi(g)(\xi)\right\rangle +\left\langle E(B)(\xi-\pi(g)(\xi)),\pi(g)(\xi)\right\rangle \right|.
\]
By the Cauchy--Schwartz inequality (and the fact $\left\Vert E(B)(\xi)\right\Vert \leqslant1$),
we get that 
\begin{align*}
\left|\left\langle E(B)(\xi),\xi-\pi(g)(\xi)\right\rangle \right| & \leqslant\left\Vert \xi-g.\xi\right\Vert ,\\
\left|\left\langle E(B)(\xi-\pi(g)(\xi)),\pi(g)(\xi)\right\rangle \right| & \leqslant\left\Vert \xi-g.\xi\right\Vert .
\end{align*}
Therefore, by the triangle inequality, we get that
\[
\left|E_{\xi}(B)-E_{g.\xi}(B)\right|\leqslant2\left\Vert \xi-g.\xi\right\Vert .
\]

Assume by contradiction that no non-zero vector of
$\mathcal{H}$ is fixed by $\mathbb{R}^{2}$. This means that $E(\left\{ (0,0)\right\} )=0$.
Therefore, the probability measures $E_{\xi_{n}}$ are probability
measures on $\widehat{\mathbb{R}^{2}}\backslash\left\{ (0,0)\right\} $.
These probability measures are asymptotically invariant. If they admitted
a limit point, it ought to have been invariant, and we would get an
$\mathrm{SL}_{2}(\mathbb{R})$-invariant probability measure on $\widehat{\mathbb{R}^{2}}$.
Since $\widehat{\mathbb{R}^{2}}\backslash\left\{ (0,0)\right\} $
is not compact, there is no reason to think that $E_{\xi_{n}}$ admit
a limit point. What we can do, however, is to push-forward these measures
to probability measures on the projective line $\mathbb{P}(\widehat{\mathbb{R}^{2}})$.
This is possible since they are defined on $\widehat{\mathbb{R}^{2}}\backslash\left\{ 0,0\right\} $.
Since the projective line is compact, so is the space of probability measures on it with respect to the weak-* topology. Then there is necessarily a
limit probability measure, which is necessarily invariant under the induced action
of $\mathrm{SL}_{2}(\mathbb{R})$. This is impossible by the next lemma.
\end{proof}
\begin{lem}
There is no probability
measure on $\mathbb{P}^{1}$ that is invariant under the natural action of $\mathrm{SL}_{2}(\mathbb{R})$.
\end{lem}
\begin{proof}
Assume by contradiction that there is such a measure $\mu$ on $\mathbb{P}^{1}$.
Then, in particular, $\mu$ is invariant under the action of $\mathrm{SO}_{2}(\mathbb{R})$.
This action is transitive, and the stabiliser of a point is $\left\{ \pm1\right\} $,
so $\mathbb{P}^{1}\cong\mathrm{SO}_{2}(\mathbb{R})/\left\{ \pm1\right\}$ and the unique invariant probability measure is the normalized Haar (or Lebesgue) measure.
However the Lebesgue measure on the projective line is not invariant under the diagonal matrix $\text{diag}(2,1/2)$.
\end{proof}

\section{Property (T) of \texorpdfstring{$\mathrm{SL}_{n}(\mathbb{R})$}{SLnR}}
We are now in a position to prove Kazhdan's celebrated theorem:

\begin{thm}[Kazhdan]
    For every $n\geqslant 3$, the group $\mr{SL}_n(\mathbb{R})$ has property (T).
\end{thm}
\begin{proof}
It is easy to see that the group $\mathrm{SL}_{2}(\mathbb{R})\ltimes\mathbb{R}^{2}$
is isomorphic to a closed subgroup of $\mathrm{SL}_{n}(\mathbb{R})$
for every $n\geqslant3$. For example, in $\mathrm{SL}_{3}(\mathbb{R})$,
the following subgroup works:
\[
\left\{ \begin{pmatrix}a & b & x\\
c & d & y\\
0 & 0 & 1
\end{pmatrix}\middle|\begin{pmatrix}a & b\\
c & d
\end{pmatrix}\in\mathrm{SL}_{2}(\mathbb{R}),\begin{pmatrix}x\\
y
\end{pmatrix}\in\mathbb{R}^{2}\right\} .
\]
A similar embedding works for $\mathrm{SL}_{n}(\mathbb{R})$ for every
$n\geqslant3$. Therefore, if $\pi:\mathrm{SL}_{n}(\mathbb{R})\to\mathcal{U}(\mathcal{H})$
is a unitary representation with almost invariant vectors, then the
restriction of this representation to the subgroup $\mathrm{SL}_{2}(\mathbb{R})\ltimes\mathbb{R}^{2}$
has almost invariant vectors as well, which means that the action
admits a nonzero vector $v\in\mathcal{H}\backslash\left\{ 0\right\} $
fixed by $\mathbb{R}^{2}$, which is a non-compact subgroup of $\mathrm{SL}_{n}(\mathbb{R})$.
Considering the action of $\mathrm{SL}_{n}(\mathbb{R})$ on $\mathcal{H}\backslash\left\{ 0\right\} $,
which is a continuous action by isometries on a metric space, we get
that the stabiliser of $v$ is non-compact. Therefore, by the Howe--Moore property, the stabiliser is all of $\mathrm{SL}_{n}(\mathbb{R})$.
That is, $v$ is a fixed vector.
\end{proof}
\begin{rem}
A Lie group is called almost simple if it has finite centre and modulo the centre it is a simple group. 
The rank of an almost simple Lie group is the maximal dimension of a subgroup which is diagonalizable in the adjoint representation. For instance the rank of $\mathrm{SL}_{n}(\mathbb{R})$ is $n-1$. Kazhdan's theorem holds for all almost simple Lie groups of rank at least two, and the proof of the general case is similar, and in particular relies on relative property (T).
It was later shown by other means that some families of rank one groups also possess property (T). This however doesn't apply of course for $\mathrm{SL}_{2}(\mathbb{R})$ (and more generally for groups of isometries of real or complex hyperbolic spaces).
\end{rem}
\section*{Further Reading}
A similar proof of the relative property of $\mb R^2$ in $\mr{SL}_2(\mb R)\ltimes \mb R^2$ can be found in \cite{bekka2008kazhdan}. The proof of the metric Howe--Moore property was established in \cite{bader2017equicontinuous}, in which it is proven in a much greater generality.

\chapter{Property FH}
\section {Isometric Actions}
In this section, we will only deal with real Hilbert spaces, so unitary transformations and representations are now orthogonal transformations and representations, respectively.

Recall that an isometry is a map between metric spaces which preserves distances. 
\begin{lem}\label{iso-is-affine}
        Let $\mc H$ be a real Hilbert space, and let $f:\mc H\to \mc H$ be an isometry (i.e., $\|f(x)-f(y)\|=\|x-y\|$ for every $x,y\in \mc H$). Then there is some orthogonal transformation $A$ of $\mc H$ and some vector $v\in \mc H$ such that $f(x)=Ax+v$ for ever $x\in\mc H$.
\end{lem}
\begin{exe}
        Prove this lemma.
\end{exe}
\begin{hint}
        Isometries send midpoints to midpoints.
\end{hint}
\begin{rem}
        This is a special case of the Mazur-Ulam theorem.
\end{rem}

Therefore, we have $$\mr{Isom}\l(\mc H\r)=O\l(\mc H\r)\ltimes \mc H.$$

Suppose a group $G$ acts by isometries on a Hilbert space $\mc H$,
$$G\x \mc H\to \mc H\quad \l(g,v\r)\mapsto g.v.$$
In the case $G$ is a topological group, we will require this map to be continuous as well. Such an action can also be viewed as a homomorphism $$G\to \mr{Isom}\l(\mc H\r)$$ sending each $g$ to the isometry defined by $v\mapsto g.v$. By the last lemma, for all $g\in G$ we have 
$$g.v=\pi\l(g\r)(v)+c\l(g\r),$$
for some $\pi:G\to O\l(\mc H\r)$, $c:G\to\mc H$. $\pi$ is a homomorphism, and $c$ is said to be a \textit{cocycle}. If $g,h\in G$, then
\begin{align*}
\l(gh\r).v&=g.\l(h.v\r) = \pi(g)\l(\pi\l(h\r) (v)+c\l(h\r)\r)+c\l(g\r)\\
&=\pi\l(g\r)\l(\pi\l(h\r)(v)\r) + \pi\l(g\r) (c\l(h\r))+c\l(g\r)=\pi\l(gh\r)(v)+\pi\l(g\r)(c\l(h\r))+c\l(g\r),
\end{align*}
while also
$$(gh).v=\pi(gh)(v)+c(gh),$$
so
$$c(gh)=\pi(g)(c(h)) + c(g).$$
This is called the \textit{cocycle equation}.

Alternately, given an orthogonal representation $\pi:G\to O(\mc H)$, a map $c:G\to\mc H$ is a \emph{\textbf{$\pi$-cocycle}} (or just a \emph{\textbf{cocycle}} if $\pi$ is clear) if it satisfies the cocycle equation relative to $\pi$, i.e.\ if
$$c(gh)=\pi(g)(c(h))+c(g)\quad\forall g,h\in G.$$
\begin{exe}
        Given an orthogonal representation $\pi$ of a group $G$, show every $\pi$-cocycle defines an action of $G$ on the underlying Hilbert space. Thus, given such $\pi$, there is a correspondence of cocycles and isometric actions.
\end{exe}
\begin{exe}
Show that an affine action $G\acts \mc H$, $g.v=\pi(g)(v)+c(g)$ is continuous if and only if both $\pi:G\to \mc U(\mc H)$ and $c:G\to \mc H$ are continuous. 
\end{exe}

If $\pi$ is the trivial representation, i.e.\ if $\pi(g)(v)=v$ for all $g\in G$, $v\in\mc H$, then a $\pi$-cocycle is just a homomorphism $G\to \mc H$.

We consider the following question: when does an action admit a global fixed point? That is, a vector $v\in\mc H$ such that $$g.v = \pi(g)(v)+c(g)=v\quad\forall g\in G,$$
or equivalently $$c(g)=v-\pi(g)v.$$
\begin{exe}
        Given an orthogonal representation $\pi$ and some $v\in\mc H$, show that the map $g\mapsto v-\pi(g)v$ is a cocycle, and that the corresponding action has a global fixed point.
\end{exe}

A cocycle of the form $c(g)=v-\pi(g)v$ for some $v\in\mc H$ is called \emph{\textbf{trivial}} (or a \emph{\textbf{co-boundary}}).

\section{Property FH}
\begin{defn}
        A topological group $G$ is said to have \emph{\textbf{property FH}} if every continuous isometric action of $G$ on a Hilbert space $\mc H$ admits a global fixed point.
\end{defn}
A topological group $G$ has property FH if and only if for any orthogonal representation $\pi$, every $\pi$-cocycle is trivial.
\begin{exe}
        Every compact group has property FH.
\end{exe}

\subsection{Actions on Trees}
The property we will define next is related to trees. A natural metric on trees is the one assigning two points the length of the (unique) path between them (which never goes through the same edge twice). With this metric, one can consider isometric actions on trees as well.
\begin{defn}
        A group is said to have \emph{\textbf{property FA}} if it admits a global fixed point whenever it acts (continuously and isometrically) on a tree.
\end{defn}

The group $\mr{Aut}(T)$ of automorphisms of a given (locally finite) tree $T$ is a topological group, which is locally compact and totally disconnected (see Example \ref{autT}). It obviously does not have property FA. An action on a tree $T$ is the same thing as a homomorphism into $\mr{Aut}(T)$.
\begin{thm}\label{thm:FA-eq}
        A discrete group $\G$ has property FA if and only if it satisfies the following:
        \begin{enumerate}[1)]
                \i $\G$ is finitely generated.
                \i There are no $A,B,C$ with $\l[A:C\r],\l[B:C\r]\ge 2$ such that $\G = A*_C B$.
                \i $\mb Z$ is not a homomorphic image of $\G$.
        \end{enumerate}
\end{thm}

Property FH is stable under taking finite index subgroups (i.e.\ if $G$ has property FH and $H\le G$ is of finite index then $H$ has property FH). On the other hand, there are examples of groups with property FA with subgroups of finite index without property FA. 

\begin{lem}
\label{fix:orb}
        If $G\acts X$ (continuously and isometrically), where $X$ is either a tree or a Hilbert space (or more generally, any uniformly convex complete metric space) then the following are equivalent:
        \begin{enumerate}[1)]
                \i There is a global  fixed point.
                \i There is a bounded orbit.
                \i All orbits are bounded.
        \end{enumerate}
\end{lem}
A \emph{\textbf{geodesic metric space}} is a metric space $X$ in which for every $x,y$ there is a path $\g:\l[a,b\r]\to X$, $\g(a)=x$, $\g(b)=y$, such that the length of $\g$, defined by
$$L\l(\g\r)=\sup\l\{\sum_{i=1}^nd\l(\g(t_{i-1}),\g(t_i)\r)\mi a=t_0<t_1<\cdots<t_n=b,n\in \mb N\r\},$$
equals $d\l(x,y\r)$. Hilbert spaces and trees (where we consider the edges as copies of the unit interval, and a part of our space) are geodesic metric spaces. In geodesic metric spaces, there is a midpoint between any two points $x,y$; namely, $\g(\fr{a+b}{2})$ for a path $\g$ between them parametrised by arc length. The space is called \textit{\textbf{uniquely geodesic}} if the geodesic $\g$ is unique up to reparametrisation. In this case, the midpoint is unique as well, and is denoted by $m(x,y)$.

A geodesic metric space $X$ is \emph{\textbf{uniformly convex}} if for every $\ep>0$ and $r>0$ there is some $\delta>0$ such that, for all $x,y,z\in X$ satisfying
\begin{align*}
d\l(x,y\r)&\le r\\
d\l(x,z\r)&\le r\\
d\l(y,z\r)&>\ep,
\end{align*}
it also holds that
$$d\l(x,m\r)\le r-\delta$$
for every midpoint $m$ between $y$ and $z$. Hilbert spaces and trees (when considered as above) are uniformly convex.
\begin{proof}[Proof of Lemma \ref{fix:orb}]
Given any bounded set $A$ in $X$, define the circum-radius of $A$ by
$$r\l(A\r)=\inf\l\{r^\pr\mi A\text{ is contained in some closed ball of radius }r^\pr\r\}.$$
\begin{lem}\label{circum}
Let $X$ be a complete, uniformly convex metric space. Let $A\ss X$ be a bounded subset. There is a unique $c\l(A\r)\in X$, called the circum-centre, such that $$\bar B\l(c\l(A\r),r\l(A\r)\r)\supseteq A.$$
\end{lem}
\begin{cor}
The circum-centre of a bounded, $G$-invariant set must be a fixed point.
\end{cor}
The corollary immediately follows from Lemma \ref{circum}, because $g.c\l(A\r)=c\l(g.A\r)$. Since orbits are $G$-invariant, Lemma \ref{fix:orb} immediately follows from the corollary. So, it only remains to prove Lemma \ref{circum}. \end{proof}
\begin{proof}[Proof of Lemma \ref{circum}]
First, we show that there can be at most one circum-centre. Assume by contradiction $p_1,p_2$ are  circum-centres with $d\l(p_1,p_2\r)>\ep$ for some $\ep>0$. Write $r=r(A)$. Then $B\l(p_1,r\r),B\l(p_2,r\r)\supseteq A$. By uniform convexity, there is a $\delta>0$ such that for all $x,y,z\in X$ satisfying
\begin{align*}
d\l(x,y\r)&\le r\\
d\l(x,z\r)&\le r\\
d\l(y,z\r)&>\ep,
\end{align*}
it also holds
$$d\l(x,m\r)\le r-\delta$$
for any midpoint $m$ between $y$ and $z$. For every $x\in A$, $\l(x,p_1,p_2\r)$ is such a triple, so, for every $x\in A$ and every midpoint $m_p$ between $p_1$ and $p_2$, we have
$$d\l(x,m_p\r)\le r-\delta,$$
i.e.\ $A$ is contained in the ball of radius $r-\fr{\delta}{2}$ around $m_p$, which contradicts the definition of $r$.

Now, let us prove existence. Take $r_n\coloneqq r+\fr{1}{n}\searrow r$, and take $x_n$ such that
$$B\l(x_n,r_n\r)\supseteq A.$$
Then $x_n$ is a Cauchy sequence (one can see this similarly to the above), and since $X$ is complete it converges to some point, which is a circum-centre.
\end{proof}

\begin{cor}
        If $\l(\pi,c\r)$ is an isometric action of $G$ on a Hilbert space $\mc H$ (where $\pi$ is the orthogonal part and $c$ the cocycle), then the following are equivalent:
        \begin{itemize}
                \i $c$ is trivial (a co-boundary).
                \i $c$ is bounded.
        \end{itemize}
\end{cor}
\begin{proof}
        $c$ is just the orbit of $0$.
\end{proof}
\begin{thm}
                \label{FH=>FA}
        Property FH implies property FA.
\end{thm}
The converse of this theorem does not hold.
\begin{proof}
        Let $G$ be a group without property FA. Then there is some tree $T$ and a (continuous isometric action) $G\acts T$ without a global fixed point. By Lemma \ref{fix:orb}, this means all the orbits of this action are unbounded. We will use this action to construct an action on a Hilbert space with unbounded orbits, thus proving $G$ doesn't have property FH.
        
        We will consider $T$ as a directed graph, where for every edge $e=\l(v_1,v_2\r)$ the edges in $T$ also contain its reversed edge $\bar e=\l(v_2,v_1\r)$. Denote $E=\l\{\text{all directed edges in }T\r\}$, and fix some $O$ as the root of $T$. Recall the orbit of $O$ is unbounded.
        
        Consider the Hilbert space $\mc H\coloneqq \ell_2(E)$. Let $\pi$ be the orthogonal representation of $G$ associated with $G\acts T$, i.e.
        $$\l(\pi\l(g\r) f\r)\l(e\r)=f\l(g^{-1} e\r).$$
        Define $c:G\to\mc H$ by 
        $$c\l(g\r)=1_{\l[O,g\cd O\r]}-1_{\l[g\cd O,O\r]},$$
        where $\l[a,b\r]$ is (the set of edges in) the unique path in $T$ from $a$ to $b$ (which doesn't go through the same edge twice). We claim $c$ is an unbounded cocycle. It's unbounded since $\l\{d\l(O,g\cd O\r)\r\}$ is unbounded (because it's the orbit of $O$) and if $d\l(O,g\cd O\r)=n$ then $\l\|c\l(g\r)\r\|=\sqrt {2n}$. In order to show it's a cocycle, we need to show
        $$c(gh)=\pi(g)c(h)+c(g).$$
        Observe we have $\pi(g)\cd 1_{\l[a,b\r]}=1_{\l[g\cd a,g\cd b\r]}$ for all $g\in G$ and $a,b\in V$ (i.e.\ $\pi(g)$ ``shifts'' paths) and $1_{\l[a,b\r]}+1_{\l[b,c\r]}=1_{\l[a,c\r]}$ (i.e.\ the sum of the paths from $a$ to $b$ and from $b$ to $c$ is just the path from $a$ to $c$), so
        $$\pi(g)c(h)=\pi(g)\cd 1_{\l[O,h\cd O\r]}-\pi(g)\cd 1_{\l[h\cd O,O\r]}=1_{\l[g\cd O,gh\cd O\r]}-1_{\l[gh\cd O,g\cd O\r]},$$
        and
        \begin{align*}
        \pi(g)c(h)+c(g)&=1_{\l[g\cd O,gh\cd O\r]}-1_{\l[gh\cd O,g\cd O\r]}+1_{\l[O,g\cd O\r]}-1_{\l[g\cd O,O\r]}\\
        &=1_{\l[O,g\cd O\r]}+1_{\l[g\cd O,gh\cd O\r]}-\l(1_{\l[gh\cd O,g\cd O\r]}+1_{\l[g\cd O,O\r]}\r)\\
        &=1_{\l[O,gh\cd O\r]}-1_{\l[gh\cd O,O\r]}=c(gh),
        \end{align*}
        as needed. More visually, $c(h)$ is the path from $O$ to $hO$, so $\pi(g)c(h)$ is the path from $gO$ to $ghO$. $c(g)$ is the path from $O$ to $gO$, so the sum $c(g)+\pi(g)c(h)$ is the sum of the path from $O$ to $gO$ and the path from $gO$ to $ghO$, i.e.\ the path from $O$ to $ghO$, which is exactly $c(gh)$.
        
        In fact, the attentive reader will have noticed we cheated a bit, since there might be some overlap between the geodesics. That is, if $a,b,c$ are vertices, the intersection $[a,b]\cap[b,c]$ might be nonempty. However, there is always some vertex $q$ such that $$1_{[a,b]}+1_{[b,c]}=1_{[a,c]}+1_{[q,b]}+1_{[b,q]},$$where it might be that $q=b$ (in which case we understand $[b,b]$ to be the empty set, so that $1_{[b,b]}=0$). The reader may take it as an exercise to show that these terms cancel out in the computation above.
\end{proof}

\section{Back to Property (T)}
It turns out property FH is almost the same as property (T). 

\begin{thm}
\label{FH:T}
        If $G$ is a locally compact $\sigma$-compact group, then $G$ has property (T) if and only if it has property FH.
\end{thm}
\begin{rem}
Recall that we always assume our locally compact groups are second-countable, and, in particular, $\sigma$-compact. We emphasised this condition in the statement of the theorem since it is actually needed: there is a locally compact (non-$\sigma$-compact) group with property FH without property (T). The other direction of the theorem is true for every locally compact group, even if it's not $\sigma$-compact or second-countable.
\end{rem}
For the proof of (FH)$\implies$(T), we need the open mapping theorem. In its general formulation (which we will use) it concerns \textit{Fr\`echet spaces}. A \emph{\textbf{Fr\`echet space}} is a topological vector space $X$ which satisfies the following:
\begin{enumerate}[1)]
        \i $X$ is Hausdorff.
        \i The topology of $X$ is induced by a countable family of semi-norms $\l\|\cd\r\|_n$, $n\in N$, $\l|N\r|\le \aleph_0$ (i.e.\ the collection $\bigcup _{n\in N}\mc U_n$ is a subbasis, where $\mc U_n$ is the collection of balls in the semi-norm $\l\|\cd\r\|_n$).
        
        A sequence or a net $\l(x_n\r)$ converges in such a topology if and only if it converges with respect to all of the semi-norms.
        \i $X$ is complete with respect to the above family of semi-norms.
\end{enumerate}The theorem says that if $X$, $Y$ are Fr\`echet spaces then any continuous surjective linear map $T:X\to Y$ is an open map.
\begin{proof}[Proof of (FH) $\implies$ (T)]\
        Assume a group $G$ has property FH. We want to show it has property (T). Recall $G$ has property (T) if every orthogonal representation $\pi$ with almost invariant vectors (i.e.\ such that for all compact $K$ and all $\ep>0$ there is a $\l(K,\ep\r)$-invariant vector) admits a globally invariant vector. Thus, let $\pi$ be an orthogonal representation with almost invariant vectors, and assume by contradiction it doesn't admit a globally invariant vector.
        
        Define $C(\pi)$ to be the set of all $\pi$-cocycles $c:G\to \mc H$. For each compact $K\ss G$ and $c\in C(\pi)$, set 
        $$\Vert c\Vert_K=\mr{max}\{c(k)|k\in K\}.$$
        This is a semi-norm on $C(\pi)$, so $C(\pi)$ with the induced topology is a Fr\`echet space.\footnote{We are using here $\sigma$-compactness: the countable family of seminorms induced by the countable family of compact subsets covering the space induce the topology.} We leave it to you to check this. 
        
        Consider the map $T:\mc H\to C(\pi)$, defined by
        $$\l(Tv\r)\l(g\r)=v-\pi\l(g\r)\cd v.$$
        That is, for every $v$, $Tv$ is the $\pi$-cocycle defined by $g\mapsto v-\pi\l(g\r)v$. $T$ is continuous and linear, and since $G$ has property FH, all cocycles are trivial and $T$ is onto. We assumed $\pi$ has no invariant vectors, so it's injective. Therefore, by the open mapping theorem, $T^{-1}$ is continuous as well.
        
        $\pi$ admits almost invariant vectors, so for all compact $K\ss G$ and for all $\ep>0$ there is some unit vector $v_{K,\ep}$ such that
        $$\l\|v_{K,\ep}-\pi\l(g\r)v_{K,\ep}\r\|\le \ep\quad\forall g\in K,$$
        i.e.\ such that
        $$\l\|Tv_{K,\ep}\r\|_K\le\ep.$$
        We can define an order on $\l\{\l(K,\ep\r)\r\}$ by setting $\l(K_1,\ep_1\r)\ge\l(K_2,\ep_2\r)$ if and only if $K_1\supseteq K_2$ and $\ep_1\le \ep_2$. Then $\l(v_{K,\ep}\r)$ and $\l(Tv_{K,\ep}\r)$ are nets, and we just showed
        $$\l\|Tv_{K,\ep}\r\|_C\too 0$$
        for all compact $C\ss G$, i.e.\ that
        $$Tv_{K,\ep}\too 0$$
        in the topology of $C(\pi)$. By applying $T^{-1}$ on both sides we get
        $$v_{K,\ep}\too T^{-1}0\eqqcolon v_0,$$
        and then $v_0$ is an invariant vector, a contradiction.
\end{proof}
For the proof of (T)$\implies$(FH) we need a few more constructions.

Let $X$ be a set. A map $K:X\x X\to\mb C$ is a \emph{\textbf{positive definite kernel}} if for all $x_1,\dots,x_n\in X$ and for all $c_1,\dots,c_n\in \mb C$ we have
$$\sum_{i,j}c_i\bar c_jK\l(x_i,x_j\r)\ge 0.$$
(In particular, this sum needs to be real). This implies $K\l(x,y\r)=\overline{K\l(y,x\r)}$ for all $x,y\in X$. We say a kernel is \emph{\textbf{strictly positive definite}} if the sum it strictly greater than zero (unless $c_i=0$ for all $i$).

This is equivalent to $K\l(x_i,x_j\r)$ being a positive-definite matrix or a strictly positive definite matrix for all $x_1,\dots,x_n\in X$ respectively.

If $G$ is group, a function $\psi:G\to\mb C$ is \emph{\textbf{positive definite}} if $K\l(x,y\r)\coloneqq \psi\l(x^{-1}y\r)$ is a positive definite kernel.

In the case $G$ is topological or $X$ is a topological space, we also require the maps to be continuous.
\begin{exa}
        Let $\mc H$ be a Hilbert space. Then $K\l(x,y\r)\coloneqq \l<x,y\r>$ is a positive definite kernel (since $\sum_{i,j}c_i\bar c_j\l<x_i,x_j\r>=\l\|\sum_ic_ix_i\r\|^2$).
\end{exa}

So, if $X$ is a set and $\psi:X\to\mc H$ is map into a Hilbert space, then $K\l(x,y\r)\coloneqq \l<\psi\l(x\r),\psi\l(y\r)\r>$ is a positive definite kernel.
\begin{thm}[The GNS Construction]
        If $X$ is a set and $K$ a positive definite kernel, then there is a Hilbert space $\mc H$ and a map $\psi:X\to \mc H$ such that 
$$K\l(x,y\r)=\l<\psi(x),\psi(y)\r> \quad \forall x,y\in X.$$
Moreover, we can take $\mc H$ such that $\overline{\mr{Sp}\l(\psi\l(X\r)\r)}=\mc H$, and then the pair $\l(\mc H,\psi\r)$ is unique (in the sense that if $({\mc H}^\pr,\psi^\pr)$ is another such pair then there is an isometry $T:\mc H\to \mc H^\pr$ such that $\psi^\pr=T\circ \psi$).
\end{thm}
\begin{proof}
\def \tmch {\t{\mc H}}
        Define $\t{\mc H}$ to be the free vector space over $X$: that is, its elements are formal (finite) sums of the form $\sum_{i=1}^nc_ix_i$ for $x_i\in X$ and $c_i\in \mb C$. In particular, we have $\tmch=\mr{Sp}\l(X\r)$.\footnote{More formally, the elements of $\tmch$ are functions $f:X\to \mb C$ with finite support (i.e., for which $\{x\in X | f(x)\neq 0\}$ is finite), with addition and scalar multiplication defined pointwise. For each $x\in X$, one defines the Dirac function $\delta_x$ by setting $\delta_x(x)=1$ and $\delta_x(y)=0$ for $y\neq x$. Then $\tilde{\mc H}=\mr{Sp}(\{\delta_x\}_{x\in X})$. What we did is simply to identify each $x\in X$ with $\delta_x$.} Define an inner product in $\tmch$ as follows:
$$\l<\sum_{i=1}^n c_ix_i,\sum_{j=1}^na_jx_j\r>=\sum_{i,j}c_i\overline {a_j}K\l(x_i,x_j\r).$$
If $K$ is strictly positive, then $\tmch$ is an inner product space and we can take $\mc H$ to be its completion. In the more general case, where $K$ is just some positive definite kernel, one need to first quotient $\tmch$ by the subspace of vectors with norm zero.
\end{proof}
\begin{prop}
If $K_1,K_2:X\x X\to \mb C$ are positive definite kernels, then
\begin{enumerate}[1)]
\i \label{plus-ker}$K_1+K_2$ is positive definite.
\i \label{multi-ker}$K_1\cd K_2$ (with multiplication pointwise) is positive definite.
\end{enumerate}
\end{prop}
For the proof, we will use tensor products of Hilbert spaces. Recall that, if $V$ and $W$ and vector spaces, then their \emph{tensor
product }is a vector space $V\otimes W$ that admits a bilinear map
\begin{align*}
V\times W & \to V\otimes W\\
(v,w) & \mapsto v\otimes w
\end{align*}
which satisfies the following property: for every bilinear map $h:V\times W\to Z$
(into some vector space $Z$), there is a unique linear map $\tilde{h}:V\otimes W\to Z$
that satisfies
\[
\tilde{h}(v\otimes w)=h(v,w).
\]
In fact, any vector space $U$ that admits a bilinear map $V\times W\to U$ satisfying this property is canonically isomorphic to $V\otimes W$. The tensor product of two vector spaces can be defined more constructively, either
by taking the free vector space over the Cartesian product $V\times W$
and then modding out some relations (namely, $\lambda(v,w)=(\lambda v,w)=(v,\lambda w)$,
$(v_{1}+v_{2},w_{1}+w_{2})=\sum_{1\leqslant i,j\leqslant2}(v_{i},w_{j})$),
or by taking bases $B_{V}$ and $B_{W}$ of $V$ and $W$ respectively
and letting $V\otimes W$ be the free vector space over $\left\{ v\otimes w\middle|v\in B_{V},w\in B_{W}\right\} $.
When $V$ and $W$ are Hilbert spaces, we define the inner product
by
$$\l<\sum_i c_i\l(v_i\otimes w_i\r),\sum_{j}a_j\l(v^\pr_j\otimes w^\pr_j\r)\r>_{V\otimes W}=\sum_{i,j}c_i\overline{a_j}\l<v_i,v_j^\pr\r>_{V}\cd\l<w_i,w_j^\pr\r>_{W},$$
for every $v_{i}\in V,w_{j}\in W,a_{i},c_{j}\in\mathbb{C}$, $i\in\left\{ 1,\dots,n\right\} $,
$j\in\left\{ 1,\dots,m\right\} $. Then the completion of this space is the tensor product of these two spaces.

\begin{proof}
We leave (\ref{plus-ker}) as an exercise, and prove (\ref{multi-ker}). By the GNS construction, there are $\psi_i:X\to \mc H_i$ such that $K_i\l(x,y\r)=\l<\psi_i\l(x\r),\psi_i\l(y\r)\r>_{\mc H_i}$ for $i=1,2$.
Consider the tensor product $\mc H_1\otimes \mc H_2$. Consider the map $\psi_1\otimes \psi_2:x\mapsto \psi_1\l(x\r)\otimes \psi_2\l(x\r)$. We have
$$K_1\cd K_2\l(x,y\r)=K_1\l(x,y\r)\cd K_2\l(x,y\r)=\l<\psi_1\otimes \psi_2\l(x\r),\psi_1\otimes\psi_2\l(y\r)\r>_{\mc H},$$
so $K_1\cd K_2$ is also a positive definite kernel.
\end{proof}

\begin{cor}
\label{pos:exp}
If $K$ is a positive definite kernel, then $e^K\coloneqq\sum_{i\in \mb N}\fr{K^n}{n!}$ is also positive definite.
\end{cor}
\begin{cor}
If $\mc H$ is a Hilbert space and $s>0$, then $e^{-s\l\|\cd\r\|^2}$ is a positive definite function (that is, $K\l(x,y\r)\coloneqq e^{-s\l\|x-y\r\|^2}$ is a positive definite kernel).
\end{cor}
\begin{proof}
We have $$e^{-s\l\|x-y\r\|^2}=e^{-s\l\|x\r\|^2-s\l\|y\r\|^2+2s\mr{Re}\l<x,y\r>}=\l(e^{-s\l\|x\r\|^2}e^{-s\l\|y\r\|^2}\r)e^{2s\mr{Re}\l<x,y\r>}.$$
The kernel $\l(x,y\r)\mapsto e^{2s\mr{Re}\l<x,y\r>}$ is positive definite by Corollary \ref{pos:exp},
and the kernel $\l(x,y\r)\mapsto e^{-s\l\|x\r\|^2}e^{-s\l\|y\r\|^2}$ is positive definite by
\begin{exe} For any function $f:X\to \mb R$, the kernel $K\l(x,y\r)\coloneqq f\l(x\r)f\l(y\r)$ is positive definite. (Observe we don't require $f$ to be positive definite, only to be real).
\end{exe}
\end{proof}

\def \mchs {\t{\mc H}_s}
Thus, for any Hilbert space $\mc H$ the kernel ${K\l(x,y\r)\coloneqq}\;e^{-s\l\|x-y\r\|^2}$ is positive definite, and there is an embedding $\psi:\mc H\to \mchs$ such that $$\l<\psi\l(x\r),\psi\l(y\r)\r>=K\l(x,y\r)=e^{-s\l\|x-y\r\|^2}.$$
Therefore, no two vectors in $\psi\l(\mc H\r)$ are orthogonal, but the further away $x,y$ are from one another, the `closer' $\psi\l(x\r),\psi\l(y\r)$ are to be orthogonal. Moreover, $\l\|\psi\l(x\r)\r\|=1$ for every $x\in \mc H$, so $\psi\l(\mc H\r)$ is a subset of the unit sphere of $\mchs$.

We are now finally ready to complete the proof of Theorem \ref{FH:T}.

\begin{proof}[Proof of (T) $\implies$ (FH)]
Let $G$ be a group with property (T) which acts (continuously) by isometries on $\mc H$. We prove it admits a global fixed point, thus proving it has property FH.

Let $\l(K,\ep\r)$ be a Kazhdan pair such that $e\in K$, $\ep<\fr{\sqrt 2}{2}$, and let $x\in\mc H$ to be some vector. Set $d \coloneqq\mr{diam}\l(Kx\r)$.

Let $(\mchs ,\psi)$ be the pair obtained by the GNS construction with $X=\mc H$ and $K:\l(x,y\r)\mapsto e^{-s\l\|x-y\r\|^2}$, such that $\mchs=\overline{\mr{Sp}\l(\psi\l(\mc H\r)\r)}$. We can take any $s>0$, so take it to be sufficiently  small so that $e^{-sd^2}>1-\ep^4/2$, which will make sense shortly.

$G$ acts on $\mc H$, so we can define an action on $\psi\l(\mc H\r)$ by 
$$g.\psi\l(x\r)=\psi\l(g^{-1}.x\r).$$
$G$ acts by isometries on $\psi\l(\mc H\r)$, because
$$\l<g.\psi\l(x\r),g.\psi\l(y\r)\r>=e^{-\|g^{-1}.x-g^{-1}.y\|^2}=e^{-\l\|x-y\r\|^2}= \l<\psi\l(x\r),\psi\l(y\r)\r>.$$
There is a unique extension of this action (of $G$ on $\psi\l(\mc H\r)$) to an orthogonal representation of $G$ on $\mchs$. We leave it to you to check this.

Set $\hat v\coloneqq \psi\l(x\r)$. Then $\hat v$ is $\l(K,\ep^2\r)$ invariant, since, for all $k\in K$, we have $\l\|kx-x\r\|\le d$ (by definition of $d$) and therefore $\l<k\hat v,\hat v\r>_{\mchs}\ge e^{-sd^2}\ge1-\ep^4/2$, and hence $\l\|k\hat v-\hat v\r\|\le \ep^2$. Thus, by Lemma \ref{clo:inv}, there is an invariant vector at distance at most $\ep$ from $\hat v$. Therefore the entire orbit of $\hat v$, $G\hat v$, is at distance at most $\ep$ from this invariant vector as well (since $G$ acts by isometries), and hence it is at distance at most $2\ep$ from $\hat v$ by the triangle inequality. That is, for all $g\in G$
$$\l\|g\cd \hat v-\hat v\r\|\le 2\ep<\sqrt 2,$$
so (since $\hat v$, $g\cd \hat v$ are  both of norm $1$) the inner product is bounded from below,
$$\l<g\cd\hat v,\hat v\r>_{\mchs}\ge c>0$$ for some $c$ which is determined by $\ep$. But
$$e^{-s\l\|g\cd x-x\r\|^2}=K\l(g\cd x,x\r)=\l<\psi\l(g\cd x\r),\psi\l(x\r)\r>_{\mchs}=\l<g\cd \hat v,\hat v\r>_{\mchs}\ge c,$$ so $e^{-s\l\|g\cd x-x\r\|^2}$ is bounded from below by $c>0$ and therefore $\l\|g\cd x-x\r\|$ is bounded from above by some $M$ which is determined by $c$.

Thus, $Gx$ is a bounded orbit, and all orbits in $\mc H$ are bounded. Therefore, by Lemma \ref{fix:orb}, $G$ admits a global fixed point.
\end{proof}

\section{The Hyperbolic Space}
The Hyperbolic space $\mb H^n$ is the simply connected Riemannian manifold of dimension $n$ and constant curvature $-1$. 

In this section we will show that property (T) implies the fixed point property for hyperbolic spaces (i.e., if a Kazhdan group acts on a hyperbolic space, then
it admits a fixed point). We will assume some familiarity with such spaces and use basic results without proving them.

First, we demonstrate two models of the hyperbolic space $\mb H^n$.
\paragraph{The Hyperboloid Model}
Consider $\mb R^{n+1}$. We define the quadratic form $Q$ by
$$Q\l(x_0,\dots,x_n\r)=x_0^2-\sum_{i=1}^nx_i^2.$$
Consider
$$\l\{\bar x\in \mb R^{n+1}\mi Q\l(\bar x\r)=1\r\}=Q^{-1}\l(\l\{1\r\}\r).$$
This is a hyperboloid, with two connected components. We let the upper one be our model of the hyperbolic space,
$$\mb H^n =\l\{\bar x\in \mb R^{n+1}\mi Q\l(\bar x\r)=1,x_0>0\r\}.$$
This is a submanifold of $\mb R^{n+1}$ because $Q:\mb R^{n+1}\to \mb R$ is a submersion. In fact, $\mb H^n$ is a Riemannian manifold, and $Q$ is positive definite when restricted to it.

\paragraph{The Upper Half Space Model}
Consider $\mb R^{n,+}\coloneqq \l\{\l(x_1,\dots,x_n\r)\mi x_n>0\r\}$. Set

$$ds^2=\fr{\sum_{i=1}^{n}dx_i^2}{x_n^2}.$$
That is, the length of a path $\g:\l[0,L\r]\to \mb H^n$ is defined by
$$\ell\l(\g\r)=\int_0^L\!  \fr{\sqrt{\sum_{i=1}^{n}\dot x_i^2(t)}}{x_n(t)}  \;dt,$$
where  $\g\l(t\r)=\l(x_1\l(t\r),\dots,x_n\l(t\r)\r)$. We define the distance between two points in this space to be the infimum of lengths of paths between them. It turns out a minimum is actually achieved in this case, i.e.\ for every two points there is a path of minimal length from one to the other, and this path is unique (up to reparametrisation). Such paths are called geodesics. Thus, this turns $\mb H^n$ into a uniquely geodesic metric space, which is moreover uniformly convex. 

It's not difficult to prove that geodesics in this space are generalised half circles orthogonal to $\partial \mb H^n=\l\{\l(x_0,\dots,x_n\r)\mi x_n=0\r\}$ (where we consider lines as ``generalised circles").

In the case $n=2$, you can think of $\mb H^2$ as $\l\{z\in \mb C\mi \mr{Im} z>0\r\}$. The connected component of the identity of the group of isometries of $\mb H^2$, $\mr{Isom}^\circ\!\l(\mb H^2\r)$, is isomorphic to $\mr{PSL}_2(\mb R)=\mr{SL}_2(\mb R)/\!\l\{\pm I\r\}$: every matrix $\bigl(\begin{smallmatrix}a &b\\c &d\end{smallmatrix}\bigr)$
acts on $z$ by the M\"obius transformation
$$z\mapsto \fr{az+b}{cz+d}\quad.$$
The area of a triangle in $\mb H^2$ is always less than $\pi$ (it's actually $\pi$ minus the sum of its angles). If you consider the case of a triangle with vertices at infinity (either points in $\mb R$ or the point $\infty\in \hat {\mb C}$) then the area is exactly $\pi$ (all lines will be tangent to one another).

The fact that triangles have bounded area implies that there exists a $\delta>0$ such that every triangle is $\delta$-thin, i.e.\ every edge is contained in a $\delta$-neighbourhood of the union of the other two edges.

In the case $n=3$, we have $\mr{Isom}^\circ\!\l(\mb H^3\r)\cong \mr{PSL}_2(\mb C)$. 

\begin{thm}
\label{T:hyp}
If $G$ has property (T), it has the fixed point property for Hyperbolic spaces (i.e.\ it admits a fixed point whenever it acts (continuously) by isometries on $\mb H^n$).
\end{thm}

First, we need to mention a few facts about the group of isometries on $\mb H^n$. Denote $G_n\coloneqq \mr{Isom}^\circ (\mb H^n)$ (the connected component of the identity); then $G_n\cong \mr{SO}\l(n,1\r)$. Reflections (over hyperspaces in $\mb H^n$ with codimension $1$) aren't in $G_n$, they are in $\mr{Isom}(\mb H^n)$. If $s$ is a reflection, then $\left\langle G_n,s\right\rangle=\mr{Isom}(\mb H^n)$.

The first fact we need is that $G_n$ acts transitively on subspaces of codimension $1$, called \textit{hyperplanes}. That is, any hyperplane can be mapped onto any hyperplane by some element in $G_n$. For $n=2$, hyperplanes are the same as geodesics, so this can be more easily visualised: $G_2$ can send any geodesic (i.e.\ generalised half-circle) onto any geodesic (i.e.\ generalised half-circle). 

Every hyperplane cuts $\mb H^n$ into two half-spaces, so can define an \textit{oriented hyperplane}, whose orientation may be either one of these two half spaces. Set $$X\ceq \l\{\text{oriented hyperplanes}\r\}=\l\{\text{half-spaces}\r\}.$$ $G_n$ acts transitively on $X$. Given a hyperplane $K\subseteq \mb H^n$, every isometry of $K$ extends to an orientation-preserving isometry on $\mb H^n$. Thus, given a half-space $K^+$, the subgroup $\mr{Stab}_{G_n}K^+=\l\{g\in G_n\mi g\l(K^+\r)=K^+\r\}$ is isomorphic to $G_{n-1}$.

As we stated before, $G_n\cong \mr{SO}\l(n,1\r)$, so it's unimodular. Since $X$ is a homogeneous space, we have
$$X=G_n/G_{n-1}.$$
Since $G_n$, $G_{n-1}$ are unimodular, $G_n/G_{n-1}$ admits a $G_n$-invariant Borel regular measure $\mu$ such that $\mu\l(C\r)<\infty$ for compact $C$ and $\mu\l(U\r)>0$ for nonempty open $U$. 

\begin{proof}[Proof of Theorem \ref{T:hyp}]
        Assume a group $\G$ acts on the hyperbolic space $\mb H^n$ without a fixed point. Observe that this means that the orbits of this action are unbounded, since the barycentre of a bounded orbit is a fixed point. We will use this fact to prove $\G$ doesn't have property (T), by constructing an action on some Hilbert space with unbounded orbits.
        
        Our Hilbert space is $\mc H\ceq L_2\l(X\r)$ where $X$ is the space of half-spaces in $\mb H^n$. $\G$ acts on $\mb H^n$, so it acts on $X$ (since $G_n$ does). We have a natural representation $\rho:\G \to U\l(\mc H\r)$ defined by
        $$\l(\rho\l(\g\r) f\r)\l(x\r)=f\l(\g^{-1}.x\r).$$
         We need to find an unbounded $\rho$-cocycle.
        
        For every $p,q\in\mb H^n$, we define $\Om_{p,q}\ss X$ to be the set of half-spaces containing $q$ but not $p$. The subset $\overline{\Omega_{p,q}}$ is compact (in dimension two, this is easy to visualise -- the choice of a point in the interval between $p$ and $q$ is compact, and the choice an angle between $0$ and $\pi$ or $\pi$ and $2\pi$ is compact -- and, more generally, this follows from compactness of closed subset of the $n$-sphere), so  $\mu\l(\Omega_{p,q}\r)<\infty$.
        
        Fix some $x_0\in\mb H^n$. Define $c:\G\to\mc H$ by
        $$c\l(\g\r)=1_{\Om_{x_0,\g\cd x_0}}-1_{\Om_{\g\cd x_0,x_0}}.$$ Checking the cocycle equation is straightforward  and is similar to the proof of Theorem \ref{FH=>FA} (it becomes clear once one draws a picture of all the dots).
        
        Since $\mu\l(\Om_{p,q}\r)$ is proportional to $d\l(p,q\r)$, and by our assumption the $\G$-orbit of $x_0$ is unbounded, we get that $\mu\l(\Om_{x_0,\g\cd x_0}\r)$ is unbounded and therefore $\l\|1_{\Om_{x_0,\g\cd x_0}}-1_{\Om_{\g\cd x_0,x_0}}\r\|_2$ is unbounded. Thus, the cocycle $c$ is unbounded, and $\G$ doesn't have property (T).
\end{proof}
\section*{Further Reading}
For a more thorough analysis of property FH and similar properties, the readers are referred to \cite{bekka2008kazhdan}.
To read more about group actions on trees (and in particular the proof of Theorem \ref{thm:FA-eq}), we warmly recommend Serre's book \cite{Trees}.

\chapter{Dynamics of Linear Groups}
In this chapter and the next we demonstrate some results in analytic group theory. 
We will prove the following theorems:
\begin{thm*}[The Jordan theorem]
For every $n\in\mb N$ there is some $m\in\mb N$ such that every finite subgroup of $\glnc$ has an abelian subgroup of index at most $m$.
\end{thm*}

\begin{thm*}[The Bieberbach theorem]
For every $n\in\mb N$ there is some $m\in\mb N$ such that every  flat, compact $n$-dimensional Riemannian manifold $M$ admits a finite covering map $\t M\to M$ of degree at most $m$, where $\t M$ is a torus.

\end{thm*}

\begin{thm*}[Zimmer's theorem]
An infinite group with property (T) cannot act conformally on the sphere $\mathbb{S}^2$.
\end{thm*}
\begin{thm*}[The Tits alternative]
    A finitely generated subgroup of $\mathrm{GL}_{n}(\mathbb{C})$ is either  virtually solvable or contains a non-abelian free subgroup.
\end{thm*}

We will deduce the first two theorems from the famous Margulis lemma, which has many other applications. In order to prove the third theorem we will make use of the tree associated to $\mr{SL}_2\l(\mb Q_p\r)$, which is the simplest example from the important Bruhat--Tits theory.

\section{The Margulis Lemma and Applications}
\subsection{Zassenhaus' Theorem}
The first result we will use is Zassenhaus' theorem. 

Consider $G\ceq\glnc.$
Take $x,y\in G$ which are close to $1$, i.e., such that $x=1+X$, $y=1+Y$ where $\l\|X\r\|<\ep$, $\l\|Y\r\|<\delta$ for some small $\ep,\delta>0$. Consider 
$$\l[x,y\r]-1=xyx^{-1}y^{-1}-1=\l(xy-yx\r)x^{-1}y^{-1}.$$
By direct computation,  $\l[x,y\r]-1=\l(XY-YX\r)x^{-1}y^{-1}$, so
\begin{align*}
\l\|\l[x,y\r]-1\r\|&=\l\|\l(XY-YX\r)x^{-1}y^{-1}\r\|\le\l\|XY-YX\r\|\cd\l\|x^{-1}\r\|\cd\l\|y^{-1}\r\|
\\&\le \l(\l\|XY\r\|+\l\|YX\r\|\r)\cd\l\|x^{-1}\r\|\cd\l\|y^{-1}\r\|.
\end{align*}
If $\ep$ and $\delta$ are small enough, we have $\|x^{-1}\|,\|y^{-1}\|<2$, so
$$\l\|\l[x,y\r]-1\r\|\le 8\ep\delta.$$

Take $\Om=\l\{x\in\glnc\mi \l\|x-1\r\|<\fr{1}{8r}\r\}$ for some $r>1$. Then 
$$\Om^{[n]}\ceq \l\{xyx^{-1}y^{-1}\mi x\in \Om, y\in \Om^{n-1}\r\}$$
approaches $\l\{1\r\}$, because, if $v\in\Om^{[n]}$, then $\l\|v-1\r\|<\fr{1}{r^n}$. Indeed, if $\l\|x-1\r\|<\fr{1}{8r}$, $\l\|y-1\r\|<\fr{1}{r^{n-1}}$, then by the above $\l\|\l[x,y\r]-1\r\|<8\fr{1}{8r}\cd\fr{1}{r^{n-1}}=\fr{1}{r^n}$.

\begin{exe}
If $\G=\l<S\r>$ and $S^{[n]}=\l\{1\r\}$, then $\G$ is nilpotent of degree at most $n$.
\end{exe}

\begin{hint}
    Show that $S^{[n-1]}$ is contained in the centre, divide by the centre, and use induction.
\end{hint}

\begin{thm}[Zassenhaus]
If $\G\le \glnc$ is discrete and generated by its intersection with $\Om$, then $\G$ is nilpotent.
\end{thm}
\begin{proof}
Set $S\ceq \G\cap \Om$. Then $S^{[k]}\overset{k\to\infty}{\too}\l\{1\r\}$ by the above, so by discreteness  there is some $n\in\mb N$ such that $S^{[n]}=\l\{1\r\}$, and hence (by the last exercise) $\G$ is nilpotent.
\end{proof}

Using Lie algebras and the exponential map, a similar reasoning provides a stronger version of this theorem, which we will need. Namely, that $\G$ is not only nilpotent, but is also contained in some connected nilpotent Lie subgroup. (In that case our restriction on $r$ will be more than just $r>1$.) We refer to \cite[Theorem 8.6]{Rag} for more details.

\subsection{The Margulis Lemma}
We now prove the Margulis lemma.

\begin{thm}[The Margulis lemma]
Let $G$ be a connected Lie group acting transitively, continuously, properly and by isometries on a connected Riemannian manifold $X$. There exist $\ep>0$, $m\in\mb N$ such that, if $\G\le G$ is discrete and
$$\G=\l<\g\in\G\mi d\l(\g.x,x\r)<\ep\r>$$
for some $x\in X$, then $\G$ has a subgroup of index at most $m$ which is contained in a connected nilpotent Lie group.
\end{thm}

This result is very similar to Zassenhaus' theorem. It seems Zassenhaus' theorem is stronger: it proved $\G$ itself is contained in a connected nilpotent Lie subgroup, and not just some subgroup of it of index at most $m$ (the fact that it's about $\glnc$ doesn't matter much, because it can be easily generalised). But in the Margulis lemma we didn't require $\G$ to be generated by a neighbourhood around the identity: we only required it  to be generated by a neighbourhood of a stabiliser, which might be big (albeit compact).

\begin{proof}
Fix some $x\in X$. Let $K\ceq G_x$ be the stabiliser of $x$, which is compact (by the propriety of the action).

Let $\Om$ be a Zassenhaus neighbourhood\ (i.e.\ such that Zassenhaus' theorem in its stronger version holds). Consider the set $G_{x,1}=\l\{g\in G\mi d\l(g.x,x\r)\le 1\r\}$. This is a compact set.

Let $U\ss G$ be a relatively compact (i.e.\ with compact closure) neighbourhood of the identity such that  $U^{-1}=U$ and $U^2\ss \Om$ (we can find such because $G$ is a topological group).

Set $$m=\l\lceil\fr{\mu\l(G_{x,1}.U\r)}{\mu\l(U\r)}\r\rceil$$
for $\mu$ a Haar measure on $G$, and set $\ep=1/m$. We will later consider the reason why $m,\ep$ do not depend on $x$.

Assume $\G\le G$ is a discrete subgroup such that $\G=\l<S\r>$ for
$$S=\l\{\g\in \G\mi d\l(\g.x,x\r)<\ep\r\}.$$
We need to prove $\G$ has a subgroup of index at most $m$  which is contained in a connected nilpotent Lie subgroup.

Set $A=\l<\G\cap \Om\r>$. By Zassenhaus' theorem, $A$ is contained in a connected nilpotent subgroup. It remains to show $\l[\G:A\r]\le m$.

Let  $Y\ceq G\l(A\!\setminus\!\G,S\r)$ be the Schreier graph (whose vertices are the cosets in $A\!\setminus\!\G$ and whose edges are cosets with representatives $\g_1,\g_2$ such that $\g_1 s=\g_2$ for some $s\in S$), and assume by contradiction $\l|A\!\setminus\!\G\r|\ge m+1$.

Consider the (closed) balls around the identity in $Y$ with radius $k$, $k=0,1,\dots, m$. Whenever $k$ grows by one, the ball must contain at least one more element, unless it already contained all of the cosets. In either case, since the ball with radius $0$ contains one element, the ball with radius $m$ must contain at least $m+1$ elements. Let $\g_1,\dots,\g_{m+1}\in \G$ be representatives of distinct vertices in this ball. We can choose them so that  $\g_i\in S^m$ for all $i$, and necessarily $\g_j^{-1}\g_i\notin A$ for all $i\ne j$.

If $\g_i\in S^m$, then $d\l(\g_i.x,x\r)\le m\cd \ep=1$, which implies $\g_i\in G_{x,1}$.

If $\g_j^{-1}\g_i\notin A$, then $\g_j^{-1}\g_i\notin \Om\supseteq U^2$, which implies $\g_j U\cap \g_i U=\varnothing$ (for all $i\ne j$).

So, we have $m+1$ translations of $U$ (the sets $\g_1 U,\dots \g_{m+1} U$) which are disjoint and contained in $G_{x,1}.U$, contradicting the definition of $m$.

Lastly, the constants $m,\ep$ are independent of $x$ since the action is transitive, and had we taken a different point, we could have just conjugated all the sets and got the same numbers. 
\end{proof}

\subsection{The Jordan Theorem}
The Margulis lemma has many applications. In this subsection and the next we shall present two elementary but classical examples.

\begin{thm}[The Jordan theorem]
Given $n\in\mb N$, there is an $m\in\mb N$ such that every finite subgroup of $\glnc$ has an abelian subgroup of index at most $m$.
\end{thm}
\begin{rem*}Observe that obviously for every finite subgroup there is such an $m$ (namely the number of its elements), but the theorem claims there is one such $m$ for \emph{all}  finite subgroups, regardless of their size.

One should also note we can find finite abelian subgroups of $\glnc$ as large as we wish (for example, rotations by $\fr{2\pi}{n}$ for all $n$). What the theorem basically says is that, ``up to abelian stuff", the size of finite subgroups of $\mr{GL}_n(\mb C)$ is uniformly bounded. 
\end{rem*}
\begin{proof}
Considering $\mb C^n$ as a real $2n$-dimensional vector space, one sees that $\glnc \le \text{GL}_{2n}(\mb R)$. This allows us to replace $\mb C$ by $\mb R$ in the statement of the theorem. Furthermore, a finite subgroup $\G$ of $\text{GL}_{n}(\mb R)$ preserves some real inner product on $\mb R^{n}$. To see this, take some real inner product $\l(\cd,\cd\r)$ on $\mb R^{n}$; then
$$\l<v_1,v_2\r>\ceq \fr{1}{\l|\G\r|}\sum _{\g\in\G}\l(\g.v_1,\g.v_2\r)$$
is an inner product preserved by $\G$\ (we're in fact using the amenability of $\G$ here). Up to replacing $\G$ by an index two subgroup, we may suppose that it is contained in $\mr {SO}_n(\mb R)$. Therefore we have shown that it is enough to prove the analogous claim for finite subgroups of $\mr {SO}_n(\mb R)$.

Being compact, $\mr {SO}_n(\mb R)$ acts properly on a point, which is a zero dimensional Riemannian manifold. Thus, 
the Jordan theorem follows from the Margulis lemma, applied to this trivial case. The fact that the subgroup of index at most $m$ we obtain is abelian, rather than just nilpotent, is a consequence of the following basic fact from Lie theory: a connected compact solvable Lie group is necessarily abelian. We will prove it in the next section.
\end{proof}
\subsection{A Little Lie Theory}
In this subsection, we recall some basic notions from Lie theory in order to prove:

\begin{lem}\label{lem:compact-sol-abelian}
A connected compact solvable Lie group is abelian.
\end{lem}

If $G$ is a Lie group and $1_G$ is its identity element, there is a natural way to make the tangent space $\mf{g}\coloneqq \mr{T}_{1_G}G$ into a \textit{Lie algebra}. That is, there is a natural way to define a bilinear map
\begin{align*}
\mathfrak{g}\times\mathfrak{g} & \to\mathfrak{g}\\
(X,Y) & \mapsto[X,Y],
\end{align*}
called the \textit{Lie bracket}, that satisfies the axioms of a Lie algebra. Namely, we have $[X,X]=0$ and $[X,[Y,Z]]+[Y,[Z,X]]+[Z,[X,Y]]=0$ for all $X,Y,Z\in\mf g$. A subspace $\mf a\subseteq \mf g$ of a Lie algebra $\mf g$ is a \textit{subalgebra} if $[\mf a,\mf a]\subseteq \mf a$, and it is an \textit{ideal} if $[\mf g,\mf a]\subseteq \mf a$. 

If $H\le G$ is a Lie subgroup,  we can identify the Lie algebra $\mf h$ of $H$ with a subalgebra of the Lie algebra $\mf g$ of $G$. It is an ideal if and only if $H$ is normal.

Let $G$ be a Lie group. Its action on itself by conjugations descends to an action on its own Lie algebra, called the \textit{adjoint action}. In other words, we have a homomorphism
\[
\mr{Ad}:G\to \mr{GL}(\mf g).
\]
One can then show that a subspace $\mf a\subseteq \mf g$ is an ideal of $\mf g$ if and only if it is $\mr{Ad}(G)$-invariant.

A Lie algebra $\mf g$ is called \textit{abelian} if $[X,Y]=0$ for every $X,Y\in\mf g$. An abelian Lie algebra is essentially just a vector space. A Lie algebra $\mf g$ is called \textit{simple} if it is not abelian and $\{0\},\mf g$ are its only ideals. In this case, we have $[\mf g,\mf g]=\mf g$, since $[\mf g,\mf g]$ is an ideal (and it cannot be $\{0\}$, or else $\mf g$ would be abelian). 

We set $\mf g^{(0)}=\mf g$ and $\mf g^{(n+1)}=[\mf g^{(n)},\mf g^{(n)}]$. A Lie algebra $\mf g$ is \textit{solvable} if $\mf g^{(n)}=\{0\}$ for some $n\in\mb N$. We set $\mf g_0=\mf g$, $\mf g_{n+1}=[\mf g,\mf g_n]$, and say $\mf g$ is \textit{nilpotent} if $\mf g_n=\{0\}$ for some $n\in\mb N$. It is not hard to show that a connected Lie group is abelian (respectively nilpotent, respectively solvable) if and only if its Lie algebra is abelian (respectively nilpotent, respectively solvable).

We are now ready to prove Lemma \ref{lem:compact-sol-abelian}.
\begin{proof}[Proof of Lemma \ref{lem:compact-sol-abelian}]
Let $G$ be a connected compact solvable Lie group, $\mathfrak{g}$
its Lie algebra. Consider the adjoint action of $G$ on $\mathfrak{g}$.
As usual, since $G$ is compact, it is easy to see that $\mathfrak{g}$
admits an inner product that is $\mathrm{Ad}(G)$-invariant: take
any inner product $(\cdot,\cdot)$ on $\mathfrak{g}$, and then define
a new inner product, $[\cdot,\cdot]$, by
\[
[X,Y]=\int_{G}[\mathrm{Ad}(g)(X),\mathrm{Ad}(g)(Y)]dg
\]
(with respect to the Haar probability measure). If $\mathfrak{a}\subseteq\mathfrak{g}$
is $\mathrm{Ad}(G)$-invariant, then its orthogonal complement $\mathfrak{a}^{\perp}$
is $\mathrm{Ad}(G)$-invariant (since $\mathrm{Ad}(G)$ preserves
the inner product $[\cdot,\cdot]$). In other words, if $\mathfrak{a}\subseteq\mathfrak{g}$
is an ideal, its orthogonal complement is an ideal. It follows by
induction that $\mathfrak{g}$ is a direct sum of ideals which are
`minimal', in the sense that they contain no nontrivial ideals themselves.
This means that they are either one-dimensional and abelian, or simple.
That is, there are ideals $\mathfrak{a}_{1},\dots,\mathfrak{a}_{n}\trianglelefteqslant\mathfrak{g}$
in $\mathfrak{g}$, each of which is either simple as a Lie algebra
or $1$-dimensional and abelian, such that
\[
\mathfrak{g}=\mathfrak{a}_{1}\oplus\cdots\oplus\mathfrak{a}_{n},
\]
a direct sum of Lie algebras. It follows that 
\[
[\mathfrak{g},\mathfrak{g}]=[\mathfrak{a}_{1},\mathfrak{a}_{1}]\oplus\cdots\oplus[\mathfrak{a}_{n},\mathfrak{a}_{n}].
\]
Since each $\mathfrak{a}_{i}$ is either simple or abelian, we either
have $[\mathfrak{a}_{i},\mathfrak{a}_{i}]=\mathfrak{a}_{i}$ or $[\mathfrak{a}_{i},\mathfrak{a}_{i}]=\left\{ 0\right\} $.
Since $G$ is solvable, $\mathfrak{g}$ is solvable as well. Therefore,
it is not possible that $[\mathfrak{a}_{i},\mathfrak{a}_{i}]=\mathfrak{a}_{i}$
for any $i$, or we would get that $\mathfrak{a}_{i}\subseteq\mathfrak{g}^{(n)}$
for all $n\in\mathbb{N}$. Therefore, $[\mathfrak{g},\mathfrak{g}]=\left\{ 0\right\} $.
This means that $\mathfrak{g}$ is abelian, so $G$ is abelian too. 
\end{proof}

\subsection{The Bieberbach Theorem}
A \textit{crystallographic manifold} is a compact flat Riemannian manifold. Equivalently, it is a compact Riemannian manifold that is locally isometric to $\mb R^n$. It is a standard fact of Riemannian covering theory that an $n$-dimensional crystallographic manifold is of the form $\G\backslash \mb R^n$, where $\G$ is a torsion-free discrete cocompact subgroup of $\mr{Isom}(\mb R^n)$. 

One obvious example of a crystallographic manifold is the torus $\mb R^n/\mb Z^n$. Hilbert's 18th problem was whether this is virtually the only example: more precisely, whether every crystallographic manifold is finitely covered by a torus. For example, the Klein Bottle is  a crystallographic manifold that is not a torus, but admits a $2$-covering by a torus. The Bieberbach theorem, which we will now prove, answers Hilbert's question in the affirmative. 

Recall that, by Lemma \ref{iso-is-affine} applied to $\mb R^n$, we have $\mr{Isom}(\mb R^n)=\mr{O}(n)\ltimes \mb R^n$. 
\begin{thm}\label{Bieberbach}
For every $n\in \mb N$ there is $m\in\mb N$ such that, if $\G$ is a discrete, torsion-free, cocompact subgroup of $\mr{Isom}(\mb R^n)=\mr{O}(n)\ltimes \mb R^n$, 
then $\G$ has a subgroup of index at most $m$ consisting of translations. 
\end{thm}
\begin{cor}[The Bieberbach theorem]
        For every $n\in\mb N$ there is $m\in \mb N$ such that every crystallographic $n$-manifold admits a covering of degree at most $m$ by a torus. 
\end{cor}
We shall deduce the Bieberbach theorem from the Margulis lemma, and the observation that the $\ep$ of the Margulis lemma in the case $G=\son\ltimes\mb R^n$ can be taken to be infinity. This is true since one can use homotheties on $\mb R^n$. A \textit{homothety} is a transformation of a Euclidean space that is a multiplication by some constant positive scalar, relative to some base point (so, if our base point is $0$, it is a linear transformation of the form $v\mapsto \lambda v$ for some $\lambda>0$). Conjugating by a homothety can be visualised as looking at $\mb R^n$ from afar. By going further and further away, we can make any finite subset of $\mb R^n$ have an arbitrarily small diameter, which means that, since the $\ep$ in the Margulis lemma is positive, it can be arbitrarily large.

We will need a few lemmas. We leave the proof of the first as an exercise.
\begin{lem}\label{disc-of-Rn}
    Let $\G\leqslant\mb R^n$ be a discrete subgroup. Show that $\G$ is isomorphic to $\mb Z^k$ for some $k\leqslant n$, and that $\G$ is cocompact if and only if $k=n$.
\end{lem}
\begin{exe}
    Prove Lemma \ref{disc-of-Rn}.
\end{exe}
\begin{lem}\label{Min}
    Every isometry $\alpha : \mb R^n\to\mb R^n$ admits a subspace $V_t$ and a point $p\in V_t^\perp$ such that $\alpha$ acts on $V_t+p$ by translations, and
    \[
    V_t+p=\{x\in \mb{R}^n|d(\alpha(x),x)=\min_{y\in \mb R^n}d(\alpha(y),y)\}.
    \]
\end{lem}
\begin{proof}
    Let $A\in\mr{O}_n(\mb R^n)$ and $w\in \mb R^n$ be the linear part and the translation part of $\alpha$, so that $\alpha(x)=Ax+v$ for all $x\in\mb R^n$. Let $V_t$ be the (possibly trivial) eigenspace of the eigenvalue $1$ of $A$, and let $V_r$ be its orthogonal complement. Then $A$ fixes $V_t$ and preserves $V_r$, and the restriction of $A$ to $V_r$ does not fix any non-zero vector. 
    
    Write $w=w_t+w_r$ for $w_t\in V_t$ and $w_r\in V_r$. Then, for every $x=x_t+x_r\in \mb R^n$, we have
    \[
        \alpha(x)=Ax_t+Ax_r+w_t+w_r=(x_t+w_t)+Ax_r+w_r.
    \]
    Defining $\alpha_t:V_t\to V_t$ by $\alpha(x_t)=x_t+w_t$ for $x_t\in V_t$ and $\alpha_r:V_r\to V_r$ by $\alpha_r(x_r)=Ax_r+w_r$ for $x_r\in V_r$, we get that $\alpha(x_t+x_r)=\alpha_t(x_t)+\alpha_r(x_r)$ for all $x_t\in V_t$ and $x_r\in V_r$.

    Now, consider the restriction of $A$ to $V_r$. Since $1$ is not an eigenvalue of $A\restriction_{V_r}$, the transformation $A-\mr{Id}$ is invertible, and we can set 
    \[
    p=(A-\mr{Id})^{-1}(-w_r)\in V_r.
    \]
    Direct computation shows that $p$ is the unique fixed point of $\alpha_r:V_r\to V_r$. 

    Thus, the affine subspace $V_t+p$ is preserved by $\alpha$, and $\alpha$ acts on it by translations of $w_t$. In the `$V_r$-direction', $\alpha$ acts by rotations. Since the restriction of $A$ to $V_r$ does not admit $1$ as an eigenvalue, it acts on it by a faithful rotation, without `fixed directions'. It is then intuitively clear that $\alpha$ moves the points in $V_t+p$ the least. That is: 
    \[
    V_t+p=\{x\in \mb{R}^n|d(\alpha(x),x)=\min_{y\in \mb R^n}d(\alpha(y),y).
    \]
    To see this formally, let $x=v_{t}+v_{r}+p$. Then 
    \[
    \alpha(x)=\alpha_{t}(v_{t})+\alpha_{r}(v_{r}+p)=v_{t}+w_{t}+\alpha_{r}(v_{r}+p)
    \]
    so (since $V_{t}\perp V_{r}$)
    \[
    \left\Vert \alpha(x)-x\right\Vert =\left\Vert w_{t}+\alpha_{r}(v_{r}+p)\right\Vert =\left\Vert w_{t}\right\Vert +\left\Vert \alpha_{r}(v_{r}+p)\right\Vert .
    \]
    Since $p$ is the unique fixed point of $\alpha_{r}$, it is clear that $\Vert \alpha(x)-x\Vert$ is minimal when $x_r=0$, i.e., when $x\in V_t+p$.
    \end{proof}
\begin{lem}\label{Abelian-sgp}
    Let $A\leqslant \mr{O}_n(\mb R^n)\ltimes \mb R^n$ be an abelian, torsion-free, discrete, cocompact subgroup. Then $A$ consists of  translations (i.e., $A\leqslant \mb R^n$) and is isomorphic to $\mb Z^n$.
\end{lem}
\begin{proof}
    We prove this by induction on the dimension. Thus, assume the lemma is true for dimensions strictly smaller than $n$.

    For $\alpha\in A$, set $d_\alpha(x)=d(\alpha.x,x)$. By Lemma \ref{Min},  $\mr{Min}(\alpha)=\{x\in\mb R^n|d_\alpha(x)=\min d_\alpha\}$ is an affine subspace, and $\alpha$ acts by translations on it.
    Now, let $\beta\in A$; since $\beta\alpha=\alpha\beta$, we get
that 
\begin{align*}
d_{\alpha}(\beta x) & =d(\beta x,\alpha\beta x)=d(\beta x,\beta\alpha x)=d(x,\alpha x)=d_{\alpha}(x),
\end{align*}
so $\mathrm{Min}(\alpha)$ is $\beta$-invariant. Therefore, $\mathrm{Min}(\alpha)$
is $A$-invariant (for every $\alpha\in A$). We want to show $\mathrm{Min}(\alpha)=\mathbb{R}^{n}$.

Consider the restriction of the action of
$A$ to $\mathrm{Min}(\alpha)$. The map $\beta\mapsto \beta\restriction_{\mathrm{Min}(\alpha)}$
is injective. To see this, suppose $\beta\neq 1$ and $\beta\restriction_{\mathrm{Min}(\alpha)}=1$.
Then $\beta$ admits some fixed point $p$. The stabiliser of $p$ inside $A$ is both discrete and compact, hence finite; therefore, $\beta\in A$ is nontrivial and of finite order, contradicting the assumption $A$ is torsion-free. 

Assume by contradiction that $\mathrm{Min}(\alpha)\neq\mathbb{R}^{n}$.
Then $k\coloneqq\dim\mathrm{Min}(\alpha)$ is smaller than $n$, so (by the induction hypothesis), $A\cong\mathbb{Z}^k$. By Exercise \ref{disc-of-Rn}, $A$ cannot be cocompact, which is a contradiction.

Thus, $\mr{Min}(\alpha)=\mb R^n$. Therefore, by Lemma \ref{Min}, $A$ acts by translations.
\end{proof}
Lastly, we will need the following fact from Lie theory.
\begin{lem}\label{nil-is-ab}
A connected nilpotent Lie subgroup of $\son\ltimes \mb R^n$ is abelian.
\end{lem}
\begin{proof}
    Let $N$ be a connected, nilpotent subgroup of $\son\ltimes \mb R^n$, and denote by $\mf n$ its Lie algebra. 

    Denote by $\pi:\mf{so}_n(\mb R)\ltimes \mb R^n\to \mf{so}_n(\mb R)$ the projection. Then $\pi(\mf n)$ is a nilpotent subalgebra of $\mf{so}_n(\mb R)$, hence abelian (by Lemma \ref{lem:compact-sol-abelian}). We claim that, furthermore, the action of $\pi(\mf n)$ on $\mf n\cap \mb R^n$ is trivial, which  will imply $\mf n$ is abelian.

    Assume by contradiction that there is an element $X\in \mf n$ and an element $Y\in \mf n\cap \mb R^n$ such that $\pi(X)(Y)=\mr{ad}(X)(Y)\neq 0$. Since $\pi(X)$ is a normal matrix, it is diagonalisable over $\mb C$, so it follows that $\mr{ad}(X)^n(Y)\neq 0$ for all $n$. But $\mf n$ is a nilpotent Lie algebra, say of degree $n$, so $\mr{ad}(X)^n(Y)=[X,[X,\dots[X,Y]\cdots]=0$, a contradiction.
\end{proof}
\begin{proof}[Proof of Theorem \ref{Bieberbach}] 
Let $\G$ be a discrete, torsion-free, cocompact subgroup of $\mr O_n(\mb R^n)\ltimes \mb R^n$. By paying a small price in the index, and taking a subgroup of $\G$ of index at most $2$, we may assume for convenience sake that $\G$ is a subgroup of $\son\ltimes\mb R^n$. The group $\G$ is finitely-generated, since it is the fundamental group of the compact manifold $\G\backslash\mb R^n$.

Let $\g_1,\dots,\g_k$ be a finite generating subset of $\G$. 
Since one can take Margulis constant $\ep$ guaranteed by the Margulis lemma to be arbitrarily large, the group $\G=\langle \g_1,\dots,\g_k\rangle$ admits a subgroup of index at most $m$ (where $m$ depends only on $n$) that is contained in a connected nilpotent subgroup of $\son\ltimes\mb R^n$. By Lemma \ref{nil-is-ab}, this subgroup is abelian, so the result follows from Lemma \ref{Abelian-sgp}.
\end{proof}

One can actually prove a more general version of Theorem \ref{Bieberbach}, not assuming co-compactness.

\begin{exe}
Show that for every $n$ there is $m$ such that, if $\G$ is a discrete torsion-free subgroup of 
$\son\ltimes\mb R^n$, then $[\G:\G\cap \mb R^n]\le m$. Show furthermore that $\G\cap \mb R^n$ is isomorphic to $\mb Z^k$ for some $k\le n$, and that $k=n$ if and only if $\G$ is co-compact.
\end{exe}

\begin{rem}
    It is also possible to remove the assumption that $\G$ is torsion-free and obtain a result that generalizes simultaneously both Bieberbach's and Jordan's theorems. 

    Adding back the co-compactness assumption, the statement become as elegant as in the theorem.
    A creature of the form $\G\backslash M$, where $M$ is a homogeneous Riemannian manifold and $\G\leqslant \mr{Isom}(M)$ is a discrete subgroup that might have torsion, is called an \textit{orbifold}. Thus, every (compact) crystallographic $n$-\textit{orbifold} admits a covering of degree at most $m'$ by a torus, where $m'$ depends only on the dimension.

\end{rem}

\section*{Further Reading}
There are many other applications of the Margulis lemma. The most famous is the thick--thin decomposition of hyperbolic manifolds (see \cite[Section 4.5]{Thurston-3d}). For applications concerning growth of linear groups, we refer to \cite{MR2200587}.

\section{Some Algebraic Number Theory}

\subsection{The \texorpdfstring{$p$}{p}-adic Numbers \texorpdfstring{$\mb Q_p$}{}}
We state a few properties of the field of $p$-adic numbers, $\mb Q_p$, and finite extensions of it. It is good to keep in mind the more basic example of $\mbm{k}=\mb Q_p$. 

A \textit{number field} is a finite extension $\mb K$ of $\mb Q$. A map $\left|\cdot\right|:\mb K\to[0,\infty)$
is called an \textbf{\emph{absolute value }}if it satisfies the following:
\begin{enumerate}
\item For $x\in \mb K$, $\left|x\right|=0$ if and only if $x=0$,
\item For $x,y\in \mb K$, $\left|xy\right|=\left|x\right|\cdot\left|y\right|$,
and
\item For $x,y\in \mb K$, $\left|x+y\right|\leqslant x+y$. 
\end{enumerate}
It follows that $\left|1\right|=1$ for every absolute value $\left|\cdot\right|$.
An absolute value defines a metric $d$ on $\mb K$, by setting $d(x,y)=\left|x-y\right|$,
and its completion $\mbm k$ is a locally compact field. An absolute value is \textbf{\textit{non-Archimedean}} if it is an ultra-metric, i.e., if it satisfies the equation $\l|x+y\r|\le \max\l\{\l|x\r|,\l|y\r|\r\}$ for every $x,y$. This is the case if and only if $\mathbb{Z}$ is bounded with respect to $\left|\cdot\right|$. In this case, $\mathcal{O}\coloneqq\left\{ x\in \mbm k\middle|\left|x\right|\leqslant1\right\} $
is a ring, called the \textbf{\emph{ring of integers}}. It is compact, and it is a maximal ring in $\mbm k$. 

We denote by $\mb Q_p$ the completion of $\mb Q$ with respect to the absolute value $\l|\cd\r|_p$, where $\l|\fr{m}{n}p^i\r|_p=p^{-i}$ when $m,n$ are co-prime to $p$. More generally, $\mbm k=\mb K_\nu$ is the completion of $\mb K$ with respect to some non-Archimedean absolute value $\l|\cd\r|_v$. The ring of integers of $\mb Q_p$ is called the \emph{\textbf{ring of $p$-adic integers}} and is denoted by $\mb Z_p$.

A uniformiser is an element $\pi$ in $\mc O_\mbm k$ with maximal absolute value which is less than $1$ (i.e.\ $\l|\pi\r|<1$ and if $\l|\alpha \r|<1$ then $\l|\alpha\r|\le \l|\pi\r|$). A uniformiser always exist. In the case $\mbm k=\mb Q_p$, $p$ is a uniformiser.

$\mc O_\mbm k$ is a local ring, i.e.\ it has a unique maximal ideal $m_\mbm k$. If $\pi$ is a uniformiser, then $m_\mbm k=\pi\mc O_\mbm k=\l(\pi\r)$. Every ideal of $\mc O_\mbm k$ is of the form $m_\mbm k^n=\l(\pi^n\r)$ for some $n\in\mb N$. In the case $\mbm k=\mb Q_p$, $m_\mbm k=p\mb Z_p$.

The field $\mc O_\mbm k/m_\mbm k\eqqcolon\mb F$ is finite. This is true because $\mc O_\mbm k$ is compact and $m_\mbm k$ is open in $\mc O_\mbm k$ 
(because of the ultra-metric equation), 
so the cosets in $\mc O_\mbm k/m_\mbm k$ are a disjoint open cover of $\mc O_\mbm k$, and hence must be finite. In the case $\mbm k=\mb Q_p$, $\mb F$ is the field with $p$ elements.\footnote{One may generalise our discussion to fields $\mbm k$ which aren't necessarily locally compact. In this case, the finiteness of $\mb F$ is equivalent to the local compactness of $\mbm k$.}

Elements $\alpha$ in $\mb Q_p$ are represented by
$$\alpha=\sum_{i\ge -N}a_ip^i,$$
where $a_i\in\{0,\dots,p-1\}$. Given this form for $\alpha$, we have $\l|\alpha \r|=p^n$ for the largest $n$ such that $a_{-n}\ne 0$.

A complex number is an \textit{\textbf{algebraic integer}} if it satisfies some monic polynomial with integer coefficients. The sum and the product of algebraic integers are algebraic integers.
If $\mb K$ is a number field, then its \textit{\textbf{ring of integers}} $\mc O_{\mb K}$ is the ring of algebraic integers contained in it. It is related to the ring of integers of completions of $\mb K$ in the following way: 
$$\mc O_{\mb K}=\mb K\cap \bigcap \mc O_\mbm k,$$ 
where $\mbm k$ ranges over all non-Archimedean local fields containing $\mb K$. That is,
$$\mc O_{\mb K}=\l\{\alpha\in\mb K\mi \l|\alpha\r|\le 1\;\forall\l|\cd\r|\text{ a non-Archimedean absolute value}\r\}.$$

\subsection{The Tree of \texorpdfstring{$\mr{SL}_2\l(\mb Q_p\r)$}{SL2Qp}}\label{SL_2-tree}

In this subsection we prove the following theorem:

\begin{thm}\label{sltq}
$\sltq$ acts on a $\l(p+1\r)$-regular tree. The stabilisers are compact, and the action admits two orbits.
\end{thm}

Let $V=\mb Q_p^2$. $\sltq$ acts on it by left multiplication.

We're going to use the following facts from commutative algebra:
\begin{enumerate}[(1)]
\i A \emph{\textbf{lattice}} of $V$ is a finitely generated $\mb Z_p$-submodule which generates $V$ over $\mb Q_p$. Every lattice $M$ of $V$ is of the form $M=\mb Z_p\cd e_1+\mb Z_p\cd e_2$ for some $e_1,e_2$ which form a basis of $V$ over $\mb Q_p$. This will be denoted by $M=\l(e_1,e_2\r)$.

\i For $M=\l(\bigl(\begin{smallmatrix}1\\0\end{smallmatrix}\bigr),\bigl(\begin{smallmatrix}0\\1\end{smallmatrix}\bigr)\r)$, $A\in\sltq$, we have
$$AM=M\iff A\in \mr{SL}_2\l(\mb Z_p\r).$$
\i Given two lattices $M,L$, there is a base $\l\{e_1,e_2\r\}$ of $V$ over $\mb Q_p$ such that $M=\l(e_1,e_2\r)$ and $L=\l(p^ae_1,p^be_2\r)$.
\end{enumerate}
You can try to prove these facts as an exercise.

We say two lattices $M_1,M_2$ are \textit{equivalent}, and denote $M_1\sim M_2$, if there is some $\alpha\in\mb Q_p$ such that $M_1=\alpha M_2$. We denote by $X$ the set of equivalence classes.

\begin{proof}[Proof of Theorem \ref{sltq}]
We define a metric on $X$. If $\overline M, \overline L\in X$ and $M\in \overline M$, $L\in\overline L$ are representatives, then set 
$$d\l(\overline M,\overline L\r)=\l|a-b\r|$$
for $a,b$ as above for $M,L$. The function $d$ is a well defined metric.

Another way of looking at $d$ is the following:
given such $M,L$, there exists $n_1$ such that $p^{n_1} L\ss M$ but $p^{n_1-1}L\nsubseteq M$ and there exists $n_2$ such that $M\ss p^{n_2}L$ but $M\nsubseteq p^{n_2+1}L$. Then $d\l(\overline M,\overline L\r)=\l|n_1-n_2\r|$ (and in fact $\l\{n_1,n_2\r\}=\l\{a,b\r\}$). This $d\ceq d\l(\overline M,\overline L\r)$ is the minimal number such that $\exists n\in\mb Z$ such that $p^{n+d}L\ss M\ss p^n L$.

We claim $\l(X,d\r)$ is a $\l(p+1\r)$-regular tree. We denote $\overline L-\overline M$ if $\overline M,\overline L$ are neighbours (i.e.\ if $d\l(\overline M,\overline L\r)=1$). This is equivalent to the condition there are $L\in \overline L,M\in\overline M$ such that 
$$M\supsetneqq L \supsetneqq pM.$$
Consider the projection $M\to M/pM$. The module $M/pM$ is actually a $2$-dimensional vector space over $\mb Z_p/p\mb Z_p\cong \mb F_p$, and $L$ projects to a subspace of it. It can't be of dimension $0$ or $2$ because $L\ne M$, $L\ne pM$, so $L$ is a $1$-dimensional subspace (a ``line"). 

Thus, the neighbours of $\overline M$\ correspond to lines in $\mb F_p^2$, and there are exactly $p+1$ such lines (since they are of the form $\mb F_p\cd \l(a,b\r)$ which is either $\mb F_p\cd \l(1,b/a\r)$ for $b/a=0,\dots,p-1$ or $\mb F_p\cd\l(0,b\r)$), so $\overline M$ has precisely $p+1$ neighbours.

We now need to show $X$ is connected and that it admits no cycles (i.e., that it is a tree).

\paragraph{Connectedness}
Given $\overline M,\overline L$, take representatives $M\in\overline M$, $L\in\overline L$ such that $M\supseteq L$ and $\exists e_1,e_2$ such that $M=\l(e_1,e_2\r)$, $L=\l(e_1,p^ne_2\r)$ for some $n\in\mb N$, $\l\{e_1,e_2\r\}$ is a basis of $V$ over $\mb Q_p$. Then
$$\overline M=\overline {(e_1,e_2)}-\overline {(e_1,pe_2)}\cdots-\overline {(e_1,p^{n-1}e_2)}-\overline {(e_1,p^ne_2)}=\overline L,$$
is a path from $\overline M$ to $\overline L$ so $X$ is connected.

\paragraph{No cycles} First we prove the following lemma:
\begin{lem}
If $M_0,M_1,M_2$ are lattices with $\overline {M_0}\neq\overline{M_2}$ and $M_0-M_1-M_2$ in the sense that
\begin{align*}
M_0&\supseteq M_1\supseteq pM_0,
\\M_1&\supseteq M_2\supseteq pM_1,
\end{align*}
then $M_1\equiv M_2\pmod {pM_0}$
\end{lem}
\begin{proof}
We have $pM_1\ss pM_0\ss M_1$ and $pM_1\ss M_2\ss M_1$, so both $pM_0$ and $M_2$ project to lines in $M_1/pM_1$ as above. These lines must be distinct, because otherwise we would get (by the correspondence of neighbours of $M_1$ and lines in $M_1/pM_1\cong \mb F_p^2$ described above) that $\overline{M_0}=\overline{M_2}$, contradictory to our assumption. But distinct lines in $2$-dimensional vector spaces generate the entire space, so we have
$$pM_0/pM_1+M_2/pM_1=M_1/pM_1,$$
and thus
$$ pM_0+M_2=M_1,$$
i.e.
$$M_2\equiv M_1\pmod{pM_0}.$$
\end{proof}
We now prove there can be no cycles in $X$. Suppose $\overline{M_0}-\cdots-\overline{M_n}$ is a path in $X$ with no backtracking; we need to prove that $\overline{M_0}\ne \overline{M_n}$. 

Take $M_i\in\overline{M_i}$ such that
$$M_0-\cdots -M_n,$$
which means
$$M_{i-1}\supsetneqq M_i\supsetneqq pM_{i-1}\quad\text{for } i=1,\dots,n$$
and $\overline{M_i}\ne\overline{M_{i+2}}$ for $i=0,\dots,n-2$.

It's enough to show $pM_0\nsupseteq M_n$. We prove this by induction on $n$, the case $n=1$ being trivial. 

  By our lemma $M_n=M_{n-1}\pmod{pM_{n-2}}$, so in particular (since $M_0\supseteq M_{n-2}$) it holds $M_n= M_{n-1}\pmod{pM_0}$. By the induction hypothesis $pM_0\nsupseteq M_{n-1}$, so $pM_0\nsupseteq M_n$, as desired.

The stabiliser of $(e_1,e_2)$ (with respect to the action $\sltq\acts X$) is $\sltz$, which is a compact subgroup.

$\mr{GL}_2(\mb Q_p)$ acts transitively on $X$, but $\mr{SL}_2(\mb Q_p)$ does not: $\l(e_1,e_2\r)$ cannot be mapped to $\l(e_1,pe_2\r)$.
But it can be mapped to $\l(e_1,p^2e_2\r)\sim \l(p^{-1}e_1,pe_2\r)$ by $\bigl(\begin{smallmatrix}p &0\\0 &1/p\end{smallmatrix}\bigr)\in\sltq$. The stabiliser of $(e_1,pe_2)$ is also a compact subgroup of $\sltq$, and these are the two orbits of its action on $X$.
\end{proof}

\begin{rem}
A discrete subgroup acting freely on a tree is free. Indeed,  trees are contractible and hence simply connected, so if a group $\G$ acts on a tree $T$ freely (on vertices as well as edges), then $\G\cong \pi_1\l(\G\backslash T\r)$, which is free since  $\G\backslash T$ is a graph and hence homotopy equivalent to a bouquet of circles.

It follows that every discrete, torsion-free subgroup of $\sltq$ is free. This is true because, if $\G\le \sltq$ is discrete, then $\G$ acts on $X= T_{p+1}$, and the stabiliser of any point is finite (being compact and discrete), and hence trivial.
\end{rem}

\begin{exe}
Let $\mbm k$ be a finite extension of $\mb Q_p$. Show that $\mr{GL}_2(\mbm k)$ acts transitively on a regular tree with compact stabilisers. What is the regularity of the tree?
\end{exe}

For a more thorough analysis of the tree constructed above we refer to Chapter 2 of \cite{Trees}.

\subsection{The Adeles}

Let $\mb K$ be a number field, and let $\left\{ \left|\cdot\right|_{\nu}\right\} _{\nu\in N}$ be the
collection of all absolute values on $\mb K$ (up to equivalence). For
$\left|\cdot\right|_{\nu}$, let $\mbm k_{\nu}$ be the completion of $\mb K$
with respect to (the metric defined by) $\left|\cdot\right|_{\nu}$.
It is a local field, and if $|\cd|_\nu$ is non-Archimedean, then the ring of integers $\mathcal{O}_{\nu}$ of $\mbm k_\nu$ is compact. The \textbf{\emph{adele ring }}of $\mb K$ is the restricted
product
\[
\mathbb{A}_{\mb K}=\tilde{\prod}_{\nu\in N}(\mbm k_{\nu},\mathcal{O}_{\nu})=\left\{ (x_{\nu})\in\prod_{\nu\in N}\mbm k_{\nu}\middle|x_{\nu}\in\mathcal{O}_{\nu}\text{ for almost every }\nu\in N\right\} .
\]
Since there
are only finitely many `infinite places' (i.e., Archimedean absolute
values on $\mb K$, up to equivalence), we don't need to worry about $\mathcal{O}_{\nu}$
in these cases.

We define multiplication and addition in $\mb A_{\mb K}$ componentwise.  We topologise it as follows: the basic open sets are of the form $\prod V_\nu$, where $V\nu\subseteq \mbm k_\nu$ is open for every $\nu\in N$ and $V_\nu=\mc O_\nu$ for almost every $\nu\in N$. These operations and topology make $\mb A_{\mb K}$ into a locally compact topological ring.

For $x\in \mb K$ we have $|x|_\nu=1$ (and hence $x\in\mathcal{O}_{\nu}$) for almost every $\nu$. This is easy to see in the case $\mb K=\mb Q$, since every rational number is composed of only finitely many primes, and its $p$-norm with respect to every other prime $p$ is exactly $1$. Therefore, we can embed $\mb K$ diagonally in $\mb A_{\mb K}$ (and we do).

Note that, if $|\cdot|_\nu$ is a non-Archimedean absolute value, then also $|\cdot|^t_\nu$ is an absolute value for every $t>0$. Two such absolute values are called \textit{equivalent}. We normalize the absolute values so that, if $\pi\in \mbm k_\nu$ is a uniformiser, then $|\pi |_\nu=[\mathcal{O}_{\nu}:\pi \mathcal{O}_{\nu}]^{-1}$. For instance, in $\mb Q_p$ we have $|p|=1/p$. 
With respect to this normalisation, we have the following:

\begin{fact}[The Product Formula]\label{product-formula}
    If $x\in\mb K\setminus\{0\}$, then $\prod_\nu |x|_\nu=1$.
\end{fact}
\begin{exe}
    Prove Fact \ref{product-formula} in the case $\mathbb{K}=\mb Q$.
\end{exe}
Observe that the product $\prod_{\nu}|x|_\nu$ is well-defined, since $|x|_\nu=1$ for almost every $\nu$.

For a finite subset $S\subseteq N$ of valuations, we denote 
\[\mc O_{\mb K}(S)=\l\{(x_\nu)\in\mb A_{\mb K}\m|x_\nu\in\mc O_\nu \quad\forall \nu\notin S\r\}.\]
This is a subring of $\mb A_{\mb K}$.
\begin{cor}
    Let $S\subseteq N$ be a finite subset of valuations. The projection of $\mc O_{\mb K}(S)$ to $\prod_{\nu\in S}\mbm k_\nu$ is discrete.
\end{cor}
\begin{proof}
    Denote by $p:\mb A_{\mb K}\to \prod_{\nu\in S}\mbm k_\nu$ the projection. By the Product Formula (Fact \ref{product-formula}), if $(x_\nu)\in \mc O_{\mb K}(S)$, then $\prod_{\nu\in S}|x_\nu|_\nu\geqslant 1$. Therefore, the intersection of  $p(\mc O_{\mb K}(S))$ with the open neighbourhood $\{(x_\nu)_{\nu\in S}\in \prod_{\nu\in S}\mbm k_\nu|\ |x_\nu|_\nu <1\}$ of $0$ is trivial. Since $p(\mc O_{\mb K}(S))$ is a ring, it follows that it is discrete. 
\end{proof}

For instance, if $\mathbb{K}=\mathbb{Q}(\sqrt{2})$ and $S=\varnothing$,
then we have a discrete embedding of $\mathcal{O}_{\mathbb{K}}=\mathbb{Z}[\sqrt{2}]$
in $\mathbb{R}\times\mathbb{R}$. If $\mathbb{K}=\mathbb{Q}$ (so
that $\mathcal{O}_{\mathbb{K}}=\mathbb{Z}$) and $S=\left\{ 7\right\} $,
then we have a discrete embedding of $\mathcal{O}_{\mathbb{K}}(S)=\mathbb{Z}[\frac{1}{7}]$
in $\mathbb{R}\times\mathbb{Q}_{7}$.
\begin{cor}
\label{cor:unbdd}Let $\mb K$ be a finite extension of $\mathbb{Q}$,
and let $P\subseteq \mb K$ be an infinite subset such that $P\ss \mc O_\nu$ for all but finitely many absolute values $|\cd |_\nu$. (That is, there is a finite subset $S\ss N$ such that $|s|_\nu\le 1$ for every $\nu\in N\setminus S$ and every $s\in P$). Then there exists an
absolute value on $\mb K$ with respect to
which $P$ is unbounded.\footnote{Observe that we assumed $P$ is contained in $\mc O_\nu$ for all but finitely many $\nu$'s. It follows from the statement that it cannot be contained in all of them: for at least one $\nu$, $P$ is unbounded, and \textit{a fortiori} not contained in $\mc O_\nu$.}
\end{cor}
\begin{proof}
Consider the embedding of $\mb K$ in $\mb A_{\mb K}$. 
By the previous corollary, we
get that $p(P)$ is discrete in $\prod_{\nu\in S}\mbm k_\nu$.
Therefore, since it is infinite, its image in $\mbm k_{\nu}$ for some $\nu\in S$ must
be noncompact: otherwise, the image of $P$ in $\prod_{\nu\in S}\mbm k_\nu$ would
be contained in the product of compact subsets, hence would be compact itself.\end{proof}

For a short and elegant exposition of the Adeles we refer to the beginning of \cite{arithmetic-groups}.

\section{Zimmer's Theorem}

The Zimmer Program is concerned with the phenomenon that groups with highly rigid structure do not act on manifolds of small dimension. One of the first results in this direction states that an infinite group with property (T) cannot act faithfully on $\mathbb{S}^{2}$ by conformal
maps. Since the group of conformal automorphisms of $\mathbb{S}^{2}$ is $\mathrm{PSL}(2,\mathbb{C})$ (acting by M\"obius transformation on the Riemann sphere), we can formulate this result as follows:

\begin{thm}
\label{thm:Zimmer}If a subgroup $\Gamma$ of $\mathrm{SL}(2,\mathbb{C})$
has property (T) as an abstract group, then it is finite.
\end{thm}

The following lemma is very useful:

\begin{lem}
\label{fact:algb-dens}Let $X$ be an algebraic variety defined over
$\mathbb{Q}$, and denote by $X(\mathbb{C})$ and $X(\bar{\mathbb{Q}})$
its $\mathbb{C}$-points and its $\bar{\mathbb{Q}}$-points respectively
(where $\bar{\mathbb{Q}}$ is the algebraic closure of $\mathbb{Q}$),
endowed with the Hausdorff topology (as subspaces of $\mathbb{C}^{N}$
for suitable $N$). Then $X(\bar{\mathbb{Q}})$ is dense in $X(\mathbb{C})$
in this topology.
\end{lem}

To prove this lemma, one can suppose that $X\subseteq \mb C^m$ and that $X$ is irreducible and defined by $k=m-\dim X$ polynomials. Given any point in $X$ one can approximate $k$ of its coordinates by rationals, and use the implicit function theorem to `correct' the other $m-k$ coordinates. We refer to \cite[Lemma 3.2]{uti} for a detailed proof.

\begin{proof}[Proof of Theorem \ref{thm:Zimmer}]
Assume by contradiction that $\Gamma$ is infinite. Since $\Gamma$ has property (T), it is finitely generated. Let $\Sigma$ be a finite set of
letters and $R$ a set of relations such that $\Gamma\cong\left\langle \Sigma\middle|R\right\rangle $,
and identify the elements of $\Sigma$ with the corresponding elements
of $\Gamma$, $\Sigma=\left\{ \gamma_{1},\dots,\gamma_{\ell}\right\} $.
Denote $G=\mathrm{SL}(2,\mathbb{C})$, and consider $\mathrm{Hom}(\Gamma,G)$,
the space of all homomorphisms from $\Gamma$ to $G$. We have an
injective map
\[
\mathrm{Hom}(\Gamma,G)\hookrightarrow G^{\ell},
\]
sending $f\in\mathrm{Hom}(\Gamma,G)$ to $(f(\gamma_{1}),\dots,f(\gamma_{\ell}))$.
Thus, we may think of $\mathrm{Hom}(\Gamma,G)$ as a subset of $G^{\ell}$.
The inclusion $\iota_0:\Gamma\hookrightarrow G$
is identified with $\left(\gamma_{1},\dots,\gamma_{\ell}\right)$.
A point $(g_{1},\dots,g_{\ell})\in G^{\ell}$ belongs to $\mathrm{Hom}(\Gamma,G)$
if, and only if, the tuple $(g_{1},\dots,g_{\ell})$ satisfies every
relation in $R$. But $G$ is an algebraic variety, and group relations
are algebraic formulas,
so $\mathrm{Hom}(\Gamma,G)$ is an algebraic
subvariety of $G^{\ell}$, which is moreover defined over $\mathbb{Q}$.
Denote this variety by $X$, so that $X(\mathbb{C})=\mathrm{Hom}(\Gamma,G)$.
Then, by Lemma \ref{fact:algb-dens}, $X(\bar{\mathbb{Q}})$ is dense
in $\mathrm{Hom}(\Gamma,G)$. Therefore, every neighbourhood of the
inclusion $\iota_0:\Gamma\hookrightarrow G$ admits a point, which is
a homomorphism from $\Gamma$ to $G$, whose
image lies in $\mathrm{SL}(2,\bar{\mathbb{Q}})$.

\begin{lem}\label{small-def}
If $\iota\in\mathrm{Hom}(\Gamma,G)$ is close enough to the inclusion
$\iota_0:\Gamma\hookrightarrow G$, then the image of $\iota$ in
$G$ is infinite. In other words, there is a neighbourhood $V$ of
the inclusion $\iota_0$ such that every homomorphism in $V$ has infinite
image.
\end{lem}
\begin{proof}[Proof of Lemma \ref{small-def}]
By the Jordan theorem, there is some $m\in\mathbb{N}$ such that
every finite subgroup of $\mathrm{SL}(2,\mathbb{C})$ admits an abelian subgroup of index at most $m$. Since $\Gamma$ is finitely
generated, it admits only finitely many subgroups of index at most
$m$,
so
\[
\Delta\coloneqq\bigcap_{N\leqslant\Gamma,[\Gamma:N]\leqslant m}N
\]
is of finite index in $\Gamma$. In particular, $\Delta$ has property (T) and
it is infinite, and hence is not abelian. Pick $\delta_{1},\delta_{2}\in\Delta$
such that $\left[\delta_{1},\delta_{2}\right]=\delta_{1}\delta_{2}\delta_{1}^{-1}\delta_{2}$
is nontrivial. There are words $w_{1},w_{2}$ on $\ell$ letters such
that $w_{i}(\gamma_{1},\dots,\gamma_{\ell})=\delta_{i}$ for $i=1,2$.
Set $w_{0}=[w_{1},w_{2}]$. By continuity, there is a neighbourhood
$V$ of the inclusion $\iota_0:\Gamma\hookrightarrow G$ such that,
if $i_{1}\in V$, then
\[
w_{0}(i_{1}(g_{1}),\dots,i_{1}(g_{\ell}))\neq1.
\]
In particular, $i_{1}(\Delta)$ is not abelian. But every subgroup
of index at most $m$ of $i_{1}(\Gamma)$ contains $i_{1}(\Delta)$ (to see this, look at the preimage in $\G$). Thus, no subgroup
of index at most $m$ of $i_{1}(\Gamma)$ is abelian. In view of Jordan's theorem, this implies  that $i_{1}(\Gamma)$
is infinite.
\end{proof}

Therefore, there is $i_{1}:\Gamma\to\mathrm{SL}(2,\bar{\mathbb{Q}})$
whose image is infinite. Denote $\tilde{\Gamma}=i_{1}(\Gamma)$.
$\tilde{\Gamma}$ is generated by the finite subset $i_1(\Sigma)$; let $F$ be the collection of all entries of all matrices in $i_1(\Sigma)$. Since $F\subseteq\mathbb{Q}$ is finite, it is contained in some finite extension $\mathbb{K}$ of $\mathbb{Q}$. Then $\mathbb{K}$ is a number field, and  $\tilde{\Gamma}\subseteq\mathrm{SL}(2,\mb K)$. Moreover, letting $A\subset \mathbb{K}$ denote the ring generated by $F$, we see that $\G\le \mathrm{SL}(2,A)$.

Consider $\tilde{\Gamma}$ as a subset of $\mb K^{4}$. 
Since $F$ is finite, $A$ is finitely generated, and there are only finitely many valuations $v$ such that the entries of $\mathrm{SL}(2,A)$ are not contained in the ring of $\mbm k_v$-integers. 
Let $S$ denote this finite set of valuations. As $\tilde\G$ is infinite, it follows from Corollary \ref{cor:unbdd} that, for some $v\in S$, the infinite group $\tilde{\Gamma}$ is unbounded in $\mathrm{SL}(2,\mbm k_v)$. 

The field $\mbm k_v$ is either $\mathbb{R},~\mathbb{C}$ or a finite extension of $\mathbb{Q}_{p}$.
Since $\tilde \G$ is unbounded in $\mathrm{SL}(2,\mbm k_v)$ we get that it acts on $\mathbb{H}^2$, $\mathbb{H}^3$ or the corresponding tree without fixed points. However, we have seen that property (T) implies fixed points for actions on hyperbolic spaces and trees. 
Therefore, this is a contradiction,
so $\G$ must have been finite all along.\end{proof}
\chapter{The Tits Alternative}

We have already seen in the first chapter the importance of the existence of free subgroups. The celebrated Tits alternative establishes the following dichotomy: every linear group either admits a free subgroup, or is virtually solvable. This result has numerous applications in various fields of mathematics. This section is devoted to it proof.

\begin{thm}[The Tits alternative]\label{thm:Tits}
Let $\Gamma$ be a finitely generated subgroup
of $\mathrm{GL}_{n}(\mathbb{C})$. Either $\Gamma$ is virtually solvable,
or it contains a non-abelian free subgroup.
\end{thm}

The assumption that $\G$ is finitely generated can be removed. One can actually show that, if $\G\le \mathrm{GL}_{n}(\mathbb{C})$  is any subgroup which is not virtually solvable, then it admits a finitely generated subgroup which is not virtually solvable. 

Theorem \ref{thm:Tits} holds for every field; for fields of positive characteristic, however, the finite generation assumption is required.
For example, $\G=\mathrm{GL}_{n}(\overline{\mathbb{F}}_p)$, where $\overline{\mathbb{F}}_p$ is the algebraic closure of the field with $p$ elements, is not virtually solvable and does not contain a non-abelian free subgroup (in fact, it is locally finite).
Although the proof in the general case is virtually the same, we suppose that the field is $\mathbb{C}$.

The idea of the proof is as follows: by taking small deformations and quotients, one may reduce to the case $\G$ is a Zariski dense, unbounded subgroup of an algebraic group (over some local field) acting strongly irreducibly on a vector space. By analysing the dynamics of the action of linear groups on the projective space, we will be able to construct a pair of elements playing ping pong, and thus generating a free subgroup.

The proof we present is taken from \cite{breuillard2003dense,breuillard2007topological}, in which stronger results are established.
\section{Reductions}

We first reduce the proof to the following case:
\begin{thm}\label{thm:dyn}
Let $\mbm k$ be a local field, and let $\G\le\mr{SL}_d(\mbm k)$ be a subgroup whose Zariski closure in $\mr{SL}_d(\mbm k)$ is Zariski connected and acts irreducibly on $\mbm k^d$. If $\G$ is unbounded with respect to the norm of $\mbm k$, then it admits a non-abelian free subgroup.
\end{thm}

Let us explain why Theorem \ref{thm:dyn} implies Theorem \ref{thm:Tits}.
Let $\G$ be a finitely generated subgroup of $\mr{GL}_n(\mb C)$ that is not virtually solvable.

First, let us make the following useful observation:
if $f:\Gamma\to\Delta$ is an epimorphism and $\Delta$ contains
a non-abelian free subgroup, then so does $\Gamma$.
Thus, it is enough to prove that some quotient of $\Gamma$ contains
a non-abelian free subgroup. Observe that, if $f:\Gamma\to\Delta$
is an epimorphism, then clearly $\Delta$ is also finitely generated.

Using this observation, we will now make a few reductions, similar to the ones we made in the proof of Zimmer's theorem.
Recall that $\mathrm{GL}_{n}(\mathbb{C})$ is an algebraic group,
and admits a Zariski topology. When referring to topological properties
of this topology, we prefix them by `Zariski' (e.g., we say a subset
is `Zariski dense' if it is dense in the Zariski topology).

Set $H=\overline{\Gamma}^{Z}$, the Zariski closure of $\Gamma$ in
$\mathrm{GL}_{n}(\mathbb{C})$. Since $H$ is an algebraic group,
its connected component $H^{\circ}$ is of finite index in it ($[H:H^{\circ}]<\infty$),
and therefore $H^{\circ}$ is open in $H$. This means that $\Gamma_{1}\coloneqq\Gamma\cap H^{\circ}$
is Zariski dense in $H^{\circ}$, and is still finitely generated.

Recall that every connected Lie group $L$ admits a maximal normal solvable subgroup $\text{Rad}(L)$ called the radical, and the quotient $L/\text{Rad}(L)$ is a centre-free semisimple Lie group, hence a direct product of simple Lie groups.
Let $R$ be the solvable radical of $H^{\circ}$, so that $H^{\circ}/R$
is semisimple. If $R=H$ then $H$ is solvable, and hence $\Gamma_{1}$
is solvable, so that $\Gamma$ is virtually solvable, contrary to
our assumptions. 
Therefore, $H^{\circ}/R$ is a nontrivial centre-free semisimple
group. 
Thus, there is some nontrivial, centre-free, simple algebraic group
$G$ and an epimorphism $f:H^{\circ}\to G$. Since $\Gamma_{1}$ is
Zariski dense in $H^{\circ}$, it follows that $f(\Gamma_{1})$ is
Zariski dense in $G$. 
Abusing notations, we
denote $\Gamma=f(\Gamma_{1})$.

Let $S$ be a finite set of letters and $R$ a set of relations such
that $\Gamma\cong\left\langle S\middle|R\right\rangle $, and identify
the elements of $S$ with the corresponding elements of $\Gamma$,
$S=\left\{ \gamma_{1},\dots,\gamma_{\ell}\right\} \subseteq\Gamma$.
Consider $\mathrm{Hom}(\Gamma,G)$, the space of all homomorphisms
from $\Gamma$ to $G$. We have an injective map 
\[
\mathrm{Hom}(\Gamma,G)\hookrightarrow G^{\ell},
\]
sending $f\in\mathrm{Hom}(\Gamma,G)$ to $(f(\gamma_{1}),\dots,f(\gamma_{\ell}))$.
As we saw in the previous chapter, this identifies $\mr{Hom}(\G,G)$ with an algebraic subvariety of $G^\ell$, which is defined over $\mb Q$. This is a complex variety, so let $X$ be such that $X(\mathbb{C})=\mathrm{Hom}(\Gamma,G)$.
By Lemma \ref{fact:algb-dens}, $X(\bar{\mathbb{Q}})$ is dense
in $\mathrm{Hom}(\Gamma,G)$ with respect to the usual (complex) topology. Therefore, every neighbourhood of the
inclusion $\iota_0:\Gamma\hookrightarrow G$ admits a point (where `points'
here are homomorphisms from $\Gamma$ to $G$) whose image lies in
the $\bar{\mathbb{Q}}$-points of $G$. Let us denote the $\bar{\mathbb{Q}}$-points
of $G$ by $G(\bar{\mathbb{Q}})$.\footnote{Technically speaking, we ought to have let $\mathbb{G}$ be the algebraic
group for which $G=\mathbb{G}(\mathbb{C})$, so that the $\bar{\mathbb{Q}}$-points
are $\mathbb{G}(\bar{\mathbb{Q}})$.}
\begin{lem}
There is a neighbourhood $V$ of the inclusion $\iota_0\in\mathrm{Hom}(\Gamma,G)$
such that every homomorphism in $V$ has a Zariski dense image.
\end{lem}

\begin{proof}
    Let $\mathfrak{g}$ be the Lie algebra of $G$ and let $\text{Ad}:G\to\mathrm{GL}(\mathfrak{g})$ be the adjoint representation.
    Consider $\mathcal{A}\coloneqq \mr{Span}(\mr{Ad}(G)) \subseteq\mathrm{End}(\mathfrak{g})$. It is closed under composition, and so it is an associative subalgebra of $\mr{End}(\mf g)$. Since $\mathfrak{g}$ has no nontrivial ideals, $\mathcal{A}$ admits no nontrivial invariant subspaces. In other words, the adjoint representation is irreducible. Since $\G=\langle\gamma_1,\ldots,\gamma_\ell\rangle$ is Zariski dense and $\mathcal{A}$ is finite dimensional, there are $d=\dim \mathcal{A}$ words, $W_1,\ldots,W_d\in F_\ell$, such that 
    \[
     \mathcal{A}=\mr{Span}\{ \mr{Ad}(W_1(\underline{\gamma})),\ldots,\mr{Ad}(W_d(\underline{\gamma}))\},
    \] 
    where $\underline{\gamma}=(\g_1,\dots,\g_\ell)$. For $\iota\in\mr{Hom}(\G,G)$, we have
    \[
    \mr{Span}\{ \mr{Ad}(\iota(W_1(\underline{\gamma}))),\ldots,\mr{Ad}(\iota(W_d(\underline{\gamma})))\}=
    \mr{Span}\{ \mr{Ad}(W_1(\iota(\underline{\gamma}))),\ldots,\iota(\mr{Ad}(W_d(\iota(\underline{\gamma})))\},
    \]
    where $\iota(\underline{\gamma})=(\iota(\g_1),\dots,\iota(\g_\ell))$. If $\iota\in\mr{Hom}(\G,G)$ is sufficiently close to the inclusion $\iota_0$, then 
    \[
    \mr{Span}\{ \mr{Ad}(W_1(\iota(\underline{\gamma}))),\ldots,\mr{Ad}(W_d(\iota(\underline{\gamma})))\}=\mc{A}.
    \]
    
 For $\iota\in\mr{Hom}(\G,G)$, consider the closure of $\iota(\G)$ in the Zariski topology, and denote by $\mathfrak{n}_{\iota}\le \mathfrak{g}$ its Lie algebra. For every $\iota$, the Lie algebra $\mf n_\iota$ is $\mr{Ad}(\iota(\G))$-invariant, and hence invariant under
\[
    \mr{Span}\{ \mr{Ad}(W_1(\iota(\underline{\gamma}))),\ldots,\mr{Ad}(W_d(\iota(\underline{\gamma})))\}.
\]
Thus, if $\iota$ is sufficiently close to the inclusion $\iota_0$, then $\mathfrak{n}_{\iota}$ is $\mc {A}$-invariant, hence an ideal in $\mathfrak{g}$. Since $\mf g$ is simple, we get that $\mf n_\iota$ is either $0$ or $\mathfrak{g}$, and that $\iota(\G)$ is either finite or Zariski dense. By Lemma \ref{small-def}, if $\iota$ is sufficiently close to the inclusion, then $\iota(\G)$ is infinite, hence Zariski dense.
\end{proof}

Therefore, there is $\hat{\iota}:\Gamma\to G(\bar{\mathbb{Q}})$
such that $\hat{\iota}(\Gamma)$ is Zariski dense in $G$. We abuse
notations again and denote $\Gamma=\hat{\iota}(\Gamma)$. Since $\Gamma$
is finitely generated, there is a finite extension $\mathbb{K}$ of
$\mathbb{Q}$ such that $\Gamma\subseteq G(\mathbb{K})$ (since the
collection of all the entries of all the matrices in the generating
set of $\Gamma$ is a finite subset of $\bar{\mathbb{Q}}$, and every finite
subset of $\bar{\mathbb{Q}}$ is contained in a finite extension of $\mathbb{Q}$).

Now, since $\Gamma$ is infinite, there is (by Corollary \ref{cor:unbdd}) some valuation $\nu$ on $\mb K$ such
that $\Gamma$ is unbounded with respect to the absolute value corresponding
to $\nu$.
Denote $\mbm k=\mathbb{K}_{\nu}$, the completion of $\mathbb{K}$
with respect to $\left|\cdot\right|_{\nu}$, which is a local field.
We get that $\Gamma$ is an unbounded subgroup of $G(\mbm k)$, and it is still Zariski dense in it. 

Since $G(\mbm k)$ is linear, we can embed it in $\mathrm{SL}(V)$
for some vector space $V$ over $\mbm k$. 
Since $G$ is simple, it is completely reducible, and $V=\bigoplus_{j=1}^{r}W_{j}$,
where the action of $G$ on $W_{j}$ is irreducible for every $j$.
For $g\in\mathrm{SL}(V)$, write $g=\bigoplus_{j}g_{j}$, with each
$g_{j}$ acting on $W_{j}$. We have $\left\Vert g\right\Vert \leqslant\sum\left\Vert g_{j}\right\Vert $,
so,
there is some $j$ such that 
the image of $G(\mbm k)$ in $\mr{SL}(W_j)$ is still unbounded, and acts irreducibly.
This completes the reduction from Theorem \ref{thm:Tits} to Theorem \ref{thm:dyn}.

\section{The Cartan Decomposition}
We will use the Cartan decomposition (or $KAK$-decomposition) of
$\mathrm{SL}_{d}(\mbm k)$, where $\mbm k$ is some local field. We define $K$ and $A$ as follows:
\begin{enumerate}
\item If $\mbm k=\mathbb{R}$, we let $K=\mathrm{SO}_{d}(\mathbb{R})$; if $\mbm \mbm k=\mathbb{C}$,
we let $K=\mathrm{SU}_{d}(\mathbb{C})$; if $k$ is non-Archimedean,
we let $K=\mathrm{SL}_{d}(\mathcal{O}_{\mbm k})$.
\item If $\mbm k$ is Archimedean, we let
\[
A=\left\{ \mathrm{diag}(\alpha_{1},\dots,\alpha_{d}):\alpha_{1}\geqslant\alpha_{2}\geqslant\cdots\geqslant\alpha_{d}>0,\quad\prod_{i=1}^{d}\alpha_{i}=1\right\} .
\]
If $\mbm k$ is non-Archimedean, we let
\[
A=\left\{ \mathrm{diag}(\pi^{n_{1}},\dots,\pi^{n_{d}}):n_{i}\in\mathbb{Z},n_{1}\leqslant\cdots\leqslant n_{d},\sum_{i=1}^{d}n_{i}=0\right\} ,
\]
where $\pi$ is a uniformiser for $\mbm k$ (for $\mbm k=\mathbb{Q}_{p}$, one
can take $\pi=p$).
\end{enumerate}
If $g\in\mathrm{SL}_{d}(\mbm k)$, then there are $k_{g}^{1},k_{g}^{2}\in K$
and $a_{g}\in A$ such that $g=k_{g}^{1}a_{g}k_{g}^{2}$. The $k_{g}^{1},k_{g}^{2}$
are not necessarily unique, but $a_{g}$ is.

\begin{notation}
    Let $d\in \mb N$ and let $\mbm k$ be a local field. For $g\in\mr{SL}_d(\mbm k)$, we denote by $a_g=a(g)$ the unique element in $A$ such that $g=k_1a_gk_2$ for some $k_1,k_2\in K$. We denote by $\alpha_1(g),\dots,\alpha_d(g)\in \mbm k$ the scalars for which
    \[
    a_g=\mr{diag}(\alpha_1(g),\dots,\alpha_d(g)).
    \]
\end{notation}

Note that, with respect to the standard operator norm, $\|g\|=\|a_g\|=|\alpha_1(g)|$. 
\section{Dynamics of Projective Transformations}
In this section, we investigate the dynamics of linear groups acting on their corresponding projective spaces, in order to finish the proof of the Tits alternative.

\subsection{Unboundedness}
We first wish to investigate the properties of a sequence of linear transformations going to infinity. 
\begin{lem}\label{lem:unboundedness}
    Let $d\in \mb N$, $\mbm k$ a local field. Let $\gamma_{m}\in\mathrm{SL}_{d}(\mbm k)$ such that $\left\Vert \gamma_{m}\right\Vert _{\mbm k}\longrightarrow\infty$. Then there is some $\ell=1,\dots,d$ such that the following holds: denoting by $\tilde \g_m$ the image of $\g_m$ in $\mr{SL}(\bigwedge ^\ell (\mbm k^d))\cong \mr{SL}_{\binom{d}{\ell}}(\mbm k)$ and  considering the $KAK$-decomposition of this group, we have 
    \[
    \l|\frac{\alpha_1(\tilde \g_m)}{\alpha_2(\tilde \g_m)}\r|\longrightarrow\infty
    \]
\end{lem}

In the proof, we will use the \textit{exterior power} of a vector space. If $V$ and $U$ are vector spaces over a field $k$, then an $n$-linear map $A:V^n\to U$ is called \textit{alternating} if, for every permutation $\sigma\in S_n$, one has
\[
    A(v_1,\dots,v_n)=\mr{sgn}(\sigma)\cdot A(v_{\sigma(1)},\dots,v_{\sigma(n)}).
\]
For every vector space $V$ and $n\in \mb N$, the exterior $n^\text{th}$ power of $V$ is a vector space $\bigwedge ^n V$, equipped with an alternating map $V^n\to \bigwedge^n V$, which is denoted by
\begin{align*}
    V\times\cdots\times V&\to \bigwedge^n V\\
    (v_1,\dots,v_n) &\mapsto v_1\wedge\cdots\wedge v_n.
\end{align*}
The exterior power satisfies the following universal property: if $T:V^n\to U$ is an alternating map into some vector space $U$, then there is a linear map $\tilde T:\bigwedge ^n V\to U$ such that \[\tilde T(v_1\wedge\cdots\wedge v_n)=T(v_1,\cdots,v_n)\] for every $v_1,\dots,v_n$. 
The exterior power is in fact the unique vector space admitting a map satisfying this property.
If $V$ is $d$-dimensional, then $\bigwedge^n V$ is $\binom{d}{n}$-dimensional; in particular, it is zero dimensional if $n>d$. If $e_1,\dots,e_d$ is a basis of $V$, then the collection of all vectors of the form $e_{i_1}\wedge\cdots\wedge e_{i_n}$ for $1\le i_1<\cdots<i_n\le d   $ is a basis of $\bigwedge^n V$. 

\begin{proof}
Denote $a_{m}=a(\g_m)$ and $\alpha_{i,m}=\alpha_i(\g_m)$ for every $m\in\mathbb{N}$, with respect to the $KAK$-decomposition of $\mr{SL}_d(\mbm k)$.
Then $|\alpha_{1,m}|\longrightarrow \infty$ as $m$ goes to infinity, and also $\prod a_{i}=1$. Thus, up to taking a  subsequence, there
is $\ell=1,\dots,d$ such that $|\alpha_{\ell,m}/\alpha_{\ell+1,m}|\longrightarrow\infty$
(as $m$ goes to $\infty$). Fix this $\ell$, and set $V\coloneqq\bigwedge^{\ell}(\mbm k^{d})$. 

For every $m\in\mb N$, we denote by $\tilde{\gamma}_{m}$ the image of $\gamma_{m}$ in $\mathrm{SL}(V)$, and we denote $\tilde a_m=a(\tilde \g_m)$, $\tilde \alpha_{i,m}=\alpha_i(\tilde\g_m)$ (for $i=1,\dots,d$) with respect to the $KAK$-decomposition of $\mr{SL}(V)\cong\mr{SL}_{\binom{d}{\ell}}(\mbm k)$. 

We have $\tilde{\alpha}_{1,m}=\alpha_{1,m}\cdots\alpha_{\ell,m}$ and $\tilde{\alpha}_{2,m}=\alpha_{1,m}\cdots\alpha_{\ell-1,m}\alpha_{\ell+1,m}$,
since these are the two biggest eigenvalues of $\tilde{a}_m$ (since
$\alpha_{1,m},\dots,\alpha_{\ell,m}$ are the $\ell$ biggest eigenvalues of
$a_m$, and $\alpha_{\ell+1,m}$ is the next one).
Thus, 
 $|\tilde{\alpha}_{1,m}/\tilde{\alpha}_{2,m}|\longrightarrow\infty$
(as $m$ goes to infinity), since $|\tilde{\alpha}_{1,m}/\tilde{\alpha}_{2,m}|=|\alpha_{\ell,m}/\alpha_{i\ell+1,m}|$,
and our assumption regarding $\ell$ above is exactly that this goes
to infinity.
\end{proof}

Therefore, by replacing the linear space on which $\G$ acts, we may assume that, if it is unbounded, then it contains elements $\g$ such that $|\alpha_1(\g)/\alpha_2(\g)|$ is as large as we wish.

\subsection{Contracting Transformations}
We fix a vector space $W$ over some local field $\mbm k$. 
Recall that we denote by $\mb P(W)$ the corresponding projective space.

We now investigate the question: what does
it mean that $|\alpha_{1}(g_m)/\alpha_{2}(g_m)|\longrightarrow\infty$  for $g_m\in\mr{SL}(W)$?

\def \dw {d}

We fix some basis $e_1,\dots,e_\dw$ of $W$, and consider the associated norm on $W$. In the Archimedean case, the norm is induced from the natural inner product corresponding to the basis; in the non-Archimedean case, there is no concept of an inner product, and one simply sets $\|\sum_{i=1}^\dw v_i e_i\|=\max_{i}\{|v_i|\}$. We consider the $KAK$-decomposition of $\mr{SL}(W)$: 
\[
\mr{SL}(W)=\mr{SL}_\dw(\mbm k)={K}{A}{K}.
\]
 Observe that $K$ preserves the norm of $W=\mbm k^\dw$ (and the inner product, in the Archimedean case). 

We will use the following natural metric on $\mb P(W)$. Consider the exterior product $\bigwedge^2W$; we have a canonical norm on $\bigwedge^2W$. In the Archimedean case, one sets
\[
    \langle u_1\wedge u_2,v_1\wedge v_2\rangle=\det \begin{pmatrix}\left\langle u_{1},v_{1}\right\rangle  & \left\langle u_{1},v_{2}\right\rangle \\
\left\langle u_{2},v_{1}\right\rangle  & \left\langle u_{2},v_{2}\right\rangle 
\end{pmatrix};
\]
in the real case, this means that $\|u\wedge v\|=\|u\|\cdot\|v\|\cdot |\sin \alpha|$, where $\alpha$ is the angle between $u$ and $v$. In the non-Archimedean case, we take the natural basis $(e_i\wedge e_j)_{i>j}$ on $\bigwedge^2W$, and take the sup-norm with respect to this basis. 
We then define the metric on $\mb P(W)$ as follows:
\[
    d([v],[w])=\frac{\|v\wedge w\|}{\|v\|\cdot \|w\|}.
\]
In the real case, we get that $d([v],[w])=|\sin \alpha|$, where $\alpha$ is the angle between $u$ and $v$.
In any case, $K$ preserves the metric on $\mb P(W)$. 

\begin{notation}
    We denote by $(\cdot)_{\varepsilon}$ the $\varepsilon$-neighbourhood of a set or a point.
\end{notation}
\begin{defn}
A projective transformation $g:\mathbb{P}(W)\to\mathbb{P}(W)$ is
called \emph{$\varepsilon$-contracting} if there is a projective
hyperplane $H\subseteq\mathbb{P}(W)$, and a projective point $v\in\mathbb{P}(W)$,
such that 
\[
    g(\mathbb{P}(W)\setminus (H)_{\varepsilon})\subseteq(v){}_{\varepsilon}.
\]
In this case
$v$ is called an \emph{attracting point} of $g$, $H$ a \emph{repelling
hyperplane}, and we say $g$ is $\varepsilon$-contracting with respect
to the pair $(v,H)$.
\end{defn}

\begin{exa}
For $W=\mb R^3$, the projective transformation corresponding to 
\[
\left(\begin{array}{ccc}
100\\
 & 1/10\\
 &  & 1/10
\end{array}\right)
\]
is contracting with respect to $(e_{1},\mathrm{Span}(e_{2},e_{3}))$.
\end{exa}

\begin{exe}\label{exe:K-contracting}
If $g$ is $\varepsilon$-contracting with respect to some $(v,H)$, and
$k_{1},k_{2}\in K$ (for the $K$ of the $KAK$-decomposition of $\mr{SL}(W)$), then $k_{1}gk_{2}$ is $\varepsilon$-contracting
with respect to $k_{1}v,k_{2}^{-1}H$.
\end{exe}

\begin{prop}\label{prop:contracting}
Let $g\in\mr{SL}(W)$. For every $\varepsilon>0$ there is $M\in\mathbb{N}$ such that, if
\[
|\alpha_{1}(g)/\alpha_{2}(g)|\geqslant M,
\]
then there is some hyperplane $H\subseteq \mb P(W)$ such that 
$g$ is $\ep$-Lipschitz outside $(H)_\ep$. In particular, $g$ is
$\varepsilon$-contracting.
\end{prop}

\begin{proof}
Since the metric is $K$-invariant, it is enough to consider a diagonal
$g\in A$,
\[
g=\mathrm{diag}(\alpha_{1},\dots,\alpha_{\dw}).
\]
It is straightforward to check that, if $|\alpha_{1}/\alpha_{2}|$
is large enough, then the result holds with respect to the hyperplane $H=\mathrm{Span}(e_{2},\dots,e_{d})$. We leave it as an exercise to the reader to determine $M(\ep)$ explicitly.
Obviously, $[e_1]$ is an attracting point.

Since the distance between any two points in $\mb P(W)$ is at most $1$, we get that $g$ is $\ep$-contracting with respect to $([e_1],H)$.
\end{proof}

\begin{cor}
Let $g_m\in\mr{SL}(W)$ for $m\in\mb N$. If $(g_{m})$ satisfies $\lim_{m\to\infty}|\alpha_{1}(g_{m})/\alpha_{2}(g_{m})|=\infty$,
then there are $\varepsilon_{m}\longrightarrow0$ such that $g_{m}$
is $\varepsilon_{m}$-contracting.
\end{cor}
 
\begin{lem}\label{Lip}
For every $g\in\mathrm{GL}(W)$ and every $\delta>0$ there is an open subset
$\Omega\subseteq \mb P(W)$ such that $g\!\restriction_{\Omega}$ is $(1+\delta)$-Lipschitz.
\end{lem}

\begin{proof}
Write $g=k_1ak_2$ using the $KAK$-decomposition. Let $a=\mr{diag}(\alpha_1,\dots,\alpha_{d}$).
It is not hard to check that $a$ is $(|\frac{\alpha_2}{\alpha_1}|+\delta)$-Lipschitz in a small enough neighbourhood of $[e_1]$, where $e_1=(1,0,\ldots,0)$.\footnote{Moreover, with respect to the natural basis, the Jacobean of $a$ at $[e_1]$ is $\mr{diag}(\alpha_2/\alpha_1,\dots,\alpha_{d}/\alpha_1)$.}
Therefore $g$ is $(|\frac{\alpha_2}{\alpha_1}|+\delta)$-Lipschitz in a small neighbourhood of $k_2^{-1}[e_1]$.
Since $|\alpha_2/\alpha_1|\le 1$, the result follows. 
\end{proof}
\begin{prop}\label{LipIsCont}
    For every $\ep>0$ there is $\delta>0$ such that, if $g\in \mr{SL}(W)$ is $\delta$-Lipschitz in some open subset $\Omega\subseteq \mb P(W)$, then it is $\ep$-contracting. 
\end{prop}
\begin{proof}
    Let $g\in\mr{SL}(W)$. 
    By Exercise \ref{exe:K-contracting}, $g$ is $\ep$-contracting if and only if $a(g)$ is $\ep$-contracting, so we may assume $g=a(g)=\mr{diag}(\alpha_1(g),\dots,\alpha_\dw(g))$.   
    Write $\alpha_i=\alpha_i(g)$. 
    
    By Proposition \ref{prop:contracting}, it is enough to show the following: for every $M>0$ there is $\delta>0$ such that, if $g$ is $\delta$-Lipschitz in some open subset $\Omega\subseteq\mb P(U)$, then $|\alpha_1/\alpha_2|>M$. This can be shown by a direct computation, as follows.

    Assume that $g\!\restriction_{\Omega}$ is $\delta$-Lipschitz for
    some open subset $\Omega\subseteq\mathbb{P}(W)$. We will now bound
    $\left|\alpha_{1}/\alpha_{2}\right|$ from below with respect to $\delta$,
    so that, by choosing $\delta$ to be small enough, we will get that
    $\left|\alpha_{1}/\alpha_{2}\right|$ is arbitrarily large. 
    
    Let $v=\sum_{i=1}^{\dw}v_{i}e_{i}$ be a nonzero vector such that
    $[v]\in\Omega$. If $x\in\mbm k$ and the norm of $x$ is small enough,
    then $[w_{1}]\coloneqq[v+xe_{1}]$ and $[w_{2}]\coloneqq[v+xe_{2}]$
    are both in $\Omega$. Since $g$ is $\delta$-Lipschitz, we have
    \[
    d([gv],[gw_{1}])\leqslant\delta d([v],[w_{1}]).
    \]
    Observe that $gv\wedge gw_{1}=gv\wedge gxe_{1}$ (since $gv\wedge gv=0$).
    Moreover, since $g$ is diagonal, we get that $ge_{1}=\alpha_{1}e_{1}$.
    Therefore, by definition, we have 
    \[
    d([gv],[gw_{1}])=\frac{\left\Vert gv\wedge gxe_{1}\right\Vert }{\left\Vert gv\right\Vert \cdot\left\Vert gw_{1}\right\Vert }=\left|x\alpha_{1}\right|\frac{\left\Vert gv\wedge e_{1}\right\Vert }{\left\Vert gv\right\Vert \cdot\left\Vert gw_{1}\right\Vert }.
    \]
    On the other hand, we have
    \[
    \delta d([v],[w_{1}])=\delta\left|x\right|\frac{\left\Vert v\wedge e_{1}\right\Vert }{\left\Vert v\right\Vert \cdot\left\Vert w_{1}\right\Vert }.
    \]
    Cancelling out $|x|$ and using the fact  $\left\Vert v\wedge e_{1}\right\Vert \leqslant\left\Vert v\right\Vert$,
    we get
    \begin{align}\label{Eq:ineq1}
            \left|\alpha_{1}\right|\frac{\left\Vert gv\wedge e_{1}\right\Vert }{\left\Vert gv\right\Vert \cdot\left\Vert gw_{1}\right\Vert }\leqslant\delta\frac{\left\Vert v\wedge e_{1}\right\Vert }{\left\Vert v\right\Vert \cdot\left\Vert w_{1}\right\Vert }
            \leqslant\delta\frac{1}{\left\Vert w_{1}\right\Vert }
            .
    \end{align}
    Since $\| g\|=\alpha_{1}$, we get that $\left\Vert gw_{1}\right\Vert \leqslant\left|\alpha_{1}\right|\left\Vert w_{1}\right\Vert $.
    In particular, 
    \[
        \frac{\left\Vert gv\wedge e_{1}\right\Vert }{\left\Vert gv\right\Vert \cdot\left\Vert w_{1}\right\Vert }
        \leqslant
        \left|\alpha_{1}\right|\frac{\left\Vert gv\wedge e_{1}\right\Vert }{\left\Vert gv\right\Vert \cdot\left\Vert gw_{1}\right\Vert }
        \leqslant
        \delta\frac{1}{\left\Vert w_{1}\right\Vert }.
    \]
    Cancelling out $\|w_1\|$ and multiplying by $\|gv\|$, we get
    \begin{align}\label{Eq:ineq2}
        \left\Vert gv\wedge e_{1}\right\Vert\leqslant\delta\left\Vert gv\right\Vert.
    \end{align}
    In the Archimedean case, 
    $
    \left\Vert gv\wedge e_{1}\right\Vert ^{2}=\left\Vert gv\right\Vert ^{2}-\left|\left\langle gv,e_{1}\right\rangle \right|^{2}
    $
    so
    \[
    \left|\left\langle gv,e_{1}\right\rangle \right|^{2}\geqslant(1-\delta^{2})\left\Vert gv\right\Vert ^{2}.
    \]
    Since $g$ is diagonal, $\left\langle gv,e_{1}\right\rangle =\alpha_{1}v_{1}$, hence
    \[
    \left|\alpha_{1}v_{1}\right|\geqslant\sqrt{1-\delta^{2}}\left\Vert gv\right\Vert .
    \]
    The non-Archimedean case is simpler: since $gv=\sum_{i=1}^{\dw}\alpha_{i}v_{i}e_{i}$,
    we get that $\left\Vert gv\right\Vert =\max_{i=1,\dots,\dw}\left|\alpha_{i}v_{i}\right|$. Since $gv\wedge e_1 =\sum_{i=2}^\dw \alpha_iv_i(e_i\wedge e_1)$, we get $\left\Vert gv\wedge e_{1}\right\Vert =\max_{i=2,\dots,\dw}\left|\alpha_{i}v_{i}\right|$. Therefore, Equation \ref{Eq:ineq2} implies (assuming $\delta<1$) that $\left|\alpha_{1}v_{1}\right|=\left\Vert gv\right\Vert $.
    
    Let us show that $\left\Vert gv\wedge e_{2}\right\Vert \geqslant\left|\alpha_{1}v_{1}\right|$.
    In the Archimedean case, this is because
    \[
    \left\Vert gv\wedge e_{2}\right\Vert ^{2}=\left\Vert gv\right\Vert ^{2}-\left|\left\langle gv\right\rangle ,e_{2}\right|^{2}=\sum_{i\neq2}\left|\alpha_{i}v_{i}\right|^{2}\geqslant\left|\alpha_{1}v_{1}\right|^2.
    \]
    In the non-Archimedean case, we again have $\left\Vert gv\wedge e_{2}\right\Vert =\max_{i\neq2}\left|\alpha_{i}v_{i}\right|$,
    so this is clear.
    
    We are almost done. We now do the same thing we did above with $e_{1}$,
    only with $e_{2}$. Equation \ref{Eq:ineq1} becomes:
    \[
    \left|\alpha_{2}\right|\frac{\left\Vert gv\wedge e_{2}\right\Vert }{\left\Vert gv\right\Vert }\leqslant\delta\frac{\left\Vert gw_{2}\right\Vert }{\left\Vert w_{2}\right\Vert }.
    \]
    This time,  $\left|\alpha_{2}\right|$ does not cancel out, because
    we only know that $\frac{\left\Vert gw_{2}\right\Vert }{\left\Vert w_{2}\right\Vert }\leqslant\left|\alpha_{1}\right|$.
    So we only get
    \[
    \left|\alpha_{2}\right|\frac{\left\Vert gv\wedge e_{2}\right\Vert }{\left\Vert gv\right\Vert }\leqslant\delta\left|\alpha_{1}\right|.
    \]
    Substituting $\left\Vert gv\wedge e_{2}\right\Vert \geqslant\left|\alpha_{1}v_{1}\right|\geqslant\sqrt{1-\delta^{2}}\left\Vert gv\right\Vert $
    in the Archimedean case, we get
    \[
    \left|\frac{\alpha_{1}}{\alpha_{2}}\right|\geqslant\frac{\sqrt{1-\delta^{2}}}{\delta}.
    \]
    Substituting $\left\Vert gv\wedge e_{2}\right\Vert \geqslant\left|\alpha_{1}v_{1}\right|\geqslant\left\Vert gv\right\Vert $
    in the non-Archimedean case, we get
    \[
    \left|\frac{\alpha_{1}}{\alpha_{2}}\right|\geqslant\frac{1}{\delta}.
    \]
    In either case, we can make $\left|\frac{\alpha_{1}}{\alpha_{2}}\right|$ arbitrarily
    large by choosing $\delta$ small enough.
\end{proof}

\subsection{Separation}
We fix again a vector space $W$ over a local field $\mbm k$, and a basis $e_1,\dots,e_\dw$.
\begin{defn}
A subset $F\subseteq\mathrm{GL}(W)$ is \emph{$(n,r)$-separating
}(for $n\in\mathbb{N}$, $r>0$) if, for any $n$ hyperplanes $H_{1},\dots,H_{n}$ in $\mb P(W)$ and $n$ points $v_{1},\dots,v_{n}\in \mb P(W)$, there is $f\in F$ such that 
\begin{align*}
d(fv_{i},H_{j})&>r~\text{and}~\\
d(f^{-1}v_{i},H_{j})&>r
\end{align*}
for all $1\le i,j\le n$.
\end{defn}

\begin{prop}\label{prop:separation}
Let $L$ be a Zariski connected algebraic subgroup of $\mr{SL}(W)$ acting irreducibly on $W$, and let $\Delta\leqslant L$ be a Zariski dense subgroup. For every $n\in\mathbb{N}$ there is some $r>0$ for which
there exists a finite $(n,r)$-separating subset $F\subseteq\Delta$.
\end{prop}

\begin{proof}
Let $n\in\mb N$. Since $L$ is Zariski connected and acts irreducibly, for every $n$
hyperplanes $H_{1},\dots,H_{n}$ in $\mathbb{P}\coloneqq\mb P(W)$ and $n$ points
$v_{1},\dots,v_{n}\in\mathbb{P}$, there is $\gamma\in\Delta$ such
that both $\gamma v_{i}\notin H_{j}$ and $\gamma^{-1} v_{i}\notin H_{j}$ for every $i,j$. To see this, set
\begin{align*}
V_{i,j}^+&=\left\{ h\in L\middle|hv_{i}\in H_{j}\right\} ,\\
V_{i,j}^-&=\left\{ h\in L\middle|h^{-1}v_{i}\in H_{j}\right\} ;
\end{align*}
these are proper algebraic subvarieties of $L$. Set $V_{i,j}=V_{i,j}^+\cup V_{i,j}^-$ for all $i,j$. Since $L$ is Zariski connected, it is irreducible
(as an algebraic variety), and hence
\[
\bigcup_{1\leqslant i,j\leqslant n}V_{i,j}\neq L.
\]
Since $\bigcup_{1\leqslant i,j\leqslant n}V_{i,j}$ is closed and $\Delta$ is Zariski dense, there is some $\gamma\in\Delta\backslash\bigcup V_{i,j}$. Consider
\[
X\coloneqq\left\{ (H_{1},\dots,H_{n},v_{1},\dots,v_{n})\right\} =\left(\mathrm{Gr}_{d-1}(\mathbb{P})\right)^{n}\times\mathbb{P}^{n}
\]
with the Hausdorff topology induced by the topology of $k$. Note that $X$ is compact.
For every $\gamma\in\Delta$, 
\[
X_{\gamma}=\left\{ (H_{1},\dots,H_{n},v_{1},\dots,v_{n})\middle|\gamma v_{i},\g^{-1}v_i\notin H_{j}\ \forall i,j\right\} 
\]
is open in the Hausdorff topology, by continuity of $\gamma$. By
the above, for every $H_{1},\dots,H_{n}$ and $v_{1},\dots,v_{n}$
there is some $\gamma\in\Delta$ such that $(H_{1},\dots,H_{n},v_{1},\dots,v_{n})\in X_{\gamma}$.
Therefore,
$$
X=\bigcup_{\gamma\in\Delta}X_{\gamma}
$$
is an open cover of $X$. By compactness, there is a finite subset $F\subseteq\Delta$
such that 
$$
X=\bigcup_{\gamma\in F}X_{\gamma}.
$$
We are almost done: we will now show that $F$ is an $(n,r)$-separating subset for some $r>0$. For $\gamma\in F$, denote 
$$
 d_\gamma (H_{1},\dots,H_{n},v_{1},\dots,v_{n})=\min_{i,j}\left\{ d(\gamma v_{i},H_{j}),d(\gamma^{-1} v_{i},H_{j})\right\}.
$$
The function
$
\max_{\gamma\in F}d_\gamma
:X\to \mb R$ is continuous, and by the above it is positive everywhere. 
Since $X$ is compact, it attains a minimum. Denoting the minimum by $r$, we get that $F$ is $(n,r)$-separating.
\end{proof}

\subsection{Constructing the Ping-Pong Players}
We can now prove Theorem \ref{thm:dyn}.

Let $\G\le\mr{SL}_d(\mbm k)$ be an unbounded subgroup whose Zariski closure in $\mr{SL}_d(\mbm k)$ is Zariski connected and acts irreducibly on $\mbm k^d$.
By Proposition \ref{prop:separation},
for some $r>0$, there exists a finite, $(2,r)$-separating
subset $F\subseteq\Gamma$.

Since the projective space is compact, every projective map is $L$-Lipschitz for some $L>0$. Let $C$ be an upper bound for the Lipschitz constant of all the maps in $F\cup F^{-1}$. Pick $\delta\in (0,r)$
small enough so that the argument below works. Pick $\gamma\in\G$
such that $\gamma$ is $\delta$-Lipschitz outside $(H)_{\delta}$ for some hyperplane $H$. 
There is no reason to believe $\g^{-1}$ is contracting as well, but we do know that there is some open subset $\Omega\subseteq \mb P\coloneqq \mb P(\mbm k^d)$ such that $\g^{-1}\!\restriction_\Omega$ is $2$-Lipschitz (by Lemma \ref{Lip}). Pick $v\in \Omega$ 
 and $f\in F$ such that $f.(\g^{-1}.v)\notin (H)_r$ and $f^{-1}.(\g^{-1}.v)\notin (H)_r$. Set $\g_0=\g f\g^{-1}$. By the composition rule for Lipschitz constants, we get that, in a small neighbourhood of $v$, both $\g_0$ and $\g_0^{-1}$
 are $\delta'$-Lipschitz for $\delta'=2C\delta$.

By Proposition \ref{LipIsCont}, both $\g_0$ and $\g_0^{-1}$ are $\ep$-contracting for $\ep=\ep(\delta')$, so by choosing $\delta$ (and hence $\delta'$) small enough, we may suppose that $\ep$ is as small as we wish. 

Let $v_{0}^{+},v_0^-,H_{0}^{+},H_0^-$ be such that $\g_0$ is $\ep$-contracting with respect to $(v_0^+,H_0^+)$ and $\gamma_{0}^{-1}$ is $\ep$-contracting
with respect to $(v_{0}^{-},H_{0}^{-})$. 
That is,
\begin{align*}
\gamma_0\left(\mathbb{P}(W)\setminus\left(H_{0}^{+}\right)_{\varepsilon}\right) & \subseteq(v_{0}^{+})_{\varepsilon},\\
\gamma_0\left(\mathbb{P}(W)\setminus\left(H_{0}^{-}\right)_{\varepsilon}\right) & \subseteq(v_{0}^{-})_{\varepsilon}.
\end{align*}
Let $f_0\in F$ such that 
\begin{align}
d(f_0v_{0}^{+},H_{0}^{+}) & >r,\label{proximal}\\
d(f_0^{-1}v_{0}^{-},H_{0}^{-}) & >r.
\end{align}
Since $f_0^{-1}$ is $C$-Lipschitz, the second requirement implies 
\begin{align}\label{proximall}
d(v_{0}^{-},f_0H_{0}^{-})>\frac{r}{C}.
\end{align}
\begin{defn}
Let $\ep_0,r_0>0$. A map $g:\mathbb{P}\to\mathbb{P}$ is \emph{$(\ep_0,r_0)$-proximal
with respect to $(v,H)$} (where $v\in\mathbb{P}$ is a point and
$H\subseteq\mathbb{P}$ a hyperplane) if it is $\ep_0$-contracting with respect to $(v,H)$, and 
\[
d(v,H)>r_0.
\]
\end{defn}

Let $x\notin(H_{0}^{+})_{\varepsilon}$. Since $f_{0}$ is $C$-Lipschitz,
we have
\[
d(f_{0}\gamma_{0}x,f_{0}v_{0}^{+})\leqslant C\cdot d(\gamma_{0}x,v_{0}^{+})\leqslant C\varepsilon.
\]
In other words, $\g_1\coloneqq f_{0}\gamma_{0}$ is $(C\varepsilon,r)$-proximal
with respect to $(f_{0}v_{0}^{+},H_{0}^{+})$.

Now, consider $\gamma_{1}^{-1}=\g_0^{-1}f_0^{-1}$, and let us show it is $(C\varepsilon,\frac{r}{C})$-proximal
with respect to $(v_{0}^{-},f_0H_{0}^{-})$. Let $x\notin(f_{0}H_{0}^{-})_{C\varepsilon}$. Then $f_{0}^{-1}x\notin(H_{0}^{-})_{\varepsilon}$
(since $f_{0}$ is $C$-Lipschitz). Thus,
\[
d(\gamma_{0}^{-1}f_{0}^{-1}x,v_{0}^{-})<\varepsilon\leqslant C\varepsilon
\]
since $\gamma_{0}^{-1}$ is $\varepsilon$-contracting with respect
to $(v_{0}^{-},H_{0}^{-})$. By Equation \eqref{proximall}, $d(v_{0}^{-},f_{0}H_{0}^{-})>\frac{r}{C}$. 

Set $v_1^+=f_0v_0^+$, $H_1^+=H_0^+$, so that $\g_1$ is $(C\ep,r)$-proximal with respect to $(v_1^+,H_1^+)$. Set $v_1^-=v_0^-$, $H_1^-=f_0H_0^-$, so that $\g_1^{-1}$ is $(C\ep,\frac{r}{C})$-proximal with respect to $(v_1^-,H_1^-)$.

Now, choose $f_{1}\in F$ such that 
\begin{align}\label{proxiend}
d(f_{1}v_{1}^{+},H_{1}^{+}\cup H_{1}^{-}) & >r,\\
d(f_{1}v_{1}^{-},H_{1}^{+}\cup H_{1}^{-}) & >r.
\end{align}
Let us show $\gamma_{2}\coloneqq f_{1}\gamma_{1}f_{1}^{-1}$ is $(C^{2}\varepsilon,\frac{r}{C})$-proximal
with respect to $(f_{1}v_{1}^{+},f_{1}H_{1}^{+})$. Let $x\notin(f_{1}H_{1}^{+})_{C^{2}\varepsilon}$.
Then $f_{1}^{-1}x\notin(H_{1}^{+})_{C\varepsilon}$ (because $f_{1}$
is $C$-Lipschitz). Thus,
\[
d(\gamma_{1}f_{1}^{-1}x,v_{1}^{+})<C\varepsilon,
\]
(since $\gamma_{1}$ is $C\varepsilon$-contracting with respect to
$(v_{1}^{+},H_{1}^{+})$), so
\[
d(f_{1}\gamma_{1}f_{1}^{-1}x,f_{1}v_{1}^{+})<C^{2}\varepsilon
\]
(since $f_{1}$ is $C$-Lipschitz). We also have $d(f_{1}v_{1}^{+},f_{1}H_{1}^{+})\geqslant\frac{1}{C}d(v_{1}^{+},H_{1}^{+})>\frac{r}{C}$
(because $f_{1}^{-1}$ is $C$-Lipschitz), so $\gamma_{2}=f_{1}\gamma_{1}f_{1}^{-1}$
is $(C^{2}\varepsilon,\frac{r}{C})$-proximal with respect to $(f_{1}v_{1}^{+},f_{1}H_{1}^{+})$.

The exact same calculations (with $\cdot^-$ in stead of $\cdot^+$) shows that $\gamma_{2}^{-1}=f_{1}\gamma_{1}^{-1}f_{1}^{-1}$ is $(C^{2}\varepsilon,\frac{r}{C^{2}})$-proximal
with respect to $(f_{1}v_{1}^{-},f_{1}H_{1}^{-})$. 

Now, set $v_{2}^{+}=f_{1}v_{1}^{+}$ and $H_{2}^{+}=f_{1}H_{1}^{+}$
so that $\gamma_{2}$ is $(C^{2}\varepsilon,\frac{r}{C})$-proximal
with respect to $(v_{2}^{+},H_{2}^{+})$, and set $v_{2}^{-}=f_{1}v_{1}^{-}$
and $H_{2}^{-}=f_{1}H_{1}^{-}$ so that $\gamma_{2}^{-1}$ is $(C^{2}\varepsilon,\frac{r}{C^{2}})$-proximal
with respect to $(v_{2}^{-},H_{2}^{-})$. By Equations \eqref{proxiend}, we get that
\begin{align*}
d(\left\{ v_{2}^{+},v_{2}^{-}\right\} ,H_{1}^{+}\cup H_{1}^{-})  >r.
\end{align*}
Similarly, we have
\[
d(\left\{ v_{1}^{+},v_{1}^{-}\right\} ,H_{2}^{+}\cup H_{2}^{-})  >\frac{r}{C}.
\]
Therefore, if we took $\varepsilon$ small enough so that $C^{4}\ep<r$,
we'd get that $\gamma_{1},\gamma_{2}$ play ping pong, and hence generate
a free subgroup.

This completes the proof of Theorem \ref{thm:dyn} and the Tits alternative,
Theorem \ref{thm:Tits}.
\cleardoublepage
\phantomsection
\addcontentsline{toc}{chapter}{Bibliography}
\bibliographystyle{abbrv}
\bibliography{mybibliography}

\end{document}